\documentclass[11pt]{article}
  \usepackage{amsmath,amsfonts,amssymb,amsthm}
  \usepackage{graphics}
  \usepackage{epsfig}
  \usepackage{xcolor}
  \usepackage{bm}
  \usepackage{hyperref}
  \usepackage{pgfplots}
  \usepackage{tikz}
  \usepackage{vmargin}
  
 \usepackage{stackrel}
  \usetikzlibrary{calc,arrows}
  
  \hypersetup{
  colorlinks,
  citecolor=blue,
  linkcolor=blue,
  urlcolor=black
  }
  \hypersetup{pageanchor=false}
   
  \usepackage[nottoc]{tocbibind}

  \newcommand{\ignore}[1]{}

  \newtheorem{proposition}{Proposition}
  \newtheorem{theorem}{Theorem}
  
  \newtheorem{lemma}{Lemma}
  
  \newcommand{\R}{\mathbb{R}}
  \newcommand{\Z}{\mathbb{Z}}
  
  \newcommand{\F}{\mathcal{F}}

  \newcommand{\pp}{\overset{\approx}{\!\phi_1}}
  \newcommand{\ppa}{\overset{\approx}{\!\phi^\alpha_1}}
  \newcommand{\dppa}{\delta\overset{\approx}{\phi^\alpha_1}}
   \newcommand{\BF}{\textnormal{\tiny \textsc{bf}}}
   \newcommand{\trho}{\widetilde{\rho}}
   \newcommand{\ttrho}{\overset{\approx}{\rho}}
   
  \def\Ra{{\rm Ra}}
  \def\Pr{{\rm Pr}}
  \def\Nu{{\rm Nu}}
  
 \newtheorem*{theorem*}{Theorem}
 
\author{Camilla Nobili and Felix Otto}

\title{
Limitations of the background field method applied
to Rayleigh-B\'enard convection
}

\begin{document}

\maketitle

  \begin{abstract}

  We consider Rayleigh-B\'enard convection as modeled by the Boussinesq equations,
  in case of infinite Prandtl number and with no-slip boundary condition.
  There is a broad interest in bounds of the upwards heat flux,
  as given by the Nusselt number ${\rm Nu}$, in terms of the forcing via the imposed temperature difference,
  as given by the Rayleigh number in the turbulent regime ${\rm Ra}\gg 1$. In several works, the background field method
  applied to the temperature field
  has been used to provide upper bounds on ${\rm Nu}$ in terms of ${\rm Ra}$. 
  In these applications, the background field method
  comes in form of a variational problem where one optimizes a stratified temperature profile subject to a certain
  stability condition; the method is believed to capture marginal stability of the boundary layer.
  The best available upper bound via this method is ${\rm Nu}$
  $\lesssim {\rm Ra}^\frac{1}{3}(\ln {\rm Ra})^\frac{1}{15}$;
  it proceeds via the construction of a stable temperature background profile that increases logarithmically
  in the bulk. In this paper, we show that the
  background temperature field method cannot provide a tighter upper bound in terms of the power of the logarithm.
  However, by another method one does obtain the tighter upper bound 
  ${\rm Nu}\lesssim {\rm Ra}^\frac{1}{3}(\ln\ln {\rm Ra})^\frac{1}{3}$,
  so that the result of this paper implies that the background temperature field method is unphysical
  in the sense that it cannot provide the optimal bound.

  \smallskip


  \end{abstract}

  \newpage
  \tableofcontents

\newpage
\section{Introduction}

In a $d-$dimensional container of height normalized to unity we consider Rayleigh-B\'enard convection  
as modeled by the Boussinesq equations, which we consider in their infinite-Prandtl-number limit:
 \begin{subequations}
  \label{RBC}
  \begin{alignat}{2}
   \partial_t T+u\cdot \nabla T&= \Delta T  \qquad & {\rm for } \quad 0<z<1\,,\\
   -\Delta u+\nabla p&=\Ra T e_z \qquad & {\rm for } \quad 0<z<1 \,,\\
   \nabla\cdot u&= 0 \qquad & {\rm for } \quad 0<z<1\,,\\
   u&=0  \qquad & {\rm for } \quad z\in\{0,1\}\,,\label{new-4}\\
   T&=1 \qquad & {\rm for } \quad z=0\,,\label{new-5}\\
   T&=0 \qquad & {\rm for } \quad z=1 \label{new-6}\,.
  \end{alignat}
  \end{subequations}
   Here $u\in \R^d$ denotes the fluid velocity, $T\in \R$ its temperature and $p\in \R$ its pressure.
  We denote with $z$ the vertical component of the $d-$dimensional
   position vector $x=(y,z)$ and  with $e_z$ the upward unit normal in the vertical direction.
  As a convenient proxy of the side-wall effect, the functions $u, T$ and $p$, which depend on the spatial variable $x$ and the 
  time variable $t$, are supposed to be periodic in the $(d-1)-$horizontal directions $y$ with period $L$, where 
  $L$ is the horizontal period. In our treatment, the dimension $d$ is arbitrary and we think of $L$ as being large.
  The first equation encodes the diffusion of the temperature, driven by the Dirichlet boundary conditions (\ref{new-5})\&(\ref{new-6}), 
  and its advection by the fluid velocity. The second equation, the Stokes equation, encodes the fact that the fluid parcels move as a reaction to the
  buoyancy force $\Ra Te_z$ (hotter parcels expand and thus experience an upwards force under gravity) and are slowed down by viscosity
  ($-\Delta u$) in conjunction with the no-slip boundary condition (\ref{new-4}). The last equation expresses the incompressibility of the fluid
  and is balanced by the pressure term $\nabla p$ acting as a Lagrangian multiplier in the previous equation.
  The parameter appearing in the Stokes equation, the Rayleigh number $\Ra$,
  expresses the relative strength of the buoyancy force and is given by
  \begin{equation}\label{Rayleigh-number}
  \Ra=\frac{\alpha g (T_{b}-T_{t}) h^3}{\nu\kappa}\,,
  \end{equation}
  where $\alpha$ is the thermal expansion coefficient, $\nu$ the kinematic viscosity, $\kappa$ the 
  thermal conductivity, $(T_{b}-T_{t})$ the temperature gap between the bottom and the top plate and 
  $h$ the height of the container before the non-dimensionalization.
  
  In (\ref{RBC}), the inertia of the fluid has been neglected, which amounts to sending the Prandtl number $\Pr=\frac \nu\kappa$ to infinity.
  Therefore $\Ra$, next to $L$, is the only non-dimensional parameter.
  The linear stability analysis identifies a critical value $\Ra_c$, the Rayleigh number at which 
   the solution of (\ref{RBC}) bifurcates from the linear conduction profile $T=1-z,\; u=0, \; p=z-\frac{z^2}{2}$, see for instance \cite{Gelting}.
  When $\Ra>\Ra_c$, the buoyancy forces trigger the formation of convection rolls. Eventually, when $\Ra\gg \Ra_c$, these convection rolls break down.
  This regime features boundary layers at the top and bottom plates, with a high vertical temperature gradient, from which small fluid parcels
  of different temperature than the ambient fluid detach and deform, the so called plumes.
  In this paper we are interested in this turbulent regime of
  $$\Ra\gg 1\,,$$
  and in the experimentally observed enhancement of the heat transport over the pure conduction state.
  An appropriate measure to quantify the vertical heat flux is the Nusselt number, 
  \begin{equation*}
   \Nu=\int_0^1\langle(uT-\nabla T)\cdot e_z\rangle\,dz\,,
  \end{equation*}
  which represents the average heat flux $uT-\nabla T$ passing through an area element. 
  The bracket $\langle\cdot\rangle$ denotes the time and horizontal space average 
  \begin{equation}\label{bra}
  \langle f \rangle:=\limsup_{t_0\uparrow\infty}\frac{1}{t_0}\int_0^{t_0}\frac{1}{L^{d-1}}\int_{[0,L)^{d-1}}f(t,y)dydt\,,
   \end{equation}
   and might be thought as a statistical average. 
   
   In the fifties Malkus \cite{Ma}, considering fluids with very high viscosity,
  performed experiments in which he noticed sharp transitions in the slope of the $\Nu-\Ra$ relation and suggested 
  the scaling $\Nu\sim \Ra^{\frac 13}$ for very high $\Ra$ numbers based on the \textit{marginal stability argument}, which
  we reproduce now.
  Since the main temperature drop happens near the boundary, we can assume that in a bottom boundary layer of thickness $\delta$ (to be determined)
   the temperature drops from $1$ to its average $\frac 12$. Thanks to the average $\langle\cdot\rangle$, we may extract 
   the Nusselt number from the boundary layer, where by the no-slip boundary condition we have $(uT-\nabla T)e_z\approx \partial_z T$
   so that $\Nu\sim \frac{1}{\delta}$. The boundary layer is assimilated to the pure conduction state of height $\delta$. Marginal 
   stability refers to the assumption that this state is borderline stable, meaning that its Rayleigh number is critical, which in 
   view of (\ref{Rayleigh-number}) means 
  $$\Ra_c=\frac{g\alpha (T_b-T_t)(\delta h)^3}{\nu \kappa}\,,$$
  from which, because of $\Ra_c\sim 1$, we infer $\delta\sim \Ra^{-\frac 13}\,.$
  Inserting this in the scaling of Nusselt number above one finds
  $$\Nu\sim \Ra^{\frac 13}\,.$$
  The same conclusion can be achieved by rescaling equation (\ref{RBC}) according to
  \begin{equation}\label{rescaling1}
  x=\Ra^{-\frac{1}{3}}\hat x, \; t=\Ra^{-\frac{2}{3}}\hat t,\; u=\Ra^{\frac{1}{3}}\hat u, \; p=\Ra^{\frac{2}{3}}\hat p \; \mbox{ and thus } \;\Nu=\Ra^{\frac{1}{3}}\widehat{\Nu}\,.
  \end{equation}
  In this way we end up with the parameter-free system  
  \begin{equation*}
   \begin{array}{rclc}
   \partial_{\hat t}  T+\hat u\cdot \hat{\nabla}  T &=& \hat{\Delta}T  \,,\\
    -\hat\Delta \hat u+\hat\nabla \hat p&=& T e_{\hat z} \,,\\
    \hat{\nabla}\cdot \hat u&=& 0 \,,\\
    \end{array}
  \end{equation*}
  which naturally lives in the half space.
  Since for the latter system, it is natural to expect that the heat flux is universal, i.\ e.\ $\widehat{\Nu}\sim 1$, we also obtain
  $\Nu\sim \Ra^{\frac{1}{3}}.$

  The scaling $\Nu\sim \Ra^{\frac 13}$ has been been confirmed by experiments at (relatively) 
  high Prandtl numbers (cf.\ \cite{Gro-Lo}, p.$30$, for a list of experimental 
  results).  
  Rigorous analyses have produced upper bounds that capture this scaling up to logarithms which we report now.
  In the sixties Howard \cite{Howard} obtained an upper bound that scales like $\Ra^{\frac 12}$, optimizing over a field of test functions satisfying physical constraints
  coming from the Navier-Stokes equation (cf.\ \cite{Howard}, Sec.$3$), while neglecting the incompressibility constraint.
   Later Busse \cite{Busse} developed the theoretical tool of \textit{multiple boundary layer solution} (multi-$\alpha$ solution) in order to solve 
   Howard's variational problem when the incompressibility constraint is taken into account (cf.\ \cite{Busse}, Sec.$2$).
  The multi-$\alpha$ solution theory inspired Chan \cite{Chan} in the seventies. He elegantly applied it, deriving an upper bound on the Nusselt number that scales like $\Ra^{\frac 13}$ when additional conditions in the asymptotic analysis are assumed.
   In the nineties, Constantin and Doering, inspired by the works of Malkus, Howard and Busse, introduced 
  the \textit{background field method} in order to bound the average dissipation rate in plane Couette flow.
  This method was already implicitly used by Hopf \cite{Hopf} in the forties in the construction of solutions to the Navier-Stokes equations (in the sense of Leray) with inhomogeneous
  boundary data. Later, Constantin and Doering applied it to the Rayleigh-B\'{e}nard convection in order to derive rigorous upper bounds for the Nusselt number $\Nu$.
  Although Howard's problem and the background field method constitute dual variational problems \cite{Ker}, the second method has the advantage to use simple test functions
  and functional estimates.  Indeed this method turned out to be very fruitful: It has been extensively used and it has produced meaningful bounds in the theory of turbulence. 
  In the context of Rayleigh-B\'enard convection, this method
  consists of decomposing the temperature field $T$ into a steady background temperature field profile $\tau=\tau(z)$ with driving boundary conditions, $\tau=1$ for $z=0$ and $\tau=0$ for $z=1$, and into temperature fluctuations $\theta$.
  As we will see in detail in the next subsection, the advantage of the background field method is to transform the problem of finding upper bounds for the Nusselt number into a purely variational problem: Find profile (test) functions $\tau$ which satisfies a certain stability condition and then select the one with minimal Dirichlet energy. The solution of this variational method produces a new number $\Nu_{\BF}$ such that 
  \begin{equation}\label{sol-bfm}\Nu\leq\Nu_{\BF}.\end{equation} 
  Experiments suggest to try a profile $\tau$ that displays a drop by $\frac 12$ in a boundary layer and it is constant in the bulk. Such a 
  profile satisfies the stability condition only if the boundary layer size $\delta$ is chosen artificially small and gives only suboptimal bounds (see \cite{DC2001}). Replacing the constant bulk by a linearly increasing profile (at the expense of making the drop in the boundary layers deeper)
  does not improve the situation.
  The idea that the ``bad'' boundary layers can be more efficiently compensated by a 
  profile that increases fast near the boundary and slowly (almost constant) away from them, brought 
  Doering, Reznikoff and the second author in $2006$ \cite{DOR} to investigate the stability of a background profile that grows logarithmically in the bulk. 
  This Ansatz indeed proved to be successful, yielding the bound $\Nu_{\BF}\lesssim\Ra^{\frac 13}(\ln \Ra)^{\frac{1}{3}}$ and therefore reproducing the scaling proposed by Malkus up to a logarithmic correction.
   Seis and the second author \cite{OS2011} in $2011$ improved the 
  last bound by reducing the logarithmic correction
  \begin{equation}\label{OS}
  \Nu_{\BF}\lesssim \Ra^{\frac 13}(\ln \Ra)^{\frac{1}{15}}.
  \end{equation}
  They used the same logarithmic construction as in \cite{DOR} with a (logarithmically) larger boundary layer thickness, which they could 
  afford using an additional estimate on the vertical velocity component $w=u\cdot e_z$ in terms of $\theta$.
 
   In the context of the Rayleigh-B\'enard convection, the background field method has also been used to study 
  the case of free-slip boundary condition for the velocity field \cite{IKP}, of an
  imposed heat flux at the boundary \cite{Otero-witt} and in the bulk \cite{Doer-Whit-IHC}, of mixed thermal boundary conditions \cite{whit-witt},
  and of rough boundaries \cite{GD}. This method has been fruitfully applied to a variety of other problems in fluid mechanics, namely plane Couette flow \cite{Shear}, pipe flow, and arbitrary Prandtl number convection.
  Nicolaenko, Scheuer and Temam \cite{nicolaenko} applied the background field method to derive an upper bound for the long-time limit of the $L^2$-norm of the solution of the Kuramoto-Sivashinsky equation. 
   
   In this paper, we address the question of the optimality of the background field method in two ways:
   \begin{itemize}
    \item What is the optimal bound (in terms of the two scaling exponents in $\Ra^{\mu}(\ln \Ra)^{\nu}$) in the 
    background field method?
    The answer given in our main result is that the construction in \cite{DOR} and \cite{OS2011} leading to (\ref{OS}) is indeed optimal.
     \item Does the background field method catch the optimal bound on $\Nu$, in other words, is (\ref{sol-bfm}) optimal (at least in terms of
     both scaling exponents)? The answer given by our main result in conjunction with the bound in (\ref{best-upper-bound}) obtained in \cite{OS2011} is no.
     \end{itemize}
     The main result of this paper is stated in the following 
  \begin{theorem*} For $\Ra\gg 1$ we have 
   \begin{equation}\label{lower-bound-preview}
  \Nu_{\BF}\gtrsim \Ra^{\frac 13}(\ln\Ra)^{\frac{1}{15}}\,.
  \end{equation}
  \end{theorem*}
  \noindent
  We refer to Theorem \ref{Th1} for the full formulation and the explanation of the notation $\gtrsim$ and $\gg$.
 This lower bound on $\Nu_{\BF}$, together with the upper bound (\ref{OS}), implies
 $$\Nu_{\BF}\sim\Ra^{\frac 13}(\ln \Ra)^{\frac{1}{15}}\,,$$ and in particular shows that the background field method
 cannot produce  any smaller logarithmic correction than $(\ln \Ra)^{\frac{1}{15}}$.
 However, a combination of the background field with another method improves the logarithmic correction in (\ref{OS}), as we shall explain now.
 This other method, which we refer to as maximal principle method, was also introduced by Constantin and Doering \cite{CD99}.
  Indeed it is easy to verify that the temperature equation satisfies the maximum principle, leading to
 \begin{equation}\label{MP}
  0\leq T\leq  1\,,
 \end{equation}
 possibly neglecting an initial layer.
 This $L^{\infty}$ bound together with a maximal regularity estimate for Stokes equation in $L^{\infty}$ yields the bound 
 $$\Nu\lesssim \Ra^{\frac 13}(\ln\Ra)^{\frac 23}\,,$$
 where, this time, the logarithm is an expression of the failure of the $L^{\infty}-$norm to be a Calder\'on-Zygmund norm.
 Recently, Seis and the second author \cite{OS2011} combined the maximum principle method with the background field method developed in \cite{DOR},
  obtaining 
 \begin{equation}\label{best-upper-bound}
 \Nu\lesssim \Ra^{\frac 13}(\ln \ln \Ra)^{\frac{1}{3}}\,,
 \end{equation}
 which, to our knowledge, is the best rigorous upper bound. 
 We observe that the combination of all the previous results yields
 $$\Nu\leq\Nu_{\BF}\stackrel{(\ref{best-upper-bound})}{\lesssim} \Ra^{\frac 13}(\ln \ln \Ra)^{\frac{1}{3}}\ll \Ra^{\frac 13}(\ln \Ra)^{\frac{1}{15}}\stackrel{(\ref{lower-bound-preview})}{\lesssim} \Nu_{\BF}\,.$$
 So indeed, the background field method is not able to capture the behavior of the Nusselt number even in terms of the two scaling exponents.
 Therefore, the optimal background temperature profile cannot carry much of a physical meaning.
  
  In $2005$ Plasting and Ierley \cite{PI-p1} considerer piecewise linear 
  profiles with $\tau'\geq 0$ in the bulk and solved the variational problem numerically, finding $\Nu\sim \Ra^{\frac{7}{20}}$. In $2006$ inspired by Chan's multi$-\alpha$ solution treatment, Ierley, Kerswell and Plasting \cite{IKP-p2}, with help of a mixture of numerical and analytical methods, improved the previous result finding 
  $$\Nu\sim c_1\Ra^{\mu}(\ln \Ra)^{\nu}\,,$$ 
  where $\mu=0.33175$ and $\nu=0.0325$. Clearly, since $0.0325<0.0\bar 6=\frac{1}{15}$,
   our result (although slightly underestimated) has been anticipated ten years ago.
   
   Incidentally, the background field method suffers a similar fate in the context of the
  Kuramoto-Sivashinsky equation: Bronski and Gambill \cite{bronski} identified the optimal scaling
  of the upper bound that can be obtained by this method, and soon later, Giacomelli and
  the second author \cite{Gia-Otto} showed that a tighter upper bound can be obtained by an alternative
  method, which has subsequently been further improved in \cite{Otto-KSE} and \cite{GJO}. 

  \subsection{Temperature background field method and main result}\label{HFV-sec}
  We start by the rescaling (\ref{rescaling1}) suggested by Malkus' marginal stability argument, that is 
  \begin{equation*}
  x\rightsquigarrow\Ra^{\frac{1}{3}} x, \; t \rightsquigarrow \Ra^{\frac{2}{3}} t,\; u\rightsquigarrow\Ra^{-\frac{1}{3}} u \mbox{ and } p\rightsquigarrow\Ra^{-\frac{2}{3}} p \;,
  \end{equation*}
  and setting $H:=\Ra^{\frac 13}$ we rewrite (\ref{RBC}) as 
  \begin{subequations}
  \label{RBC2}
  \begin{alignat}{2}   
   \partial_t T+u\cdot \nabla T &=\Delta T  \qquad & {\rm for } \quad 0<z<H\,,\label{Ad-Diff}\\
   -\Delta u+\nabla p&= T e_z \qquad & {\rm for } \quad 0<z<H \,,\label{Stokes-eq}\\
   \nabla\cdot u&= 0 \qquad & {\rm for } \quad 0<z<H\,,\label{Incompress}\\
   u&=0  \qquad & {\rm for } \quad z\in\{0,H\}\,,\label{bc-for-u}\\
   T&=1 \qquad & {\rm for } \quad z=0\, ,\label{bc1-for-T}\\
   T&=0 \qquad & {\rm for } \quad z=H.\label{bc2-for-T}
  \end{alignat}
  \end{subequations}
   Notice that in this non-dimensionalization of the equation, the only parameter appearing is the height $H$ of the container
  and Malkus' scaling $\Nu\sim\Ra^{\frac 13}$ corresponds to $\Nu\sim 1$. 
  
  We recall from the previous section that the Nusselt number is defined as 
  \begin{equation}\label{Nu-1}
   \Nu=\frac{1}{H}\int_0^H\langle(uT-\nabla T)\cdot e_z\rangle\,dz\,,
  \end{equation}
  and now derive some useful representations starting from the equation for the temperature in (\ref{Ad-Diff}):
  Applying $\langle\cdot\rangle$ to the equation (\ref{Ad-Diff}) and qualitatively using the bound (\ref{MP}) on $T$ given by the maximum principle 
   it is easy to show that the upward heat flux is constant in the vertical direction, 
  \begin{equation}\label{heat-flux-property}
  \Nu=\langle Tw -\partial_z T\rangle \quad \mbox{ for } z\in (0,H)\,.
  \end{equation}
  Testing the  equation with $T$, appealing to incompressibility (\ref{Incompress}) and using (\ref{heat-flux-property}) for $z=0$, 
  we obtain (see \cite{DC2001}) the alternative representation
  \begin{equation}\label{Nusselt-nabla-T}
  \Nu=\int_0^H\langle |\nabla T|^2\rangle dz\,.
  \end{equation}
 
  For the convenience of the reader we sketch the derivation of the background field method, see \cite{CD96} for more details.
  The background field method consists of decomposing the temperature field $T$ into a steady background temperature profile $\tau$ which depends only on the vertical variable $z$ and satisfies the inhomogeneous (driving) boundary conditions, $\tau=1$ for $z=0$ and $\tau=0$ for $z=H$, and into temperature fluctuations $\theta$, with homogeneous boundary conditions $\theta=0$ for $z\in\{0,H\}$. 
  Inserting this decomposition
  \begin{equation}\label{decomposition}
  T=\tau+\theta\,
  \end{equation}
   in the equation for the temperature (\ref{Ad-Diff}) we find that the fluctuations $\theta$ evolve
  according to $$\partial_t\theta+u\cdot\nabla \theta-\Delta\theta=\frac{d^2\tau}{dz^2}-w\frac{d\tau}{dz}\,.$$
  From the incompressibility condition (\ref{Incompress}) we obtain by testing with $\theta$
  $$\int_0^H\langle|\nabla \theta|^2\rangle\, dz=-\int_0^H\frac{d\tau}{dz}\langle\partial_z\theta\rangle\, dz-\int_0^H\frac{d\tau}{dz}\langle\theta w\rangle\, dz\,.$$
  Together with (\ref{Nusselt-nabla-T}) this yields the final representation 
  \begin{equation}\label{Nusselt-new}
  \Nu=\int_0^H\left(\frac{d\tau}{dz}\right)^2dz-\int_0^H\left\langle2\frac{d\tau}{dz}\theta w+|\nabla\theta|^2\right\rangle dz\,.
  \end{equation}
  Applying the divergence to the Stokes equation (\ref{Stokes-eq}) we find that the pressure satisfies $\Delta p= \partial_z T$. Inserting $\Delta p$
  into the equation $\Delta$ (\ref{Stokes-eq})$\cdot e_z$, we 
  find the direct relationship between $\theta$ and the vertical velocity component $w:=u\cdot e_z$: 
  \begin{equation}\label{w<-->theta}
  \begin{array}{rclc}
  \Delta^2 w&=&-\Delta_{y} \theta   \qquad & {\rm for } \quad 0<z<H\,,\\
   w&=&\partial_z w =0  \qquad & {\rm for } \quad z\in\{0,H\}\,.\\
   \end{array}
 \end{equation}
  The representation (\ref{Nusselt-new}) shows: Any $\tau=\tau(z)$ that satisfies the driving boundary conditions and is stable in the sense that the 
  following quadratic form is non-negative  
  \begin{equation}\label{stab-cond}
   \int_0^H\left\langle2\frac{d\tau}{dz}\theta w+|\nabla\theta|^2 \right\rangle dz\geq 0\,, \\
  \end{equation}
  for every $\theta$ satisfying homogeneous boundary conditions (and $w$ defined through (\ref{w<-->theta})),
  yields an upper bound for the Nusselt number:
  \begin{equation*}\Nu\leq \int_0^H\left(\frac{d\tau}{dz}\right)^2 \; dz\,.\end{equation*}
  Note that in (\ref{stab-cond}), we may disregard the time variable, which is only a parameter in (\ref{w<-->theta}), so that
  $\langle\cdot\rangle$ in (\ref{stab-cond}) reduces to the horizontal average. 
  This motivates to define the Nusselt number associated to the background field method as
  \begin{equation}\label{tilde-Nu}
   \Nu_{\BF}:=\inf_{\substack{\tau:(0,H)\rightarrow \R, \\ \tau(0)=1, \tau(H)=0}}\left\{\int_0^H\left(\frac{d\tau}{dz}\right)^2dz|\;\tau \mbox{ satisfies }(\ref{stab-cond})\right\}\,,
  \end{equation}
  which in view of (\ref{Nusselt-new}) satisfies 
  $$\Nu\leq \Nu_{\BF}\,.$$
  
  Our objective is to derive an Ansatz-free lower bound for $\Nu_{\BF}$, trying to extract
  local information on $\tau$ from the completely non-local stability condition (\ref{stab-cond}). 
  The full formulation of the main result, already stated in the previous section, is contained in the following
  \begin{theorem}\label{Th1}
  Suppose that 
  \begin{equation}\label{IC}
  \tau(0)=1, \quad\tau(H)=0\,,
  \end{equation}
  and $\tau$ satisfies (\ref{stab-cond}) for all $(\theta,w)$ related by (\ref{w<-->theta}) in its Fourier transformed version (\ref{OSCinFourier-reduct}).
  Then for $H\gg 1$ $$\int_0^H\left(\frac{d\tau}{dz}\right)^2\,dz\gtrsim (\ln H)^{\frac{1}{15}}\,.$$ In particular
  for the Nusselt number associated to the background field method
  we have the lower bound
  \begin{equation}\label{lower-bound}
  \Nu_{\BF}\gtrsim (\ln H)^{\frac{1}{15}}\,.
  \end{equation}
  \end{theorem}
  \noindent
  Here and in the sequel, $\lesssim$ stands for $\leq C$ for some generic universal constant $C<\infty$.
  Likewise, $H\gg 1$ means that there exists a universal constant $C<\infty$ such that the statement holds for $H\geq C$.
  
  Besides implying the non-optimality of the background field method, this theorem offers some insights. 
  Indeed the proof is based on a characterization of profiles that satisfy the stability condition (\ref{stab-cond}) (see Section \ref{Lemmas}). This characterization is motivated by the analysis of a reduced form of the stability condition (\ref{RSC}) which indicates that long-wave length stability implies (approximate) logarithmic growth of $\tau$ in $z$, while short wave-length stability implies (approximate) monotonicity (see Proposition \ref{prop1} in the next section).

   It is convenient to introduce the slope $\xi:=\frac{d\tau}{dz}$ of the background temperature profile.
   With this convention the stability condition (\ref{stab-cond}) can be rewritten explicitly as follows 
 \begin{equation}\label{stab-cond-re}
 2\int_0^H \xi\langle w\theta\rangle\, dz+\int_0^H\langle |\nabla_{y}\theta|^2\rangle\, dz+\int_0^H \langle|\partial_z \theta|^2\rangle \, dz \geq 0,
 \end{equation}
  for all functions $\theta$ (and $w$ related to $\theta$ via the fourth-order boundary value problem (\ref{w<-->theta})) that  vanish at $z\in \{0,H\}$.
  A major advantage of the background field method is that it is amenable to (horizontal) Fourier transform:
  Indeed, denoting by $k\in \frac{2\pi}{L}\Z^{d-1}, k\neq 0$, the horizontal wavenumber, (\ref{w<-->theta}) turns into
  \begin{equation}\label{w<-->theta-in-Fourier}
  \begin{array}{rclc}
  \left(|k|^2-\frac{d^2}{dz^2}\right)^2\F w&=&|k|^2 \F\theta  \qquad & {\rm for } \quad 0<z<H\,,\\
  \F w&=&\frac{d}{d z}\F w=0 \qquad & {\rm for } \quad z\in\{0,H\}\,,\\
  \end{array}
  \end{equation}
  whereas (\ref{stab-cond-re}) assumes the form
  \begin{equation}\label{stab-cond-fourier-temp}
  2\int_0^H \xi  \F{w}\overline{\F{\theta}}\, dz+\int_0^H |k|^2|\F\theta|^2 dz+\int_0^H \left|\frac{d}{d z}\F\theta\right|^2dz\geq 0\,,
  \end{equation}
 where the bar denotes complex conjugation.
  Using equation (\ref{w<-->theta-in-Fourier}) we can eliminate $\theta$ from the stability condition
   (\ref{stab-cond-fourier-temp}), obtaining
  \begin{eqnarray}\label{OSCinFourier}
  &&2\int_0^H\xi \F w \left(-\frac{d^2}{dz^2}+|k|^2\right)^2\overline{\F w}\,dz\\
  &&+\int_0^H|k|^{-2}\left|\frac{d}{dz}\left(-\frac{d^2}{dz^2}+|k|^2\right)^2\F w\right|^2\,dz
  +\int_0^H\left|\left(-\frac{d^2}{dz^2}+|k|^2\right)^2\F w\right|^2\,dz\;\ge\;0\,,\notag
  \end{eqnarray}
   which has to be satisfied for all $k\in\frac{2\pi}{L}\mathbb{Z}^{d-1}\setminus\{0\}$ and all (complex valued) functions $\F w(z)$
   satisfying the three boundary conditions
   \begin{equation}\label{bc}
   \F w\;=\;\frac{d}{dz}\F w\;=\;\left(-\frac{d^2}{dz^2}+|k|^2\right)^2\F w\;=\;0
   \quad\mbox{for}\;z\in\{0,H\}.
   \end{equation}
   We now introduce a further simplification by letting $L\uparrow \infty$ so that (\ref{OSCinFourier})
   has to hold for all $k\in \R^d$. This strengthening of the stability condition has the additional 
   advantage that it becomes independent of the dimension $d$: We will henceforth say that $\xi$ satisfies
   the stability condition if 
  \begin{eqnarray}
  &&2\int_0^H\xi  w \left(-\frac{d^2}{dz^2}+k^2\right)^2 \overline{ w}\,dz\nonumber\\
  &&+\int_0^H k^{-2}\left|\frac{d}{dz}\left(-\frac{d^2}{dz^2}+k^2\right)^2 w\right|^2\,dz+\int_0^H\left|\left(-\frac{d^2}{dz^2}+k^2\right)^2 w\right|^2\,dz\;\ge\;0\,,\label{OSCinFourier-reduct}
  \end{eqnarray}
  holds true for all $k\in\mathbb{R}$ and all (complex valued) functions $ w(z)$
  satisfying the three boundary conditions (\ref{bc}) with $\F w$ replaced by $w$.
  The analysis of the stability condition (\ref{OSCinFourier-reduct}) imposed for all values of $k\in \R$ (which corresponds to assume $L\uparrow\infty$)
  amounts to consider profiles $\tau$ that are stable even under perturbation that have horizontal wave-length much larger than $H$.
  In \cite{Thesis-Nobili} it is shown that at least Proposition \ref{prop1} still holds true if the lateral size $L$ is of order $H$ . 
  
  The rest of the paper is organized as follow: In Section \ref{heuristics} we study a reduced stability condition (obtained by retaining only the indefinite term in the stability condition)
   and show that  in this case a stable profile must be increasing and logarithmically growing. This result is obtained by exploring the limit 
  for small and large wavelengths in the reduced stability condition, respectively. 
  In Section \ref{Lemmas}, when working with the original stability condition we  can no longer pass to the limit for small/large wavenumbers $k$ to infer the positivity for $\xi=\frac{d\tau}{dz}$ and the logarithmic growth for $\tau$. 
  Nevertheless, by subtle averaging of the stability condition we construct a non-negative convolution kernel $\phi_0$ with help of which we can express the positivity on average approximately in the bulk (see Lemma \ref{Lemma1}).
  Likewise we recover logarithmic growth, at least approximately in the bulk, on the level of the construction of $\xi_0$, see Lemma \ref{Lemma2}.
  Finally in Lemma \ref{Lemma3} we connect the bulk with the boundary layers. 
  The main result (Theorem \ref{Th1}) is proved in Section \ref{proof-M-T} and it consists of combining all the results contained in the lemmas together with an estimate that connects $\xi_0$ to $\xi$, and in particular to $\int_0^H\xi\, dz=-1$ (see Lemma \ref{Lemma4}, estimate (\ref{CLAIM-1})).
  
  In the rest of the paper we omit the constant factor $2$ in front of the indefinite term in (\ref{OSCinFourier-reduct}), which is legitimate since we are interested in the scaling of the Nusselt number.
   
  \section{Characterization of stable profiles}\label{original-stability-section}
   
  \subsection{Reduced stability condition}\label{heuristics}
  We note that the stability condition (\ref{OSCinFourier-reduct})\&(\ref{bc}) is invariant
  under the following transformation
  \begin{equation}\label{invariance-under-scaling}
  z=L\hat z \; \mbox{ and thus }\; k=\frac{1}{L}\hat k,\; H=L\hat H \;\mbox{ and }\; \xi=L^{-4}\hat{\xi}.
  \end{equation}
  Hence in the bulk ($z\gg 1$ and $H-z\gg 1$) we expect that the first term in (\ref{OSCinFourier-reduct}) dominates.
  This motivates to consider the \textit{reduced stability condition}
 \begin{equation}\label{RSC}
  \int_0^H\xi w\left(-\frac{d^2}{dz^2}+k^2\right)^2 \overline{w} \, dz\geq 0,
  \end{equation}
  for all $k\in \R$ and all (complex valued) functions $w(z)$ satisfying the three boundary conditions (\ref{bc}).
  
  The following proposition is independent of the main result (Theorem \ref{Th1}) but it serves as a preparation: 
  The ideas developed in the proof (cf.\ Section \ref{PP1}) will be adapted to the more challenging  full stability condition (\ref{OSCinFourier-reduct}).
  \begin{proposition}\label{prop1}
   Let  $\xi=\xi(z)$ be such that for all $k\in \R$ and for all $w(z)$ satisfying (\ref{bc}),
   it satisfies the reduced stability condition (\ref{RSC}) .
   Then
  \begin{equation}\label{red-positivity}\xi\geq 0\,,  \qquad \mbox{ and }\end{equation}
   \begin{equation}\label{red-log-growth}\int_{1/e}^{1} \xi dz \lesssim  \frac{1}{\ln H}\int_1^{H} \xi \,dz \,.\end{equation}
 
  \end{proposition}
  \noindent
  We notice that while (\ref{red-positivity}) means that $\tau$ is an increasing functions, the second statement (\ref{red-log-growth})
  corresponds to a logarithmic growth of $\tau$ (see Figure $1$). Hence somewhat surprisingly, monotonicity is not sufficient for stability.
  \begin{figure}[ht!]

  \centering
  \begin{tikzpicture}

  \draw[->] (0,0) -- (10,0) node[above] {$z$};
  \draw[->] (0,-1) -- (0,4) node[right] {$\tau$};

  \draw (1/e+0.3,-1) node[above] {$\frac{1}{e}$};
  \draw  (3,-1)    node[above]{$1$};
  \draw  (9.5,-1)   node[above]{$H$};

  \draw [blue,domain=0.5:9.5,smooth] plot (\x, {ln(\x+1/2)+1});
  \draw[dashed] (9.5,3) -- (9.5,0);
  \draw[dashed] (3,2.25) -- (3,0);
  \draw[dotted] (3.1,2.25) -- (9.5,2.25);
  \draw[dotted] (3,1.155)--(1/e+0.3, 1.155);
  \draw[dashed] (1/e+0.3, 1.155) -- (1/e+0.3,0);
  \end{tikzpicture}
  \caption{Logarithmic growth of the background profile $\tau$ as expressed in Proposition \ref{prop1}.}
  \end{figure}
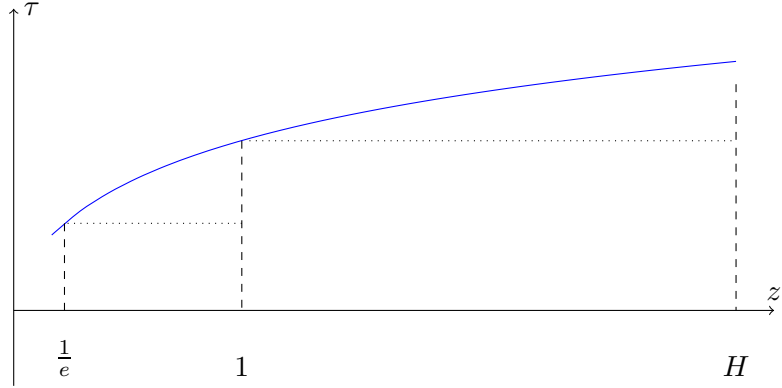
  \subsection{Original stability condition: statement of lemmas}\label{Lemmas}

  In the following lemmas, we derive properties of those profiles $\tau$ that, in terms of their slope $\xi=\frac{d\tau}{dz}$, 
  satisfy the original stability condition (\ref{OSCinFourier-reduct})
  and (for the last lemma) the driving boundary conditions (\ref{IC}).
  These four lemmas are the (only) ingredients of the main theorem.
  They all are formulated on the level of the logarithmic variables $s=\ln z$ and 
  $\hat\xi=z\xi=\frac{d\tau}{ds}$, cf.\ (\ref{new-variables}). Lemma \ref{Lemma1} establishes approximate positivity
  of the slope $\hat\xi$ in the bulk, and thus is the generalization of (\ref{red-positivity}) in Proposition \ref{prop1},
  replacing the stricter reduced stability
  condition (\ref{RSC}) there by the original stability condition (\ref{OSCinFourier-reduct}) here. It does so in terms of a suitable convolution
  $\hat\xi_0$ of $\hat\xi$ in the logarithmic variable $s$.
  Lemma \ref{Lemma2} establishes approximate logarithmic growth of the profile in the bulk, again on the level of $\hat\xi_0$,
  and amounts to the generalization of (\ref{red-log-growth}) in Proposition \ref{prop1}. Lemma \ref{Lemma3} is the most subtle and 
  shows that the convolved slope $\hat\xi_0$ cannot be too negative in the boundary layer $-s\gg 1$ provided
  it is sufficiently small in the transition layer $|s|\lesssim 1$.
  Lemma \ref{Lemma4} translates the driving boundary conditions (\ref{IC}) on $\tau$ in form of $\int_0^H\xi\, dz=-1$
  from the slope $\xi$ to its logarithmic-variable convolution $\hat\xi_0$.
  \begin{lemma} \label{Lemma1}\ \\
  There exists a $\phi_0$, which will play the role of a convolution kernel, with the properties
  \begin{equation}\label{convolution-kernel}
   \phi_0(s)\geq 0, \quad \int_{-\infty}^{\infty}\phi_0(s) \,ds=1,\quad \mathrm{supp}\, \phi_0(z)\subset \left(\frac{1}{4},\frac{3}{4}\right), \quad \phi_0\left(\frac 12-z\right)=\phi_0\left(\frac 12+z\right)\,,
  \end{equation}
  such that, for all $s'\leq \ln H$
  \begin{equation}\label{1}
   \hat{\xi}_0(s') \gtrsim -\exp(-3s')\,,
  \end{equation}
  where 
   \begin{equation}\label{convolution}
  \hat\xi_0(s'):=\int_{-\infty}^{\infty}\hat\xi(s+s')\phi_0(s)\,ds\,.
  \end{equation}
  \end{lemma}
%
  \begin{lemma}\label{Lemma2}\ \\
  For $S_1\gg 1$ we have 
  \begin{equation}\label{2}
  \int_{-1}^{0}\hat \xi_0 ds\lesssim\frac{1}{S_1}\int_{0}^{S_1}\hat \xi_0\, ds+1\,.
 \end{equation}
  \end{lemma}
  \begin{lemma}\label{Lemma3}\ \\
  For all $S_2\gg 1$ and $\varepsilon \leq 1$ we have 
  \begin{equation}\label{3}
 \int_{-S_2}^{-1}\hat \xi_0 ds\gtrsim -\left(\frac{1}{\varepsilon}\int_{-1}^{0}\hat\xi_0\, ds+\frac{1}{\varepsilon}+\int_{-S_2}^{-S_2+1}|\hat\xi_0|ds+\varepsilon \exp(5S_2)\right)\,.
  \end{equation}
  \end{lemma}
%
  \begin{lemma}\label{Lemma4}\ \\
  Suppose that the slope $\xi$ of the profile $\tau$ satisfies $\int_0^H\xi\, dz=-1$ and $\int_0^H\xi^2\, dz\lesssim (\ln H)^{\frac{1}{15}}$. Then
 \begin{equation}\label{CLAIM-1}
 \int_{-\infty}^{\ln H}\hat \xi_0 \,ds\lesssim -1\,.
 \end{equation}
  \end{lemma}

\subsection{Proof of Proposition \ref{prop1}}\label{PP1}
  \underline{Argument for (\ref{red-positivity}):}
  
  \noindent
  Letting $k\uparrow \infty$, (\ref{RSC}) reduces to 
  $$\int_{0}^{H}\xi|w|^2\, dz\geq 0\, $$
  for all compactly supported $w$,
  from which we infer (\ref{red-positivity}).
 
   \noindent
  \underline{Argument for (\ref{red-log-growth}),} heuristic version:
  
  \noindent
  Letting $k\downarrow 0$, (\ref{RSC}) reduces to
  \begin{equation}\label{SSC2}
  \int_0^{H}\xi w\frac{d^4}{dz^4} \overline{w} dz\geq 0\,
  \end{equation}
  for all functions $w(z)$ satisfying the three boundary conditions 
  \begin{equation}\label{bc-lg}
   w=\frac{dw}{dz}=\frac{d^4 w}{d z^4}=0 \qquad \mbox{ for } z\in \{0,H\}\,.
  \end{equation}
  In fact, besides Subsection \ref{Proof-Lemma3}, we will work with $w$ compactly supported in $z\in (0,H)$, so that the boundary condition
  (\ref{bc-lg}) are trivially satisfied.
  Focusing on the lower half of the container
  we make the following Ansatz
  $$w=z^2\hat w\,,$$
  where $\hat w(z)$ is a real function with compact support  in $(0,H)$. 
  The merit of this Ansatz is that in the new variable $\hat w$, the multiplier in (\ref{SSC2})
  can be written in the scale-invariant form 
  \begin{equation}\label{poly}
  \phi=w\frac{d^4}{dz^4}\bar w=
  \hat w z^2\frac{d^4}{dz^4}z^2\hat w=\hat w\left(z\frac{d}{dz}+2\right)\left(z\frac{d}{dz}+1\right)z\frac{d}{dz}\left(z\frac{d}{dz}-1\right)\hat w\,.
  \end{equation}
  Note that the fourth-order polynomial in $z\frac{d}{dz}$ appearing on the r.\ h.\ s.\ of (\ref{poly}) may be inferred, 
  without lengthy calculations, from the fact that $z^2\frac{d^4}{dz^4}z^2$ annihilates $\{\frac{1}{z^2},\frac{1}{z},1,z\}$.
  This suggests to introduce the new variables 
  \begin{equation}\label{new-variables}
   s=\ln z \, \mbox{ and } \, \xi=z^{-1}\hat \xi,
  \end{equation}
  for which the stability condition turns into 
  \begin{equation}\label{stability-conditions-in-new-variables}
   \int_{-\infty}^{\ln H}\hat \xi\, \phi \,ds\, \geq 0\,,
   \end{equation} where
  $$ \phi=\hat w\left(\frac{d}{ds}+2\right)\left(\frac{d}{ds}+1\right)\frac{d}{ds}\left(\frac{d}{ds}-1\right)\hat w\,,$$
 for all functions $\hat w$ with compact support in $z\in (0,H)$.
 Here comes the heuristic argument for (\ref{red-log-growth}):
 For $H\gg1$, we may think of test functions $\hat w$ that vary slowly in the logarithmic variable $s$.
 For these $\hat w$ we have 
 \begin{equation}\label{heu}
 \phi=
 \hat w\left(\frac{d}{ds}+2\right)\left(\frac{d}{ds}+1\right)\frac{d}{ds}\left(\frac{d}{ds}-1\right)\hat w\approx -2\hat w\frac{d}{ds}\hat w=
 -\frac{d}{ds}\hat w^2\,,
 \end{equation}
 which in particular implies
 $$0\leq \int_{-\infty}^{\ln H}\hat \xi\,\phi\, ds\,\approx -\int_{-\infty}^{\ln H}\hat \xi\frac{d}{ds}\hat w^2 \,ds=\int_{-\infty}^{\ln H}\frac{d\hat \xi}{ds}\hat w^2\, ds\,.$$
 Since $\hat w$ was arbitrary besides varying slowly in $s$, it follows 
 $$\frac{d\hat \xi}{ds}\geq 0\,,$$
  approximately on large $s-$scales.
 We expect that this implies that for any $1\ll S_1\leq \ln H$:
 \begin{equation}\label{LOG-gro}
  \int_{-1}^{0}\hat \xi ds \lesssim  \frac{1}{S_1}\int_0^{S_1}\hat \xi \,ds \,,
 \end{equation}
 which in the original variables (\ref{new-variables}), for $S_1$ turns into (\ref{red-log-growth}).
  We now establish rigorously that 
 (\ref{SSC2}) and (\ref{red-positivity}) imply (\ref{LOG-gro}).
                            
  \noindent
  \underline{Argument for (\ref{LOG-gro}),} rigorous version:

  \noindent
  We start by noticing that because of translation invariance in $s$,  
  (\ref{stability-conditions-in-new-variables}) can be reformulated as follows: For any function
  $\hat w(s)$ supported in $s\le 0$, and any $s'\le \ln H$ we have
  \begin{equation}\label{trasl-inv}
  \int_{-\infty}^\infty\hat\xi(s'')\phi(s''-s')\,ds''\;=\;
  \int_{-\infty}^{\infty}\hat\xi(s+s')\phi(s)\,ds\;\ge\;0,
  \end{equation}
  where the multiplier $\phi$ is defined as in (\ref{stability-conditions-in-new-variables}).
  We note that (\ref{LOG-gro}) follows from (\ref{trasl-inv})
  once for given $S_1$ we construct
  \begin{itemize}
  \item a family ${\mathfrak F}=\{w_{s'}\}_{s'}$
  of smooth functions $w_{s'}$ parameterized by $s'\in\mathbb{R}$
  and compactly supported in $z\in(0,1)$ (i.\ e.\ $s\in(-\infty,0]$) and
  \item a measure $\rho(ds')=\rho(s')\,ds'$ supported in $s'\in(-\infty,\ln H]$,
  \end{itemize}
  such that the corresponding convex combination of multipliers $\{\phi_{s'}\}_{s'}$
  shifted by $s'$, i.\ e.\
  \begin{equation}\label{phi-first}
  \phi_1(s'')\;:=\;\int_{-\infty}^\infty\phi_{s'}(s''-s')\,\rho(s')ds'\,,
  \end{equation}
  satisfies
  \begin{equation}\label{phi1-upper-bound}
  \phi_1(s'') \leq
  \left\{\begin{array}{lll}
  -1 & \mbox{for } -1\leq s''\leq 0\, \\
  \frac{C}{S_1} & \mbox{for } \;\;\;\; 0\leq s''\leq S_1\, \\
  0 & \mbox{else} \,
  \end{array}\right\}\,,
  \end{equation}
  for a (possibly large) universal constant $C$.
  Indeed, using (\ref{trasl-inv}), (\ref{phi-first}), and (\ref{phi1-upper-bound}) in conjunction with the positivity 
  (\ref{red-positivity}) of the profile $\hat{\xi}$ we have
  \[0\leq \int_{-\infty}^{\infty}\hat{\xi}\phi_{1}ds''\leq -\int_{-1}^{0}\hat{\xi}ds''+\frac{C}{S_1}\int^{S_{1}}_{0}\hat{\xi}ds'',\]
which implies (\ref{LOG-gro})\,.

  We first address the form of the family ${\mathfrak F}$.
  The heuristic observation (\ref{heu}) motivates the change of variables
  \begin{equation}\label{change-s}
  s\;=\;\lambda\hat s\quad\mbox{with}\quad\lambda\geq 1,
  \end{equation}  
  our ``(logarithmic) length-scale'', to be chosen sufficiently large. 
  We fix a smooth, real-valued and compactly supported ``mask'' $\hat w_0(\hat s)$; 
  it will be convenient to restrict its support to $\hat s\in(-1,0]$, say
  \begin{equation}\label{support}
  \hat w_0^2>0 \quad \mbox{ in }\left(-\frac{1}{2},0\right) \quad \mbox{ and } \quad\hat w_0=0 \quad\quad \mbox{ else },
  \end{equation} 
  and, in order to justify the language of ``mollification by convolution''
  we impose the normalization $\int\hat w_0^2d\hat s=1$.
  By mask we mean that in (\ref{stability-conditions-in-new-variables}) we choose
  \begin{equation}\label{mask-w0}
  \hat w(\lambda \hat s)\;=\;\lambda^{-1/2}\hat w_0(\hat s)
  \end{equation}
  (see Figure $2$).
   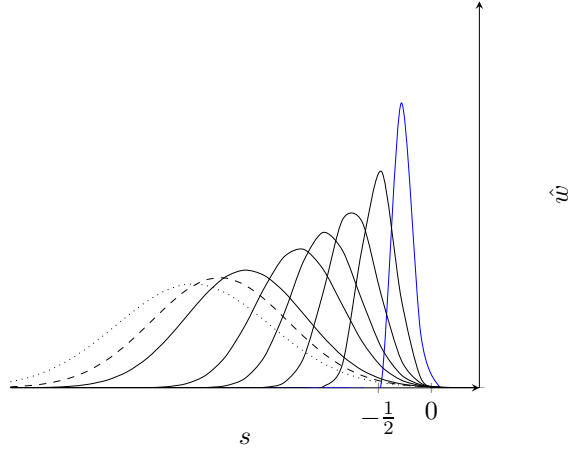
\begin{figure}[ht!]
   \centering
   \begin{tikzpicture}[scale=0.9]
  \begin{axis}[ mark=none,
  axis x line=left,
  axis y line=right,
     xlabel={$ s$},
     ylabel={$\hat{w}$},
     domain=-35:4, smooth, xmin=-35, xmax=4, ymin=0, ymax=0.7, xtick={-4.4,0},ytick={0}, yticklabels={$ $}, xticklabels={$-\frac{1}{2}$,$0$}]
  \addplot[mark=none, blue]{1/(sqrt(pi))*exp(-((x+2.2)^2))};
  \addplot[mark=none]{(1/sqrt(2))*(1/(sqrt(pi)))*exp(-((0.5*x+2.2)^2))};
  \addplot[mark=none]{(1/sqrt(3))*(1/(sqrt(pi)))*exp(-((1/3)*x+2.2)^2))};
  \addplot[mark=none]{(1/sqrt(4))*(1/(sqrt(pi)))*exp(-((0.25*x+2.2)^2))};
   \addplot[mark=none]{(1/sqrt(5))*(1/(sqrt(pi)))*exp(-(((1/5)*x+2.2)^2))};
   \addplot[mark=none]{(1/sqrt(7))*(1/(sqrt(pi)))*exp(-(((1/7)*x+2.2)^2))};
    \addplot[mark=none, dashed]{(1/sqrt(8))*(1/(sqrt(pi)))*exp(-(((1/8)*x+2.2)^2))};
     \addplot[mark=none, dotted]{(1/sqrt(9))*(1/(sqrt(pi)))*exp(-(((1/9)*x+2.2)^2))};
  \end{axis}
  \end{tikzpicture}
  \caption{ Construction of the family of test functions starting from the mask $\hat w_0$ (blue line).}
  \end{figure}
  With this change of variables, the multiplier can be rewritten as follows
  \begin{eqnarray}\nonumber
  \lefteqn{\phi_{\lambda}\;\stackrel{(\ref{heu})}{=}\;\hat w\left(\frac{d}{ds}+2\right)\left(\frac{d}{ds}+1\right)
  \frac{d}{ds}\left(\frac{d}{ds}-1\right)\hat w}\nonumber\\
  &\stackrel{(\ref{change-s})\&(\ref{mask-w0})}{=}&\frac{1}{\lambda}\hat w_0\left(\frac{1}{\lambda}\frac{d}{d\hat s}+2\right)
  \left(\frac{1}{\lambda}\frac{d}{d\hat s}+1\right)
  \frac{1}{\lambda}\frac{d}{d\hat s}
  \left(\frac{1}{\lambda}\frac{d}{d\hat s}-1\right)\hat w_0
  \nonumber\\
  &=&\hat{w}_{0}\left(\frac{1}{\lambda^5}\frac{d^4}{d\hat{s}^4}+
  \frac{2}{\lambda^4}\frac{d^3}{d\hat{s}^3}-
  \frac{1}{\lambda^3}\frac{d^2}{d\hat{s}^2}-
  \frac{2}{\lambda^2}\frac{d}{d\hat{s}}\right)\hat{w}_{0}\nonumber\,,
  \end{eqnarray}
  and reordering the terms we have
  \begin{equation}\label{phi-in-lambda}
  \phi_{\lambda}=
  -\frac{2}{\lambda^2}\hat{w}_{0}\frac{d}{d\hat{s}}\hat{w}_{0}-
  \frac{1}{\lambda^3}\hat{w}_{0}\frac{d^2}{d\hat{s}^2}\hat{w}_{0}+
  \frac{2}{\lambda^4}\hat{w}_{0}\frac{d^3}{d\hat{s}^3}\hat{w}_{0}+
  \frac{1}{\lambda^5}\hat{w}_{0}\frac{d^4}{d\hat{s}^4}\hat{w}_{0}\,.
  \end{equation}
  Heuristically, for $\lambda\gg 1$ the multiplier $\phi_{\lambda}$ can be approximated by the first term on the r.\ h.\ s.\ 
  $$\phi_{\lambda}(s)\approx -\frac{d}{ds}\left(\frac{1}{\lambda}\hat w_0^2\left(\frac{s}{\lambda}\right)\right).$$
  Inserting this approximation into the definition (\ref{phi-first}) of $\phi_1$ we have 
   \begin{eqnarray}\label{after-45}
  \phi_{1}(s'')
  &=&\int_{-\infty}^{\infty}\phi(s''-s')\rho(s')ds'=\int_{-\infty}^{\infty}\phi(s)\rho(s''-s)ds\nonumber\\
  &=&-\int_{-\infty}^{\infty}\frac{d}{d s} \left(\frac{1}{\lambda}\hat w_0^2\left(\frac{s}{\lambda}\right)\right)\rho(s''-s)ds
  \approx-\int_{-\infty}^{\infty}\frac{1}{\lambda}\hat w_0^2\left(\frac{s}{\lambda}\right)\frac{d\rho}{d s'}(s''-s)ds\,.
   \end{eqnarray}
  For $\lambda$ smaller than the characteristic scale on which $\rho$ varies, 
  we may think of $\frac{1}{\lambda}\hat w_0^2\left(\frac{s}{\lambda}\right)\approx \delta_0(s)$, in view of our normalization.
   This yields 
  \begin{equation}\label{phi1-behav}
  \phi_1\approx -\frac{d\rho}{ds'}\,,
  \end{equation}
  which in view of (\ref{phi1-upper-bound}) suggests that $\rho$ should have the form
  \begin{equation}\label{form}
  \rho(s')=
   \left\{\begin{array}{lll}
   s'+1 & \mbox{ for } -1\leq s'\leq 0\,\\
   1-\frac{s'}{S_1} & \mbox{ for } 0\leq s'\leq S_1\,
   \end{array}\right\}\,,
  \end{equation}
  (see Figure $3$).
  We will now argue that we have to modify the Ansatz for both $\phi_{\lambda}$ and $\rho$.
      \begin{figure}[ht!]
    \centering
    \begin{tikzpicture}[scale=0.8]
   \begin{axis}[ 
   axis x line=left,
   axis y line=left,
      xlabel={$s'$},
      ylabel={$\rho$},
      domain=-3:10, smooth, xmin=-1, xmax=8, ymin=0, ymax=1.3,xtick={-1,0,8},ytick={0},yticklabels={$ $}, xticklabels={$-1$,$0$,$S_1$}]
   \addplot[mark=none, domain=-1:0]{x+1};
   \addplot[mark=none, domain=0:8]{1-x/8};
    \end{axis}
   \end{tikzpicture}
   \caption{The measure $\rho$ suggested by the heuristic argument. }
   \end{figure}
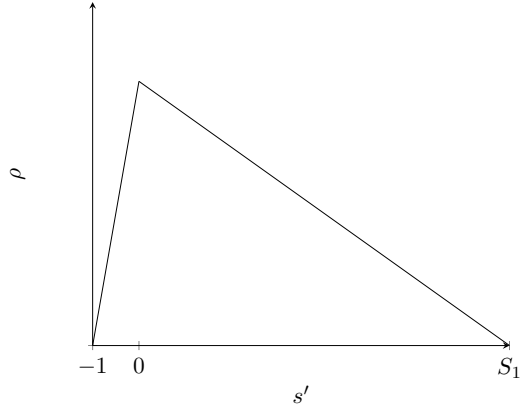
  To this purpose, we go through the above heuristic argument (heuristically) assessing the error terms.
  Expanding $\rho$ in a Taylor series around $s''$
  \[\rho(s''-s)\;\approx\; \rho(s'')-\frac{d\rho}{ds'}(s'')s+\frac{1}{2}\frac{d^2\rho}{ds'^2}(s'')s^2\,,\]
  we may write
  \begin{eqnarray}\nonumber
  \phi_{1}(s'')&\stackrel{(\ref{phi-first})}{=}&\int_{-\infty}^{\infty}\phi_{\lambda}(s)\rho(s''-s)ds\;\; \nonumber\\
  &\approx& \rho(s'')\int_{-\infty}^{\infty} \phi_{\lambda}\,ds-\frac{d\rho}{ds'}(s'')\int_{-\infty}^{\infty} s\phi_{\lambda}\,ds+\frac{1}{2}\frac{d^2\rho}{ds'^2}(s'')\int_{-\infty}^{\infty} s^2 \phi_{\lambda}\,ds\,.\nonumber
  \end{eqnarray}
  We note that the first term in (\ref{phi-in-lambda}), i.\ e.\ $-\frac{2}{\lambda^2}\hat{w}_{0}\frac{d\hat w_0}{d\hat s}=-\frac{1}{\lambda^2}\frac{d\hat{w}_{0}^2}{d\hat{s}}$, gives the leading-order contribution to the first and the second moment
  \begin{equation*}
  \int_{-\infty}^\infty s\phi_{\lambda}\,ds
  \approx
  \int_{-\infty}^\infty s\left(-\frac{1}{\lambda^2}\frac{d\hat w_0^2}{d\hat s}\right)\,ds
  \stackrel{(\ref{change-s})}{=}\int_{-\infty}^\infty \hat s\left(-\frac{d\hat w_0^2}{d\hat s}\right)\,d\hat s
  \;=\;\int_{-\infty}^\infty \hat{w}_0^2\,d\hat s=1
  \end{equation*}
  and 
  \begin{eqnarray*}
  \int_{-\infty}^\infty s^2\phi_{\lambda}\,ds\approx
  \int_{-\infty}^\infty s^2\left(-\frac{1}{\lambda^2}\frac{d\hat w_0^2}{d\hat s}\right)\,ds
  \stackrel{(\ref{change-s})}{=} -\lambda \int_{-\infty}^\infty\hat s^2\left(\frac{d\hat w_0^2}{d\hat s}\right)\,d\hat s=
  \lambda \int_{-\infty}^\infty 2\hat s\hat w_0^2\,d\hat s\,,
  \end{eqnarray*}
  while the second term in (\ref{phi-in-lambda}) gives the leading-order contribution to the zeroth moment of the multiplier $\phi$:

  \begin{equation*}
   \int_{-\infty}^{\infty}\phi_{\lambda}\,ds\approx
  \int_{-\infty}^{\infty}\left(-\frac{1}{\lambda^3}\hat w_0\frac{d^2\hat w_0}{d\hat s^2}\right)\,d s
  \stackrel{(\ref{change-s})}{=} \frac{1}{\lambda^2}\int_{-\infty}^{\infty}\left(\frac{d\hat w_0}{d\hat s}\right)^2\, d\hat s\,.
  \end{equation*}
  Hence we obtain the following specification of (\ref{after-45})
  \begin{eqnarray}\label{K0}
  \phi_{1}(s'')&=&\int_{-\infty}^{\infty}\phi_{\lambda}(s)\rho(s''-s)ds\; \nonumber\\
  &\approx& \frac{1}{\lambda^2}\rho(s'')\int_{-\infty}^{\infty}\left(\frac{d\hat w_0}{d\hat s}\right)^2\, d\hat s
  -\frac{d\rho}{ds'}(s'')\;
  +\lambda \frac{d^2\rho}{ds'^2}(s'')\int_{-\infty}^\infty\hat{s} \hat{w}_0^2\,d\hat s\,.
  \end{eqnarray}

  Our goal is to specify the choice (\ref{form}) of $\rho$ such that (\ref{phi1-upper-bound}) is satisfied.
  This shows a dilemma: On the one hand, in the ``plateau region'' $s''\sim S_1$, we would need $\lambda^2\gg S_1$
  so that the first r.\ h.\ s.\ term in (\ref{K0}) does not destroy the desired $\frac{1}{S_1}-$behavior. On the other hand 
  in the ``foot region'' $s''\in[0,1]$, we would need $\lambda\lesssim 1$ so that the last term does not destroy the 
  effect of the middle term.
 This suggests that $\lambda$ should be chosen to be small in the foot regions and large on the plateau region.
 Therefore it is natural to choose 
 \begin{equation}\label{lambda=s'}
 \lambda=s'\,,
 \end{equation}
 so that $\phi_{s'}$ in (\ref{phi-first}) indeed acquires a dependency on $s'$ besides the translation.
%
  For our choice of (\ref{lambda=s'}), (\ref{phi-in-lambda}) assumes the form  
  \begin{equation}\label{phi-in-s-prime}
  \phi_{s'}=
  -\frac{2}{(s')^2}\hat{w}_{0}\frac{d}{d\hat{s}}\hat{w}_{0}-
  \frac{1}{(s')^3}\hat{w}_{0}\frac{d^2}{d\hat{s}^2}\hat{w}_{0}+
  \frac{2}{(s')^4}\hat{w}_{0}\frac{d^3}{d\hat{s}^3}\hat{w}_{0}+
  \frac{1}{(s')^5}\hat{w}_{0}\frac{d^4}{d\hat{s}^4}\hat{w}_{0}\,.
  \end{equation}
  Note that with the choice (\ref{lambda=s'}) and $s=s''-s'$, (\ref{change-s}) turns into the \textit{nonlinear} change of variables
  \begin{equation}\label{nonlinear-change}
  \hat s=\frac{s''-s'}{s'}=\frac{s''}{s'}-1\Rightarrow s'=\frac{s''}{1+\hat s}\,.
  \end{equation}
  We consider this as a change of variables between $s'$ and $\hat s$ (with $s''$ as a parameter);
  thanks to the support restriction (\ref{support}) on $\hat w_0$, it is invertible in 
  the relevant range $\hat s\in[-\frac 12,0]$:
  $\frac{d}{d\hat s}=-\frac{s''}{(1+\hat s)^2}\frac{d}{ds'}=-\frac{(s')^2}{s''}\frac{d}{ds'}$ and 
  $ds'=\frac{s''}{(1+\hat s)^2}d\hat s$. From (\ref{phi-first}) and (\ref{phi-in-s-prime}) we thus get the first representation
  \begin{eqnarray*}
  \phi_1(s'')&=&-\frac{1}{s''}\int_{-\infty}^\infty\frac{d\hat w_0^2}{d\hat s}\; \rho \;d\hat s
  -\frac{1}{(s'')^2}\int_{-\infty}^\infty(1+\hat s)\,\hat w_0\frac{d^2\hat w_0}{d\hat s^2  }
  \rho\,d\hat s\\
  &+&\frac{2}{(s'')^3}\int_{-\infty}^\infty(1+\hat s)^2\,\hat w_0\frac{d^3\hat w_0}{d\hat s^3  }
  \rho\,d\hat s
  +\frac{1}{(s'')^4}\int_{-\infty}^\infty(1+\hat s)^3\,\hat w_0\frac{d^4\hat w_0}{d\hat s^4  }
  \rho\,d\hat s\,.
  \end{eqnarray*}
  An approximation argument in $\hat w_0$ below necessitates 
  a second representation that involves $\hat w_0$ only up to second derivatives.
  For this purpose, we rewrite (\ref{phi-in-s-prime}) in terms of the three quadratic quantities
  $\hat w_0^2$, $(\frac{d\hat w_0}{d\hat s})^2$, and $(\frac{d^2\hat w_0}{d\hat s^2})^2$:
  \begin{eqnarray}\label{phi-in-s-prime-2}
  \phi_{s'}
  &=&
  -\frac{1}{(s')^2}\frac{d\hat w_0^2}{d\hat s  }
  +\frac{1}{(s')^3}\left[\left(\frac{d\hat w_0}{d\hat s}\right)^2
                    -\frac{1}{2}\frac{d^2\hat w_0^2}{d\hat s^2}\right]
  +\frac{1}{(s')^4}\left[-3\frac{d}{d\hat s}\left(\frac{d\hat w_0}{d\hat s}\right)^2
                    +\frac{d^3\hat w_0^2}{d\hat s^3}\right]\nonumber\\
  &+&\frac{1}{(s')^5}\left[ \left(\frac{d^2\hat w_0}{d\hat s^2}\right)^2
                    -2\frac{d^2}{d\hat s^2}\left(\frac{d\hat w_0}{d\hat s}\right)^2
                    +\frac{1}{2}\frac{d^4\hat w_0^2}{d\hat s^4}\right]\nonumber\\
  &=&
  \left(-\frac{1}{(s')^2}\frac{d  }{d\hat s  }
      -\frac{1}{2}\frac{1}{(s')^3}\frac{d^2}{d\hat s^2}
      + \frac{1}{(s')^4}\frac{d^3}{d\hat s^3}
      +\frac{1}{2}\frac{1}{(s')^5}\frac{d^4}{d\hat s^4}\right)\hat w_0^2\nonumber\\
  &+&\left(\frac{1}{(s')^3}
         -3\frac{1}{(s')^4}\frac{d  }{d\hat s  }
         -2\frac{1}{(s')^5}\frac{d^2}{d\hat s^2}\right)\left(\frac{d\hat w_0}{d\hat s}\right)^2
         +\frac{1}{(s')^5}\left(\frac{d^2\hat w_0}{d\hat s^2}\right)^2\,.
  \end{eqnarray}
  Now in this formula, using the change of variables (\ref{nonlinear-change}),
  we want to substitute the derivations
  $\frac{1}{(s')^m}\frac{d^n}{d\hat s^n}$ by linear combinations of derivations of the form
  $\frac{1}{(s'')^{m-k}}\frac{d^k}{d s'^k}(1+\hat s)^{m-n-k}$ for $k=0,\cdots,n$.
  The reason why this can be done is explained in Appendix \ref{Appendix1}, where also the linear 
  combinations are explicitly computed.
  The formulas (\ref{tranformation-formula}), (\ref{A1}) \& (\ref{A2}) allow to rewrite (\ref{phi-in-s-prime-2}) as follows
  \begin{eqnarray}\label{new-star}
  \phi_{s'}
  &=&\frac{1}{s''}\left(\frac{d  }{ds'}-\frac{1}{2}\frac{d^2}{ds'^2}\frac{1}{(1+\hat s)}-\frac{d^3}{ds'^3}\frac{1}{(1+\hat s)^2}+\frac{1}{2}\frac{d^4}{ds'^4}\frac{1}{(1+\hat s)^3}\right)\hat w_0^2\nonumber\\
  &+&\bigg[  \left(\frac{1}{(s'')^3}+\frac{6}{(s'')^4}-\frac{12}{(s'')^5}\right)(1+\hat s)^3+\left(\frac{3}{(s'')^3}-\frac{8}{(s'')^4}\right)\frac{d  }{ds'  }(1+\hat s)^2\nonumber\\
  &-& \frac{2}{(s'')^3}\frac{d^2}{ds'^2}(1+\hat s)\bigg]\left(\frac{d\hat w_0}{d\hat s}\right)^2+\frac{1}{(s'')^5}(1+\hat s)^5\left(\frac{d^2\hat w_0}{d\hat s^2}\right)^2\,.
  \end{eqnarray}
  The advantage of this form is that integrations by part in $s'$ become easy, so that we obtain
  \begin{eqnarray}
  \phi_1
  &=&\frac{1}{s''}\int_{-\infty}^\infty\hat w_0^2\left(-\frac{d  \rho}{ds'  }-\frac{1}{2}\frac{1}{ 1+\hat s}\frac{d^2\rho}{ds'^2}+\frac{1}{(1+\hat s)^2}\frac{d^3\rho}{ds'^3}+\frac{1}{2}\frac{1}{(1+\hat s)^3}\frac{d^4\rho}{ds'^4} \right)\,ds'\nonumber\\
  &+&\left(\frac{1}{(s'')^3}+\frac{6}{(s'')^4}-\frac{12}{(s'')^5}\right)\int_{-\infty}^\infty(1+\hat s)^3\left(\frac{d\hat w_0}{d\hat s}\right)^2 \rho\,ds'\nonumber\\
  &-&\left(\frac{3}{(s'')^3}-\frac{8}{(s'')^4}\right)\int_{-\infty}^\infty(1+\hat s)^2\left(\frac{d\hat w_0}{d\hat s}\right)^2\frac{d  \rho}{ds'  }\,ds'\nonumber \\
  &-&\frac{2}{(s'')^3}\int_{-\infty}^\infty(1+\hat s)  \left(\frac{d\hat w_0}{d\hat s}\right)^2\frac{d^2\rho}{ds'^2}\,ds'\nonumber\\
  &+&\frac{1}{(s'')^5}\int_{-\infty}^\infty(1+\hat s)^5\left(\frac{d^2\hat w_0}{d\hat s^2}\right)^2\rho\,ds'\,. \nonumber
  \end{eqnarray}
  Finally, using the substitution $\frac{ds'}{s''}=\frac{d\hat s}{(1+\hat s)^2}$, the last
  formula turns into the desired second representation
  \begin{eqnarray}
  \phi_1
  &=&\int_{-\infty}^\infty\hat w_0^2\left(-\frac{1}{(1+\hat s)^2}\frac{d  \rho}{ds'  } -\frac{1}{2}\frac{1}{(1+\hat s)^3}\frac{d^2\rho}{ds'^2}+\frac{1}{(1+\hat s)^4}\frac{d^3\rho}{ds'^3}+\frac{1}{2}\frac{1}{(1+\hat s)^5}\frac{d^4\rho}{ds'^4} \right)\,d\hat s\nonumber\\
  &+&\left(\frac{1}{(s'')^2}+\frac{6}{(s'')^3}-\frac{12}{(s'')^4}\right)\int_{-\infty}^\infty(1+\hat s)\left(\frac{d\hat w_0}{d\hat s}\right)^2 \rho\,d\hat s \nonumber\\
  &-&\left(\frac{3}{(s'')^2}-\frac{8}{(s'')^3}\right)\int_{-\infty}^\infty\left(\frac{d\hat w_0}{d\hat s}\right)^2\frac{d  \rho}{ds'  }\,d\hat s\nonumber\\
  &-&\frac{2}{(s'')^2}\int_{-\infty}^\infty\frac{1}{1+\hat s}\left(\frac{d\hat w_0}{d\hat s}\right)^2\frac{d^2\rho}{ds'^2}\,d\hat s \nonumber\\
  &+&\frac{1}{(s'')^4}\int_{-\infty}^\infty(1+\hat s)^3\left(\frac{d^2\hat w_0}{d\hat s^2}\right)^2\rho\,d\hat s\,.\label{phi_1-formula}
  \end{eqnarray}
  From this representation we learn the following:
  If we assume that  $\rho(s')$ varies on large length-scales only, so that 
   $\frac{d\rho}{ds'}$,$\frac{d^2\rho}{ds'^2},\cdots \ll \rho$
  and  $\frac{d^2\rho}{ds'^2}$,$\frac{d^3\rho}{ds'^3},\cdots \ll \frac{d\rho}{ds'}$
   then for $s''\gg 1$, we obtain to leading order from the above
  \begin{equation*}
  \phi_1\approx-\int_{-\infty}^\infty\frac{1}{(1+\hat s)^2}\hat w_0^2\,\frac{d\rho}{ds'}\,d\hat s+\frac{1}{(s'')^2}\int_{-\infty}^\infty(1+\hat s) \left(\frac{d\hat w_0}{d\hat s}\right)^2\rho \,d\hat s\,.
  \end{equation*}
  If $\rho(s')$ varies slowly even on a logarithmic scale (so that e.\ g.\
  $s'\frac{d\rho}{ds'}$ is negligible with respect to $\rho$), the above further reduces to
  \begin{equation}\label{reduction}
  \phi_1\approx
  -\frac{d\rho(s'')}{ds'}\int_{-\infty}^\infty\frac{1}{(1+\hat s)^2}\hat w_0^2\,d\hat s
  +\frac{\rho(s'')}{(s'')^2}
  \int_{-\infty}^\infty(1+\hat s)\left(\frac{d\hat w_0}{d\hat s}\right)^2\,d\hat s\,,
  \end{equation}
 which should be compared with (\ref{K0}). We see that the first, negative, 
  right-hand-side term of (\ref{reduction}) dominates the second positive term provided
  $$\frac{d\rho}{ds'}\;\gg\;\frac{1}{(s')^2}.$$
  This is satisfied if $\rho$ is of the form
  $\rho(s')=1-\frac{S_{0}}{s'-S_{0}}$
  for some $S_0\gg 1$ to be chosen later;
  indeed $\frac{d\rho}{ds'}=\frac{S_0}{(s'-S_0)^2}\approx \frac{S_0}{(s')^2}\gg\frac{1}{(s')^2}$  for $s'\gg S_0\gg1$.
  
  Disregarding for a couple of pages the fact that $\rho$ needs to be supported in $s'\in (-\infty,\ln H)$, which will be 
  achieved by cutting off at scales $s'\sim S_1$, we define as our intermediate goal to 
   construct a measure $0\leq\trho(s')\leq 1$ (with infinite support) such that 
   \begin{equation}\label{reduced-goal}
  \widetilde{\phi}_1(s''):=\int \phi_{s'}(s''-s')\trho(s')\, ds'
                                 \left\{\begin{array}{ll}
                                  =0 & s''\leq \frac{S_0}{2}\,\\[1mm]
                                  <0 & s''>\frac{S_0}{2}\,                              
                                   \end{array}\right\}\,,
   \end{equation}
  for some universal $S_0$, to be chosen later.
  The above considerations motivate the following Ansatz for $\trho$: 
  We fix a smooth mask $\trho_0(\hat s')$ such that
  \begin{equation}\label{rho0-mask}
  \trho_0\;=\;0\;\;\mbox{for}\;\hat s'\le 0,\quad
  \frac{d\trho_0}{d\hat s'}\;>\;0\;\;\mbox{for}\;0<\hat s'\le 2,\quad
  \trho_0\;=\;1-\frac{1}{\hat s'}\;\;\mbox{for}\;2\le \hat s'\,,
  \end{equation}
     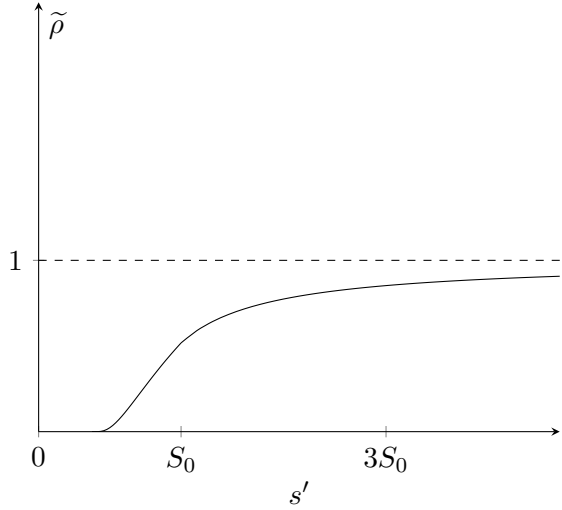
\begin{figure}[ht!]
    \centering
    \begin{tikzpicture}
     \begin{axis}[ mark=none,
     axis x line=left,
      axis y line=center,
      xlabel={$ s'$},
      ylabel={$\trho$},
      domain=0:35,xmin=0, xmax=30, ymin=0, ymax=1, smooth, xtick={0,8.2,20},ytick={0,0.399}, xticklabels={$0$,$S_0$,$3S_0$},yticklabels={$0$,$1$}]
     \addplot[mark=none, domain=0:8.2]{(exp(-5/(x-3)^(0.7))};
     \addplot[mark=none, domain=8.2:30]{0.399-1/(x-3)};
     \addplot[mark=none, dashed, domain=0:15]coordinates {(0,0.399) (30,0.399)};
      \end{axis}
     \end{tikzpicture}
     \caption{The function $\tilde{\rho}$ of the variable $s'=S_0(\hat s'+1).$}
     \end{figure}
  and consider the rescaled version 
  \begin{equation}\label{rescale-rho}
  \trho(S_0(\hat s'+1))\;=\;\trho_0(\hat s'),\quad\mbox{i.\ e.\ the change of variables}\;s'\;=\;S_0(\hat s'+1)\,
  \end{equation}
   with $S_0\gg 1$ to be fixed later (see Figure $4$).
%
%
  It is convenient to rescale $s''$ accordingly:
  \begin{equation}\label{rescale-s''}
  s''\;=\;S_0\hat s''.
  \end{equation}
  With this new rescaling and the choice of the measure $\trho$ (see (\ref{rho0-mask}) and (\ref{rescale-rho})), (\ref{phi_1-formula}) turns into
  \begin{eqnarray}
  \widetilde{\phi}_1
  &=&-\frac{1}{ S_0  }\int_{-\infty}^\infty
  \frac{\hat w_0^2}{(1+\hat s)^2}\frac{d  \trho_0}{d\hat s'  }\,d\hat s 
   -\frac{1}{2S_0^2}\int_{-\infty}^\infty
  \frac{\hat w_0^2}{(1+\hat s)^3}\frac{d^2\trho_0}{d\hat s'^2}\,d\hat s\nonumber\\
  && +\frac{1}{ S_0^3}\int_{-\infty}^\infty
  \frac{\hat w_0^2}{(1+\hat s)^4}\frac{d^3\trho_0}{d\hat s'^3}\,d\hat s
   +\frac{1}{2S_0^4}\int_{-\infty}^\infty
  \frac{\hat w_0^2}{(1+\hat s)^5}\frac{d^4\trho_0}{d\hat s'^4}\,d\hat s
  \nonumber\\
  &&+\left(\frac{1}{S_0^2}\frac{1}{(\hat s'')^2}+\frac{1}{S_0^3}\frac{6}{(\hat s'')^3}
   -\frac{1}{S_0^4}\frac{12}{(\hat s'')^4}\right)
  \int_{-\infty}^\infty(1+\hat s)\left(\frac{d\hat w_0}{d\hat s}\right)^2 \trho_0\,d\hat s
  \nonumber\\&&
  -\left(\frac{1}{S_0^4}\frac{3}{(\hat s'')^3}-\frac{1}{S_0^5}\frac{8}{(\hat s'')^4}\right)
  \int_{-\infty}^\infty  \left(\frac{d\hat w_0}{d\hat s}\right)^2\frac{d \trho_0}{d\hat s'  }\,d\hat s
  \nonumber\\&&
  -\frac{1}{S_0^5}\frac{2}{(\hat s'')^3}
  \int_{-\infty}^\infty\frac{1}{1+\hat s}\left(\frac{d\hat w_0}{d\hat s}\right)^2\frac{d^2\trho_0}{d\hat s'^2}\,d\hat s
  \nonumber\\
  &&+\frac{1}{S_0^4}\frac{1}{(\hat s'')^4}\int_{-\infty}^\infty(1+\hat s)^3\left(\frac{d^2\hat w_0}{d\hat s^2}\right)^2
  \trho_0\,d\hat s\,.\label{phi_1-new}
  \end{eqnarray}
  Since in the integrals in formula (\ref{phi_1-formula}), 
  the argument of $\trho$ was given by $s'=\frac{s''}{1+\hat s}$,
  cf.\ (\ref{nonlinear-change}), it follows from (\ref{rescale-rho}) and (\ref{rescale-s''}) that
  the argument of $\trho_0$ is given by
  \begin{equation}\label{argument-of-new-rho}
  \hat s'\;=\;\frac{\hat s''}{1+\hat s}-1\,.
  \end{equation}
  Thus all the integrals in (\ref{phi_1-new}) just depend on $\hat s''$, not on $S_0$.
  Hence (\ref{phi_1-new}) makes the dependence of $\widetilde{\phi}_1$ on $S_0$ explicit.
  We are now ready to show that the construction of the family $w_{s'}$ (cf.\ (\ref{mask-w0}) and (\ref{lambda=s'})) 
  and of the measure $\trho$ (cf.\ (\ref{rho0-mask}) and (\ref{rescale-rho})) yield the intermediate goal (\ref{reduced-goal}) when $S_0\gg 1$.
%
%
   In order to establish (\ref{reduced-goal}) it is convenient to distinguish three regions 
  (note that if $s''\in (\infty,\frac{S_0}{2}]$ all the integrals in (\ref{phi_1-new}) vanish
  because the supports of $\hat w_0$ and $\trho_0$ do not intersect, see below): 

  \noindent
  \underline{The range of large $s''$:}

  \begin{equation}\label{large-s''}
  s''\ge 3S_0\quad\mbox{or equivalently}\quad\hat s''\ge 3.
  \end{equation}	
  Note that because of our support condition (\ref{support}) on $\hat w_0$,
  all integrals in (\ref{phi_1-new}) are supported in $\hat s\in[-\frac{1}{2},0]$.
  Together with our range (\ref{large-s''}), this yields for the argument 
  $\hat s'\stackrel{(\ref{argument-of-new-rho})}{=}\frac{\hat s''}{1+\hat s}-1$ of
  $\trho_0$ and its derivatives that $\hat s'\ge 2$.
  Because of $\frac{d\trho_0}{d\hat s'}=\frac{1}{(\hat s')^2}$ for $\hat s'\ge 2$, cf.\
  our Ansatz (\ref{rho0-mask}), the first integral in (\ref{phi_1-new}) reduces to
  \begin{multline}\label{first-term}
  \int_{-\infty}^\infty\frac{\hat w_0^2}{(1+\hat s)^2}\frac{d  \trho_0}{d\hat s'  }\,d\hat s 
  =\int_{-\infty}^\infty\frac{\hat w_0^2}{(1+\hat s)^2}\frac{1}{(\frac{\hat s''}{1+\hat s}-1)^2}\,d\hat s\\
  =\int_{-\infty}^\infty\frac{\hat w_0^2}{(\hat s''-(1+\hat s))^2}\,d\hat s
  \approx\frac{1}{(\hat s'')^2}\int_{-\infty}^\infty \hat w_0^2\,d\hat s\, ,
  \end{multline}
  for $\hat s''\gg 1$, whereas
  all the other integrals in (\ref{phi_1-new}) are $O(\frac{1}{(\hat s'')^2})$ or smaller in $\hat s''\gg 1$
  (because at least one derivation falls on $\widetilde\rho_0$)
  or have pre-factors $\frac{1}{(\hat s'')^2}$ or smaller. Since only the
  term in (\ref{phi_1-new}) coming from integral (\ref{first-term}) has pre-factor $\frac{1}{S_0}$ 
  while all the other terms have pre-factors $\frac{1}{S_0^2}$ or smaller (for $S_0\gg 1$),
  the first term in (\ref{phi_1-new}) {\it uniformly} dominates all other terms for $S_0\gg 1$:
  \begin{equation}\label{dominant-term-in-small-range}
  \widetilde{\phi}_1\;\approx\;
  -\frac{1}{S_0}\int_{-\infty}^\infty\frac{\hat w_0^2}{(\hat s''-(1+\hat s))^2}\,d\hat s
  \quad\mbox{uniformly in}\;\hat s''\ge 3\quad\mbox{for}\;S_0\gg 1\,.
  \end{equation}
  In conclusion we have 
  \begin{equation}\label{negativity-in-big-range}
   \widetilde{\phi}_1\sim -\frac{1}{S_0}\frac{1}{(\hat s'')^2}<0 \quad \mbox{ in the range } \hat s''\ge 3 \quad \;\mbox{ for }\;S_0  \gg 1\,.
  \end{equation}

  \noindent
  \underline{The range of intermediate $s''$:}

  \begin{equation}\label{intermediate-s''}
  s''\;\in\;\left[\frac{3}{4}S_0,3S_0\right]\quad\mbox{or equivalently}\quad\hat s''\;\in\;\left[\frac{3}{4},3\right].\end{equation}
  Again, we consider the first integral in (\ref{phi_1-new}). Now we use that
  $\frac{\hat w_0^2}{(1+\hat s)^2}\ge 0$ is {\it strictly positive} in $\hat s\in\left(-\frac{1}{2},0\right)$,
  cf.\ (\ref{support}), and that $\frac{d\trho_0}{d\hat s'}\ge 0$ is strictly positive
  in $\hat s'>0$, cf.\ (\ref{rho0-mask}), that is, in $\hat s<\hat s''-1$, cf.\ (\ref{argument-of-new-rho}).
  We note that the two $\hat s-$intervals $(-\frac{1}{2},0)$ and $(-\infty,\hat s''-1)$
  intersect for $\hat s''>\frac{1}{2}$. Hence by continuity of the first integral
  in (\ref{phi_1-new}) in its parameter $\hat s''$, there exists a universal constant $C$ such that
  \begin{equation*}
  \int_{-\infty}^\infty\frac{\hat w_0^2}{(1+\hat s)^2}\,\frac{d\trho_0}{d\hat s'}\,d\hat s
  \;\ge\;\frac{1}{C}\quad\mbox{for}\;\hat s''\in\left[\frac{3}{4},3\right]\,.
  \end{equation*}
  Hence also in this range the first term in (\ref{phi_1-new}) dominates all other terms:
  \begin{equation}\label{dominant-term-in-itermediate-range}
  \widetilde{\phi}_1\;\approx\;
  -\frac{1}{S_0}\int_{-\infty}^\infty\frac{\hat w_0^2}{(1+\hat s)^2}\,\frac{d\trho_0}{d\hat s'}\,d\hat s
  \quad\mbox{uniformly in}\;\hat s''\in\left[\frac{3}{4},3\right]\quad\mbox{for}\;S_0\gg 1\,,
  \end{equation}
  and we may conclude that 
  \begin{equation}\label{negativity-in-intermediate-range}
  \widetilde{\phi}_1\sim-\frac{1}{S_0}<0 \quad \mbox{ in the range }  s''\in\left[\frac{3}{4}S_0,3S_0\right]\quad \;\mbox{ for } S_0 \; \gg 1\,.
  \end{equation}
  Note that the above discussion on supports also yields that $\widetilde{\phi}_1$ is supported
  in $\hat s''\in\left[\frac{1}{2},\infty\right)$.

  \noindent
  \underline{The range of small $s''$:}
  \begin{equation}\label{small-s''}
  s''\;\in\;\left(\frac{1}{2}S_0,\frac{3}{4}S_0\right)\quad\mbox{or equivalently}\quad
  \hat s''\;\in\;\left(\frac{1}{2},\frac{3}{4}\right)\,.
  \end{equation}
  We would like $\widetilde{\phi}_1$ to be strictly negative in this range for $S_0\gg 1$.
  Here, we encounter the second difficulty: No matter how large $\lambda=s'$
  in (\ref{phi-in-lambda}) is, the behavior of $\phi_{s'}$ near the left edge $-\frac{1}{2}$ of
  its support $\left[-\frac{1}{2},0\right]$ (and also at its right edge $0$, but there
  we don't care), is dominated by the $\frac{1}{\lambda^5}\hat w_0\frac{d^4\hat w_0}{d\hat s^4}$-term
  and thus automatically is {\it strictly positive}. Taking the $\trho(s')\,ds'$-average
  of the shifted $\phi_{s'}(s''-s')$ does not alter this behavior as long as
  $\trho$ is non-negative in $[S_0,\infty)$, cf.\ (\ref{rho0-mask}): $\widetilde{\phi}_1$ is
  strictly positive near the left edge $\frac{S_0}{2}$ of its support. 
  The way out of this problem is to give  {\it give up smoothness} of $\hat w_0$ near the
  left edge $-\frac{1}{2}$ of its support $\left[-\frac{1}{2},0\right]$. In fact, we shall first
  assume that $\hat w_0$ satisfies in addition
  \begin{equation}\label{w_0-near-support}
  \hat w_0\;=\;\frac{1}{2}\left(\hat s+\frac{1}{2}\right)^2\quad\mbox{for}\;\hat s\in\left[-\frac{1}{2},-\frac{1}{4}\right]\,.
  \end{equation}
  This means that $\hat w_0$ has a bounded but discontinuous second derivative (i.\ e.\ $\hat w_0\in H^{2,\infty}$). This
  is the main reason why in (\ref{phi_1-new}) we expressed $\widetilde{\phi}_1$ only in terms of up to second derivatives
  of $\hat w_0$. We argue that the so defined $\widetilde{\phi}_1$ is, as desired, strictly 
  negative on $s''\in(\frac{1}{2}S_0,\frac{3}{4}S_0]$ for all $S_0$. 
%
  Indeed, in view of (\ref{phi-in-s-prime}), the form (\ref{w_0-near-support}) implies, in terms of $s=s''-s'$,
  \begin{equation}\label{phi-s'-with-choice-w0-near-support}
  \phi_{s'}\;=\;-\frac{1}{(s')^2}\left(\frac{s}{s'}+\frac{1}{2}\right)^3
              -\frac{1}{2}\frac{1}{(s')^3}\left(\frac{s}{s'}+\frac{1}{2}\right)^2\;<\;0
  \quad\mbox{for}\;s\in\left(-\frac{s'}{2},-\frac{s'}{4}\right]\,.
  \end{equation}
  In view of (\ref{rho0-mask}) \& (\ref{rescale-rho}), $\trho\ge 0$ is strictly positive for $s'\in(S_0,\infty)$.
  On the other hand, it follows from (\ref{phi-s'-with-choice-w0-near-support}) that
  $s'\mapsto \phi_{s'}(s''-s')$ is strictly negative for 
  $s''-s'\in(-\frac{s'}{2},-\frac{s'}{4}]$, that is, for $s'\in[\frac{4}{3}s'',2s'')$
  (and supported in $s'\in[s'',2s'']$). Hence, by (\ref{phi-first}), $\widetilde{\phi}_1$ is strictly negative for
  $S_0\in[\frac{4}{3}s'',2s'')$, that is, for $s''\in(\frac{1}{2}S_0,\frac{3}{4}S_0]$, for any value of $S_0>0$. 

  Now we approximate $\hat w_{0}$ with a sequence of smooth functions $\hat w_{0}^{\nu}$ in $H^{2,2}$ and we call
  $\widetilde{\phi}_{1}^{\nu}$ the associated multiplier.
  Since $\widetilde{\phi}_{1}$ involves $\hat w_{0}$ only up to second derivatives (cf.\ (\ref{phi_1-new}))
  then $\phi_{1}^{\nu}$ converge uniformly in $\hat{s}''$ to 
  $\widetilde{\phi}_{1}$. We conclude that
  \begin{equation}\label{negativity-in-small-range}
  \widetilde{\phi}_1<0 \quad \mbox{ for } s''\;\in\;\left(\frac{1}{2}S_0,\frac{3}{4}S_0\right)\;\mbox{ and }\;S_0>0\,.
  \end{equation}
  Finally we fix a sufficiently large but universal $S_0$, 
  so that (\ref{negativity-in-small-range}) together with (\ref{negativity-in-big-range}) and
  (\ref{negativity-in-intermediate-range}) imply our intermediate goal (\ref{reduced-goal}).

  The choice (\ref{rho0-mask})\&(\ref{rescale-rho}) of $\trho$ is not admissible, since 
  $\rho$ should be supported in $(-\infty,\ln H)$, which will be achieved by cutting off at scales $s'\sim S_1$.
  Only this cut-off will ensure (\ref{phi1-upper-bound}) in the range $s''\geq S_1$. More precisely, 
  in (\ref{phi-first}) we choose $\rho$ to be $\trho(s')\eta(\frac{s'}{S_1})$ where 
  $\trho$ is defined in (\ref{rescale-rho}), $\eta$ is a mask for a smooth cut-off function with
  $$\eta(\hat{\hat s}')=\left\{\begin{array}{ll}
                         1 &   \hat{\hat s}'\leq \frac 12\,\\[1mm]
                         0 &   \hat{\hat s}'\geq 1\, 
                          \end{array}\right\}\,          $$ 
   (see Figure $5$ for an illustration of $\rho$).
   \begin{figure}[ht!]
    \centering
    \begin{tikzpicture}
   \begin{axis}[ mark=none,
   axis x line=left,
   axis y line=left,
      xlabel={$s'$},
      ylabel={$\rho$},
      domain=0:35,xmin=0, xmax=50, ymin=0, ymax=1, smooth, xtick={0,8.2,20, 40},ytick={0,0.399}, xticklabels={$0$,$S_0$,$\frac{S_1}{2}$, $S_1$},
      yticklabels={$0$,$1$}
     ]
   \addplot[mark=none, domain=0:8.2]{(exp(-5/(x-3)^(0.7))};
   \addplot[mark=none, domain=8.2:20.2]{0.399-1/(x-3)};
   \addplot[mark=none, dashed, domain=0:35] coordinates {(0,0.399) (40,0.399)};
    \addplot[mark=none, smooth, red] coordinates {   (20.2, 0.34) (22, 0.32) (25, 0.25) (35.9, 0.02) (48, 0)};
   \node[red] at (300,20) [anchor=south] {$\eta$};
    \end{axis}
   \end{tikzpicture}
   \caption{The measure $\rho$ is constructed from $\widetilde{\rho}$ (see Figure $4$) by cutting off at scales $s'\sim S_1$ .}
   \end{figure}
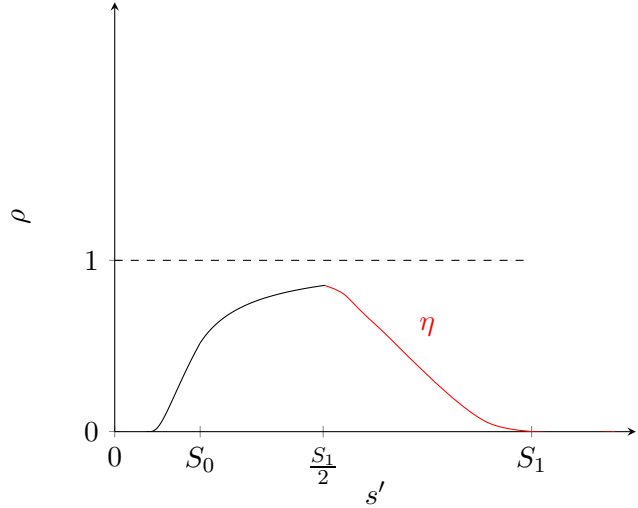
    We argue that $\phi_1(s'')$ satisfies (\ref{phi1-upper-bound}) for $S_1\gg 1$.
  This will follow immediately from the three claims
  \begin{enumerate}
   \item $\phi_1(s'')=0$ for $s''\geq S_1$,
   \item $\phi_1(s'')=\widetilde{\phi}_1(s'')$ for $s''\leq \frac{S_1}{4}$,
   \item $|\phi_1-\widetilde{\phi}_1|\lesssim\frac{1}{S_1}$ for $\frac{S_1}{4}\leq s''\leq S_1$.
  \end{enumerate}
  Claims $1$, $2$ and $3$ together with the bound (\ref{reduced-goal}) on $\widetilde{\phi}_1$ imply
    \begin{equation}\label{76}
   \phi_1(s'') 
    \left\{\begin{array}{llll}
   &=0 & \mbox{for } &  s''\leq \frac{S_0}{2} \\[2mm]
    &\lesssim -1 & \mbox{ for } & \frac{S_0}{2}+1\leq s''\leq \frac{S_0}{2}+2\\[2mm]
    & \leq 0 & \mbox{ for } & \frac{S_0}{2}\leq s''\leq \frac{S_1}{4}\\[2mm]
   &\lesssim \frac{1}{S_1} & \mbox{ for } &\frac{S_1}{4}\leq s''\leq S_1\\[2mm]
   &=0 & \mbox{ for } & s''\geq S_1 \,
   \end{array}\right\}.
   \end{equation}
  Note that with help of the scale invariance (\ref{invariance-under-scaling}), which turns into a shift invariance in the logarithmic
  variables (\ref{new-variables}), we may shift $\phi_1$ and its estimate (\ref{76}) by $\frac{S_0}{2}+1$ to the left (and redefine $S_1\gg S_0$
  noting that the universal constant $S_0$ was fixed in the previous step)
  recovering the desired form (\ref{phi1-upper-bound}).
  Let us now establish the claims $1$, $2$, and $3$. We start by noting that by the change of variables (\ref{nonlinear-change}) in conjunction
  with the support condition (\ref{support}) on $\hat s$ we have the following relation between
  the argument $s'$ of $\rho$ and the variable $s''$ of $\phi_1$ in the representation (\ref{phi_1-formula}) 
  \begin{align}\label{fo1}
  \left\{\begin{array}{ccccc}
  &&s''&\le&\frac{S_1}{4}\,\\
  \frac{S_1}{4}&\le&s''&\le&S_1\,\\
  S_1&\le&s''& &
  \end{array}\right\}
  \quad\Rightarrow\quad
  \left\{\begin{array}{ccccc}
  &&s'&\le&\frac{S_1}{2}\,\\
  \frac{S_1}{4}&\le&s'&\le&2S_1\,\\
  S_1&\le&s'& &\,
  \end{array}\right\}.
  \end{align}
  Hence Claims $1$ and $2$ are immediate consequences of the defining properties of the cut-off function $\eta$.
  We now turn to Claim $3$ and appeal to the representation (\ref{phi_1-formula}) for $\phi_1$, which we also use for $\widetilde{\phi}_1$ and
  thus obtain a representation of $\widetilde{\phi}_1-\phi_1$ with $\trho$ replaced by $(1-\eta)\trho$. 
  Clearly, all the terms bearing a pre-factor of $\frac{1}{{s''}^2}$ or smaller are of higher order in the range $s''\ge\frac{S_1}{4}$.
  Likewise, all the terms where at least one of the derivative on the product $(1-\eta)\trho$ 
  falls on the cut-off $1-\eta$, thereby producing a factor $\frac{1}{S_1}$,
  are of the desired order or smaller. We are left with the terms
  \begin{align}\label{fo2}
  \int_{-\infty}^\infty\hat w_0^2(1-\eta)\left(-\frac{1}{(1+\hat s)^2}\frac{d\trho}{ds'}
  -\frac{1}{2}\frac{1}{(1+\hat s)^3}\frac{d^2\trho}{{ds'}^2}
  +\frac{1}{(1+\hat s)^4}\frac{d^3\trho}{{ds'}^3}
  +\frac{1}{2}\frac{1}{(1+\hat s)^5}\frac{d^2\trho}{{ds'}^4}\right)d\hat s.
  \end{align}
  In our range $\frac{S_1}{4}\le s''\le S_1$ the integrand is supported in $\frac{S_1}{4}\le s'\le 2 S_1$, cf.\ (\ref{fo1}).
  For these arguments we have by choice (\ref{rho0-mask}) \& (\ref{rescale-rho}) of the averaging function $\trho(s')=1-\frac{1}{\frac{s'}{S_0}-1}
  \approx 1-\frac{S_0}{s'}$. Hence also the terms (\ref{fo2}) are at least of order $\frac{1}{S_1^2}$ and thus
  of higher order. This establishes Claim $3$ and thus (\ref{76}).

\section{Proof of the main theorem}\label{proof-M-T}

  {\it Proof of Theorem 1}.\ \\
  We start by combining Lemmas \ref{Lemma2} and \ref{Lemma3}. We first note that by rearranging (\ref{2}) in Lemma \ref{Lemma2} and then
  adding $\int_{-1}^0\hat\xi_0\,ds$ to both sides, while using $S_1\gg 1$, 
  we obtain by rearranging
  \begin{align*}
  -\int_{-1}^{S_1}\hat\xi_0\,ds\le S_1\left(-\frac{1}{C_1}\int_{-1}^0\hat\xi_0\,ds+1\right),
  \end{align*}
  where we momentarily retain the value $C_1$ of the universal constant.
  We now add this to (\ref{3}) in Lemma \ref{Lemma3} in form of
  \begin{align*}
  -\int_{-S_2}^{-1}\hat\xi_0\,ds\le C_2\left(\frac{1}{\epsilon}\int_{-1}^0\hat\xi_0\,ds+\frac{1}{\epsilon}
  +\int_{-S_2}^{-S_2+1}|\hat\xi_0|\,ds+\epsilon\exp(5S_2)\right)
  \end{align*}
  and adjust $\epsilon=\frac{C_1 C_2}{S_1}$ so that the pre-factor of
  the transition term $\int_{-1}^0\hat\xi_0\,ds$ vanishes; this choice of $\epsilon$
  is consistent with $\epsilon\le 1$ because of $S_1\gg 1$. We end up with
  \begin{align*}
  -\int_{-S_2}^{S_1}\hat\xi_0\,ds\lesssim S_1+\int_{-S_2}^{-S_2+1}|\hat\xi_0|\,ds+\frac{1}{S_1}\exp(5S_2),
  \end{align*}
  to which we add $-\int_{-\infty}^{-S_2}\hat\xi_0\,ds\le\int_{-\infty}^{-S_2}|\hat\xi_0|\,ds$, obtaining
  \begin{align}\label{fo3}
  -\int_{-\infty}^{S_1}\hat\xi_0\,ds\lesssim \int_{-\infty}^{-S_2}|\hat\xi_0|\,ds+S_1+\frac{1}{S_1}\exp(5S_2),
  \end{align}
  where we replaced $S_2$ by $S_2+1$ (at the expense of changing the multiplicative constant hidden in $\lesssim$).

  In order to provide ourselves with an additional parameter $S_0$ next to $S_1,S_2\gg 1$ to optimize in,
  we now appeal to a scaling argument: The change of variables (\ref{invariance-under-scaling}) that leaves the stability
  condition (\ref{OSCinFourier-reduct}) invariant assumes the form of a shift
  \begin{align*}
  s=S_0+\hat s,\quad\hat\xi=\exp(-3S_0)\hat{\hat\xi}\quad\mbox{and thus}\quad
  \hat\xi_0=\exp(-3S_0)\hat{\hat{\xi_0}}.
  \end{align*}
  on the level of the logarithmic variables (\ref{new-variables}).
  We apply (\ref{fo3}), which only relied on the stability condition, to the rescaled variables
  $\hat s$, $\hat{\hat{\xi_0}}$ (with the upper integral boundaries $\hat S_1$ and $-\hat S_2$), 
  and then revert to the original variables:
  \begin{align*}
  -\int_{-\infty}^{\hat S_1+S_0}\hat\xi_0\,ds\lesssim\int_{-\infty}^{-\hat S_2+S_0}|\hat\xi_0|\,ds+
  \exp(-3S_0)\left(\hat S_1+\frac{1}{\hat S_1}\exp(5\hat S_2)\right).
  \end{align*}
  This estimate holds provided $\hat S_1,\hat S_2\gg 1$.
  In terms of the original integral boundaries $S_1=\hat S_1+S_0$ and $-S_2=-\hat S_2+S_0$ this reads
  \begin{equation}\label{fo8}
  -\int_{-\infty}^{S_1}\hat\xi_0\,ds
  \lesssim\int_{-\infty}^{-S_2}|\hat\xi_0|\,ds+
  \exp(-3S_0)\left(S_1-S_0+\frac{1}{S_1-S_0}\exp(5(S_2+S_0))\right),
  \end{equation}
  which is valid provided $S_1-S_0=\hat S_1\gg 1$, $S_2+S_0=\hat S_2\gg 1$, and $S_1\le\ln H$. Now is the
  moment to optimize in $S_0$ by choosing it such that the two last terms are balanced,
  which is achieved by $S_1-S_0=\exp(\frac{5}{2}(S_2+S_0))$. In our regime $S_2+S_0\gg 1$ we have 
  $\exp(\frac{5}{2}(S_2+S_0))\approx\exp(\frac{5}{2}(S_2+S_0))-(S_2+S_0)$ so that the above
  choice means $S_1+S_2\approx\exp(\frac{5}{2}(S_2+S_0))$, which implies $\exp(-3S_0)$ $
  \approx\exp(3S_2)(S_1+S_2)^{-\frac{6}{5}}$. Hence with our choice, the entire second
  term in (\ref{fo8}) is $\approx\exp(3S_2)(S_1+S_2)^{-\frac{1}{5}}$:
  \begin{align}\label{fo9}
  -\int_{-\infty}^{S_1}\hat\xi_0\,ds
  \lesssim\int_{-\infty}^{-S_2}|\hat\xi_0|\,ds+\exp(3S_2)(S_1+S_2)^{-\frac{1}{5}}.
  \end{align}
  This estimate is valid provided $S_1\le\ln H$ and $S_1+S_2\gg 1$, since the latter by
  $S_1+S_2\approx\exp(\frac{5}{2}(S_2+S_0))$ ensures $S_2+S_0\gg 1$ and by
  $S_1-S_0=\exp(\frac{5}{2}(S_2+S_0))$ then also $S_1-S_0\gg 1$. 

  We now make use of Lemma \ref{Lemma1} and Lemma \ref{Lemma4}. Note that w.\ l.\ o.\ g.\ we may assume that
  our background profile $\tau$ satisfies on the level of its slope 
  $\int_0^H\xi^2\,dz\lesssim(\ln H)^\frac{1}{15}$ next to $\int_0^H\xi\, dz=-1$,
  since if such a profile would not exist, the statement of Theorem 1 would
  be trivially true. Hence we are in the position to apply Lemma \ref{Lemma4}.
  By Lemma \ref{Lemma1} we have $-\int_{S_1}^{\ln H}\hat\xi_0\,ds$ $\lesssim\exp(-3S_1)$. Adding
  this to (\ref{fo9}), we obtain 
  \begin{align*}
  -\int_{-\infty}^{\ln H}\hat\xi_0\,ds\lesssim \int_{-\infty}^{-S_2}|\hat\xi_0|\,ds
  +\exp(3S_2)(S_1+S_2)^{-\frac{1}{5}}+\exp(-3S_1),
  \end{align*}
  so that by (\ref{CLAIM-1}) in Lemma \ref{Lemma4} we get
  \begin{align*}
  1\lesssim \int_{-\infty}^{-S_2}|\hat\xi_0|\,ds
  +\exp(3S_2)(S_1+S_2)^{-\frac{1}{5}}+\exp(-3S_1).
  \end{align*}
  Clearly, the optimal choice of $S_1$ is given by saturating the constraint in form of $S_1=\ln H$, 
  which by $H\gg 1$ turns into
  \begin{align}\label{fo4}
  1\lesssim \int_{-\infty}^{-S_2}|\hat\xi_0|\,ds+\exp(3S_2)(\ln H)^{-\frac{1}{5}}
  \end{align}
  and holds provided $S_2\ge 0$.

  We finally connect this to the Nusselt number ${\rm Nu}_{\BF}$, 
  which on the level of the slope $\xi$ and
  the logarithmic variables turns into
  \begin{align}\label{fo5}
  {\rm Nu}_{\BF}\ge\int_0^H(\frac{d\tau}{dz})^2\,dz=\int_0^H\xi^2\,dz=\int_{-\infty}^{\ln H}\exp(-s)\hat\xi^2\,ds.
  \end{align}
  This allows us to estimate the first r.\ h.\ s.\ term in (\ref{fo4}): By definition (\ref{convolution}) of the convolution
  $\hat\xi_0$ and the support property (\ref{convolution-kernel}) of the kernel $\phi_0$ we have
  \begin{align}\label{fo6}
  \int_{-\infty}^{-S_2}|\hat\xi_0|\,ds\lesssim\int_{-\infty}^{-S_2}|\hat\xi|\,ds
  \lesssim\exp(-\frac{1}{2}S_2)\left(\int_{-\infty}^{-S_2}\exp(-s)\hat\xi^2\,ds\right)^\frac{1}{2},
  \end{align}
  where we used the Cauchy-Schwarz inequality in the last step. Inserting (\ref{fo5}) into (\ref{fo6})
  and then into (\ref{fo4}), we obtain
  \begin{align*}
  1\lesssim \exp(-\frac{1}{2}S_2){\rm Nu}_{\BF}^\frac{1}{2}+\exp(3S_2)(\ln H)^{-\frac{1}{5}}.
  \end{align*}
  Finally choosing $S_2\ge 0$ such that $\exp(3S_2)$ is a small multiple of $(\ln H)^\frac{1}{5}$ so
  that the last term can be absorbed into the l.\ h.\ s.\ and to the effect of
  $\exp(-\frac{1}{2}S_2)\sim (\ln H)^{-\frac{1}{30}}$, the above turns into the desired
  $1\lesssim (\ln H)^{-\frac{1}{30}}{\rm Nu}_{\BF}^\frac{1}{2}$.

 \section{Proofs of lemmas}
 In this section we will give the proofs of the four lemmas stated in Subsection \ref{Lemmas}.

  \subsection{Approximate positivity in the bulk: proof of Lemma \ref{Lemma1}}\label{APB}
%
  Much of the effort
  of this construction will consist in designing the kernel $\phi_0$ in such
  a way that it is both non-negative and compactly supported. Non-negativity
  of $\phi_0(s)$ and its fast decay for $s\downarrow-\infty$ and support
  in $s\le 0$ will be crucial in Subsections \ref{Proof-Lemma2} and \ref{Proof-Lemma3}, where we will work
  with the convolution (\ref{convolution}).
  In order to infer non-negativity, we can no longer let $k\uparrow\infty$
  in (\ref{OSCinFourier-reduct}) (as in the proof of Proposition \ref{prop1}), since the two last terms would blow up.
  To quantify this qualitative observation, we restrict to $k>0$
  and introduce the change of variable 
  \begin{equation}\label{k-change}
 z=\frac{\hat z}{k}\,,
  \end{equation}
  so that the stability condition (\ref{OSCinFourier-reduct}) turns into
  \begin{multline}\label{OSC-new}
  \int_0^{kH}\hat\xi\left(\frac{\hat z}{k}\right) w \left(-\frac{d^2}{d\hat z^2}+1\right)^2\bar w\,\frac{d\hat z}{\hat z}\\
  +k^3\int_0^{kH}\left|\frac{d}{d\hat z}\left(-\frac{d^2}{d\hat z^2}+1\right)^2w\right|^2\,d\hat z
  +k^3\int_0^{kH}\left|\left(-\frac{d^2}{d\hat z^2}+1\right)^2w\right|^2\,d\hat z\;\ge\;0\,.
  \end{multline}
  We shall restrict ourselves to  $k$ with $kH\ge 1$
  and real functions $w(\hat z)$ compactly supported
  in $\hat z\in(0,1]$ 
  so that the boundary conditions (\ref{bc})
  are automatically satisfied. 
  In particular, an integration by parts (based on $(-\frac{d^2}{dz^2}+1)^4=\frac{d^8}{dz^8}-4\frac{d^6}{dz^6}+6\frac{d^4}{dz^4}-4\frac{d^2}{dz^2}+1$)
  in the two last terms of (\ref{OSC-new}), yielding
  \begin{align}\nonumber
  &\int_0^{kH}\left|\frac{d}{d\hat z}\left(-\frac{d^2}{d\hat z^2}+1\right)^2w\right|^2\,d\hat z
  +\int_0^{kH}\left|\left(-\frac{d^2}{d\hat z^2}+1\right)^2w\right|^2\,d\hat z\nonumber\\
  &= \int_0^{\infty}\left[\left(\frac{d^5w}{d\hat z^5}\right)^2+5\left(\frac{d^4w}{d\hat z^4}\right)^2
  +10\left(\frac{d^3w}{d\hat z^3}\right)^2+10\left(\frac{d^2w}{d\hat z^2}\right)^2+5\left(\frac{d  w}{d\hat z}\right)^2
  + w^2\right]\,d\hat z\,,\label{correct-terms}
  \end{align}
  shows that there are no fortuitous cancellations: Provided the multiplier
  $\phi:= w (-\frac{d^2}{d\hat z^2}+1)^2 w$
  of $\hat\xi$ is non-negligible in the sense of 
  $\int_0^\infty\phi\,\frac{d\hat z}{\hat z}=O(1)$,
  the two last terms of (\ref{OSC-new}) are at least of $O(k^3)$.
  Hence we are forced to work with $k\ll 1$ and thus, as expected, the conclusion is only effective for
  $z=\frac{\hat z}{k}\gg 1$.

  In anticipation of (\ref{convolution}) we introduce the logarithmic variables
   \begin{equation}\label{her.2}
  s\;=\;\ln\hat z\quad\mbox{and}\quad s'\;=\;-\ln k\,,
  \end{equation}
   so that the first term in the stability condition (\ref{OSC-new}) can
  be rewritten as follows
  \begin{equation*}
  \int_0^{kH}\hat\xi\left(\frac{\hat z}{k}\right) w\left(-\frac{d^2}{d\hat z^2}+1\right)^2\bar w\,
  \frac{d\hat z}{\hat z}=
  \int_{-\infty}^{\infty}
  \hat\xi(s+s') w \left(-\frac{d^2}{d\hat z^2}+1\right)^2 w\,ds\,.
  \end{equation*}
  In view of this and (\ref{correct-terms}), the stability condition (\ref{OSC-new})
  turns into: For all $s'\le\ln H$ we have
  \begin{equation}\label{start}
  \int_{-\infty}^{\infty}\hat\xi(s+s') w\left(-\frac{d^2}{d\hat z^2}+1\right)^2 w\,ds
  \ge-\exp(-3s')\,\int_0^{1}\left[\left(\frac{d^5w}{d\hat z^5}\right)^2+\cdots  +w^2\right]\,d\hat z\,.
  \end{equation}
  Let us consider the l.\ h.\ s.\ of (\ref{start}) in more detail.
  In order to derive a result of the type of (\ref{1}),
  it would be convenient to have a smooth {\it compactly supported}
  $w$ such that the 
  multiplier $\phi=w (-\frac{d^2}{d\hat z^2}+1)^2 w$ is {\it non-negative}.
  Although we don't have an argument, we believe that such
  a $w$ does not exist.
  Instead, we will construct
  \begin{itemize}
  \item[$-$] a family ${\mathfrak F}$
  of smooth functions $w$
  supported in $\hat z\in(0,1]$ 
  \item[$-$] and a probability measure $\rho(dw)$ on ${\mathfrak F}$ which is 
  invariant under the symmetry transformation $w\rightarrow \hat w$ defined through 
  the change of variables $\hat w(\frac 12+z)=w(\frac 12-z)$
  \end{itemize}
  such that the convex combination of the multipliers
  \begin{equation}\label{convex-combin}
  \phi_0(\hat z)\;:=\;\int_{\mathfrak{F}}\phi(\hat z)\,\rho(dw)\,,
  \end{equation}
  where
  $$\phi\;:=\;w\left(-\frac{d^2}{d\hat z^2}+1\right)^2 w\,,$$
  is non-negative, supported in $[\frac 14,\frac 34]$ (and non-trivial) --- and thus satisfies (\ref{convolution-kernel})
  after normalization.
  Roughly speaking, the reason why this can be achieved is the following:
  For any (non-trivial) smooth, compactly supported $w$ we have 
  \begin{itemize}
  \item $\phi=w\left(-\frac{d^2}{d\hat z^2}+1\right)^2 w$ is positive on average:
  \begin{equation*}
  \int_{0}^{1}\phi\,d\hat z\;=\;\int_{0}^{1}\left[\left(\frac{d^2w}{d\hat z^2}\right)^2
  +2\left(\frac{dw}{d\hat z}\right)^2+w^2\right]\,d\hat z\,.
  \end{equation*}
  \item  $\phi=w\frac{d^4w}{dz^4}+\cdots+w^2$ is positive near the edge of the support of $w$
  (incidentally this would {\it not} be true for the positive {\it second} order operator
  $-\frac{d^2}{d\hat z^2}+1$). 
  \end{itemize}
 
  Before becoming much more specific let us address the error term 
  stemming from the r.\ h.\ s.\ of (\ref{start}) for our construction, that is
  \begin{equation}\label{rhs-corr}
  \int_{\mathcal F}\int_0^1\left[\left(\frac{d^5w}{d\hat z^5}\right)^2+\cdots+w^2\right]\,d\hat z\,\rho(dw)\,.
  \end{equation}
  The functions in our family ${\mathfrak F}$ will be of the form
  \begin{equation}\label{w-family}
  w_{\ell,\hat z'}(\hat z)\;:=\;\left(\sqrt{\ell}\right)^{3}w_0\left(\frac{\hat z-\hat z'}{\ell}\right)\,,
  \end{equation}
  that is, translations and rescalings of a ``mask'' $w_0$.
  The mask $w_0$ is some compactly supported smooth function that we fix now, say
  \begin{equation}\label{mask}
  w_0(\hat{\hat z})\;:=\;\left\{\begin{array}{ccc}
  \frac{1}{\sqrt{C_0}}
  \exp\left(-\frac{1}{(1-\hat{\hat z}^2)^2}\right)&\mbox{for}&\hat{\hat z}\in(-1,1)\\
  0&\mbox{else}&\end{array}\right\}\,,
  \end{equation}
  and the normalization constant $C_0$ chosen such that  
  \begin{equation}\label{normalize}
  \int\left(\frac{dw_0}{d\hat{\hat z}^2}\right)^2\, d\hat{\hat z}=1\,.
  \end{equation}
  Provided 	
  \begin{equation}\label{scaling-transl-restr}
  \ell\;\le\;\frac{1}{8}\quad\mbox{and}\quad \hat z'\;\in\;\left(\frac{3}{8},\frac{5}{8}\right)\,,
  \end{equation}
  then $w_{\ell,\hat z'}$ is, as desired, supported in $\hat z\in[\frac 14, \frac 34]$.
  If we choose the length-scale to be bounded away from zero, i.\ e.\
  \begin{equation}\label{scaling-restr2}
  \ell\;\ge\;\frac{1}{C},
  \end{equation}
  then the error term (\ref{rhs-corr}) is clearly finite, so that (\ref{1}) follows from (\ref{start})
  via integration w.\ r.\ t. $\rho(dw)$.

  It thus remains to construct a probability measure $\rho$ in $\ell$ and $\hat z'$ with
  (\ref{scaling-transl-restr}) \& (\ref{scaling-restr2}) such that (\ref{convex-combin}) is non-negative (and non-trivial).
  Note that $w_{\ell,\hat z'}$ in (\ref{w-family}) is scaled such that the
  corresponding multipliers satisfy
  \begin{equation}\label{resc-mult}
  \phi_{\ell,\hat z'}(\hat z)\;=\;\left(\frac{1}{\ell}w_0\frac{d^4}{d\hat{\hat z}^4}w_0
  -2\ell w_0\frac{d^2}{d\hat{\hat z}^2}w_0+\ell^3 w_0^2\right)\left(\frac{\hat z-\hat z'}{\ell}\right)\,,
  \end{equation}
  and $w_0$ is normalized in (\ref{normalize}) in such a way that
  $\int w_0\frac{d^4}{d\hat{\hat z}^4}w_0\,d\hat{\hat z}=1$. Hence
  for all $\hat z'$ 
  we have convergence as $\ell\downarrow 0$ 
  \begin{equation}\label{weak-convergence}
  \phi_{\ell,\hat z'}(\hat z)\;\rightharpoonup\;\delta(\hat z'-\hat z)\quad
  \mbox{when tested against smooth functions of}\;\hat z'\,.
  \end{equation}
  On the other hand we note the following: Writing $w_0=\frac{1}{\sqrt{C_0}}\exp(I)$ with 
  $I=-\frac{1}{(1-\hat{\hat z}^2)^2}$, we have
  \begin{eqnarray}\nonumber
  \frac{d^2w_0}{d\hat{\hat z}^2}&=&
  \frac{1}{\sqrt{C_0}}\left[\left(\frac{dI}{d\hat{\hat z}}\right)^2+\frac{d^2I}{d\hat{\hat z}^2}\right]
  \exp(I),\nonumber\\
  \frac{d^4w_0}{d\hat{\hat z}^4}&=&
  \frac{1}{\sqrt{C_0}}\left[\left(\frac{dI}{d\hat{\hat z}}\right)^4
  +6\left(\frac{dI}{d\hat{\hat z}}\right)^2\frac{d^2I}{d\hat{\hat z}^2}
  +3\left(\frac{d^2I}{d\hat{\hat z}^2}\right)^2
  +4\frac{dI}{d\hat{\hat z}}\frac{d^3I}{d\hat{\hat z}^3}
  + \frac{d^4I}{d\hat{\hat z}^4}
  \right]\exp(I)\,.\nonumber
  \end{eqnarray}
  Since near the edges $\{-1,1\}$ of the support $[-1,1]$ of $w_0$, 
  i.\ e.\ for $1-|\hat{\hat z}|\ll 1$, the term
  $\left(\frac{dI}{d\hat{\hat z}}\right)^4\gg 1$ dominates the other terms
  thanks to the {\it{quadratic}} blow up of $I$ near the edges, we have, according
  to (\ref{resc-mult})
  \begin{equation*}
  \phi_{\ell,\hat z'}(\hat z)
  \;=\;\frac{1}{\ell}w_0\frac{d^4}{d\hat{\hat z}^4}w_0
  -2\ell w_0\frac{d^2}{d\hat{\hat z}^2}w_0+\ell^3 w_0^2\;
  \approx\;\frac{1}{C_0}\frac{1}{\ell}\left[\left(\frac{dI}{d\hat{\hat{z}}}\right)^4\exp(2I)\right]\left(\frac{z-\hat z'}{l}\right)\,.
  \end{equation*}
  Hence in particular for $\ell=\frac{1}{8}$ and $\hat z'=\frac{1}{2}$, $w_{\frac{1}{8},\frac{1}{2}}$
  and  thus $\phi_{\frac{1}{8},\frac{1}{2}}$ are supported in $[\frac{3}{8},\frac{5}{8}]$, 
  $\phi_{\frac{1}{8},\frac{1}{2}}$
  is positive near the edges of the support (and thus bounded away from zero
  at some small distance of the edges of the support), and trivially bounded
  away from $-\infty$ in the support. The universal constants $\delta_0>0$, $\delta_1>0$,
  and $0<C_1<\infty$ are to quantify this:
  \begin{equation}\label{phi-pos-near-edges-supp}
  \phi_{\frac{1}{8},\frac{1}{2}}\;
  \left\{\begin{array}{lcl}
  = 0&\mbox{for}&\hat z\not\in(\frac{3}{8},\frac{5}{8})\\[2mm]
  > 0&\mbox{for}&\hat z\in(\frac{3}{8},\frac{3}{8}+\delta_0]\cup
  [\frac{5}{8}-\delta_0,\frac{5}{8})\\[2mm]
  > \delta_1&\mbox{for}&\hat z\in[\frac{3}{8}+\delta_0,\frac{3}{8}+3\delta_0]
  \cup[\frac{5}{8}-3\delta_0,\frac{5}{8}-\delta_0]\\[2mm]
   > -C_1&\mbox{for}&\hat z\in[\frac{3}{8}+3\delta_0,\frac{5}{8}-3\delta_0]
  \end{array}\right\}\,.
  \end{equation}
  We now choose a universal smooth $\rho_0(\hat z')$ such that
  \begin{equation}\label{rho-0}
  \rho_0\;
  \left\{\begin{array}{lcl}
          =0&\mbox{for}&\hat z'\not\in(\frac{3}{8}+2\delta_0,\frac{5}{8}-2\delta_0)\\[1mm]
          \geq 0 &\mbox{for}&\hat z'    \in[\frac{3}{8}+2\delta_0,\frac{5}{8}-2\delta_0]\\[1mm]
         =2C_1&\mbox{for}&\hat z'    \in[\frac{3}{8}+3\delta_0,\frac{5}{8}-3\delta_0]\\[1mm]
  \end{array}\right\}\,,
  \end{equation}
  and that is even w.\ r.\ t. $\hat z'=\frac{1}{2}$.
  Since $\rho_0$ is smooth in $\hat z'$ we have according to (\ref{weak-convergence}) 
  \begin{equation*}
  \int_{-\infty}^\infty\phi_{\ell,\hat z'}(\hat z)\rho_0(\hat z')\,d\hat z'
  \;\rightarrow\;\rho_0(\hat z)\quad\mbox{uniformly in}\;\hat z\;
  \mbox{as}\;\ell\downarrow0\,.
  \end{equation*}
  In view of the properties (\ref{rho-0}), there thus exists (a possibly small) $\ell_0>0$ such that
  \begin{equation}\label{to-add}
  \int_{-\infty}^\infty\phi_{\ell_0,\hat z'}(\hat z)\rho_0(\hat z')\,d\hat z'\;
  \left\{\begin{array}{lcl}
  = 0 &\mbox{for}&\hat z\not\in(\frac{3}{8}+\delta_0,\frac{5}{8}-\delta_0)\\[2mm]
  \ge -\delta_1&\mbox{for}&\hat z    \in[\frac{3}{8}+\delta_0,\frac{5}{8}-\delta_0]\\[2mm]
  \ge C_1&\mbox{for}&\hat z    \in[\frac{3}{8}+3\delta_0,\frac{5}{8}-3\delta_0]
  \end{array}\right\}\,.
  \end{equation}
  In view of (\ref{phi-pos-near-edges-supp}), the properties (\ref{to-add}) are made to ensure that
  \begin{equation*}
  \phi_0(\hat z)\;:=\;\phi_{\frac{1}{8},\frac{1}{2}}(\hat z)
  +\int_{-\infty}^\infty\phi_{\ell_0,\hat z'}(\hat z)\rho_0(\hat z')\,d\hat z'\;
  \left\{\begin{array}{lcl}
  = 0 &\mbox{for}&\hat z\not\in(\frac{3}{8},\frac{5}{8})\\[1mm]
  >  0 &\mbox{for}&\hat z    \in(\frac{3}{8},\frac{5}{8})
  \end{array}\right\}\,
  \end{equation*}
  defines a $\phi_0$ that is strictly positive in its support and that is
  of the form (\ref{convex-combin}) (after a gratuitous normalization to
  obtain a probability measure).
  
  An inspection of our construction shows that $\phi_0$ satisfies the symmetry property 
  in (\ref{convolution-kernel}).

  \subsection{Approximate logarithmic growth: proof of Lemma \ref{Lemma2}}\label{Proof-Lemma2}
  In this subsection, we return to the approximate logarithmic growth
  of $\tau$ worked out in case of the reduced stability condition in
  Section \ref{heuristics}. Compared to Section \ref{heuristics}, we have to work
  with the mollified version $\hat\xi_0$ of $\hat\xi$, cf.\ (\ref{convolution}) in Lemma \ref{Lemma1},
  since only for the former we have approximate positivity in the
  bulk according to Lemma \ref{Lemma1}. As stated in Lemma \ref{Lemma2}, we shall show that for 
 $S_1\gg 1$
  we have
  \begin{equation}\label{approx-log-growth-OSC}
  \int_{-1}^0\hat\xi_0\,ds\;\lesssim\frac{1}{S_1}\int_{0}^{S_1}\hat\xi_0\,ds+1.
  \end{equation}
  We start the proof recalling 
  \begin{itemize}
 \item  The starting point for Subsection \ref{APB}, that is (\ref{start}), which we rewrite as
  \begin{multline}\label{here-ref}
  \int_{-\infty}^{\infty}
  \hat\xi(s+s'+s'') w\left(-\frac{d^2}{d\hat z^2}+1\right)^2 w\,ds\\ \ge-\exp(-3s'-3s'')\int_0^{1}\left(
  \left(\frac{d^5w}{d\hat z^5}\right)^2+\cdots+w^2\right)\,d\hat z\,,
  \end{multline}
  for all $s'\le\ln H$, $s''\le 0$ and all smooth $w$ compactly supported in $\hat z\in(0,1]$.
  
  \item The outcome of Subsection \ref{APB}, that is (\ref{1}):
  \begin{equation}\label{approx-positiv-bulk-2}
  \hat\xi_0(s')\;=\;\int_{-\infty}^\infty\hat\xi(s'+s'')\phi_0(s'')\,ds''\,
  \;\gtrsim\;-\exp(-3s')
  \end{equation}
  for all $s'\le\ln H$.
  \end{itemize}
  Since the kernel $\phi_0(s'')$ is non-negative and compactly supported in $s''\in (-\infty,0]$,
  we obtain by testing the inequality in $(\ref{here-ref})$ with $\phi_0(s'')\,ds''$:
  \begin{eqnarray}\label{1)+2)}
  \int_{-\infty}^{\infty}\hat\xi_0(s'')\phi(s''-s')\,ds''
  &=&\int_{-\infty}^{\infty}\hat\xi_0(s+s')\phi(s)\,ds\\
  &\gtrsim&-\exp(-3s')\,\int_0^{1}\left[\left(\frac{d^5w}{d\hat z^5}\right)^2+\cdots +w^2\right]\,d\hat z\nonumber\,,
 \end{eqnarray}
%
%
  where we continue to use the abbreviation $\phi$ for the multiplier
  corresponding to the generic $w$:
  $\phi\;=\;w\left(-\frac{d^2}{d\hat z^2}+1\right)^2 w\,.$

  We recall that in terms of $w=\hat z^2\hat w$ and $s=\ln \hat z$, the multiplier assumes the form
  $\phi=\hat w\left(\frac{d^4}{d\hat s^4}+2\frac{d^3}{ds^3}-\frac{d^2}{ds^2}-2\frac{d}{ds}\right)\hat w$, cf.\ (\ref{stability-conditions-in-new-variables}).
  The structure of the argument is similar to the one for (\ref{red-log-growth}) in Section \ref{heuristics}. We seek
  \begin{itemize}
  \item a family ${\mathfrak F}=\{w_{s'}\}_{s'}$ 
  of smooth functions $w_{s'}$ parametrized by $s'\in\mathbb{R}$ 
  and compactly supported in $\hat z\in(0,1]$, that is
  $s=\ln z\in(-\infty,0]$, and
  \item a probability measure $\rho(ds')$ supported in $s'\in(-\infty,\ln H]$,
  \end{itemize}
  such that the corresponding convex combination of multipliers
  shifted by $s'$, i.\ e.\ 
  \begin{equation}\label{phi_1}
  \phi_1(s'')\;:=\;\int_{-\infty}^\infty\phi_{s'}(s''-s')\,\rho(ds'),
  \quad\mbox{where}\quad \phi_{s'}\;:=\;\hat w_{s'}\left(\frac{d^4}{ds^4}+2\frac{d^3}{ds^3}-\frac{d^2}{ds^2}-2\frac{d}{ds}\right) \hat w_{s'}\,,
  \end{equation}
%
  is estimated by above as follows
  \begin{eqnarray}\label{upper-bound-OSC}
  \phi_1(s'')&\le&\left\{\begin{array}{ccc}
  -1&\mbox{for}&-1\le s''\le 0\\
  \frac{C_1}{S_1}&\mbox{for}&0\le s''\le S_1\\
  0&\mbox{else}&\end{array}\right\}\,,
  \end{eqnarray}
  where $C_1(\le S_1)$ is some universal constant, the value of which we want to
  remember momentarily, and estimated by below 
  \begin{eqnarray}\label{weak-lower-bound}
  \phi_1(s'')&\gtrsim&-\left\{\begin{array}{ccc}
  \exp(\frac{s''}{S_2})&\mbox{for}&s''\le 0\\
  1&\mbox{for}&s''\ge 0
  \end{array}\right\}\,,
  \end{eqnarray}
  where the exponential rate $\frac{1}{S_2}\gg 1$ could be replaced by any rate 
  larger than 3.
  Furthermore, we need that
  \begin{equation}\label{error-control}
  \int_{-\infty}^\infty\exp(-3s')\int_0^{1}\left[\left(\frac{d^5w_{s'}}{d\hat z^5}\right)^2+\cdots+w_{s'}^2\right]\,d\hat z\,\rho(ds')\;\lesssim\;1\,.
  \end{equation}
%
  
  For later use we note that in terms of $w_{s'}=\hat z^2\hat w_{s'}$, (\ref{error-control}) turns into 
  $$\int_{-\infty}^{\infty}\exp(-3s')\int_0^1\left[\hat z^4\left(\frac{d^5 \hat w_{s'}}{d\hat z^5}\right)^2+\hat z^2\left(\frac{d^4 \hat w_{s'}}{d\hat z^4}\right)^2+\left(\frac{d^3 \hat w_{s'}}{d\hat z^3}\right)^2+\cdots+\hat w_{s'}^2\right]\, d\hat z \rho(ds')\lesssim 1\,.$$
  In terms of $s=\ln \hat z$, this means 
  \begin{equation}\label{error-control-bis}
   \int_{-\infty}^{\infty}\exp(-3s')\int_0^1\exp(-5s)\left[\left(\frac{d^5\hat w_{s'}}{d\hat z^5}\right)^2+\cdots+\hat w_{s'}^2\right]\, ds \rho(ds')\lesssim 1\,.
  \end{equation}  
  It is almost obvious how (\ref{upper-bound-OSC}), (\ref{weak-lower-bound}) \& (\ref{error-control}) 
  allow to pass
  from (\ref{1)+2)}) to (\ref{approx-log-growth-OSC}) by substituting $w$ with $w_{s'}$  and
  integrating in $\rho(ds')$. We just need to show how (\ref{upper-bound-OSC}) 
  \& (\ref{weak-lower-bound}) yield
  \begin{equation*}
  \int_{-\infty}^\infty\hat\xi_0\phi_1\,ds''
  \;\le\;
  -\int_{-1}^0\hat\xi_0\,ds''
  +\frac{C_1}{S_1}\int_{0}^{S_1}\hat\xi_0\,ds''
  +C\,.
  \end{equation*}
  Indeed, we write 
  \begin{eqnarray}\nonumber
  &&\int_{-\infty}^\infty\hat\xi_0\phi_1\,ds''
  +\int_{-1}^0\hat\xi_0\,ds''
  -\frac{C_1}{S_1}\int_{0}^{S_1}\hat\xi_0\,ds''\nonumber\\
  &&=
  \int_{-\infty}^\infty(-\hat\xi_0)
  \left(\left\{\begin{array}{ccc}
  -1&\mbox{for}&-1\le s''\le 0\\
  \frac{C_1}{S_1}&\mbox{for}&0\le s''\le S_1\nonumber\\
  0&\mbox{else}&\end{array}\right\}-\phi_1\right)\,ds''\\
  &\stackrel{(\ref{upper-bound-OSC}),(\ref{approx-positiv-bulk-2})}{\lesssim}&
  \int_{-\infty}^\infty\exp(-3s'')
  \left(\left\{\begin{array}{ccc}
  -1&\mbox{for}&-1\le s''\le 0\\
  \frac{C_1}{S_1}&\mbox{for}&0\le s''\le S_1\nonumber\\
  0&\mbox{else}&\end{array}\right\}-\phi_1\right)\,ds''\\
  &\stackrel{(\ref{weak-lower-bound})}{\lesssim}&
  \int_{-\infty}^\infty\exp(-3s'')
  \left\{\begin{array}{ccc}
  \exp(\frac{s''}{S_2})&\mbox{for}&s''\le 0\\
         1&\mbox{for}&s''\ge 0
  \end{array}\right\}
  \,ds''\;\stackrel{S_2<\frac 13}{\lesssim}\;1\,.\nonumber
  \end{eqnarray}

  The proof of (\ref{upper-bound-OSC}) is the same as the one for (\ref{phi1-upper-bound}) in Section \ref{heuristics} until the
  point where we used an approximation argument to allow for the non-smooth choice (\ref{w_0-near-support}) of $\hat w_{s'}(s)$.
  This approximation argument in $H^{2,\infty}$ was sufficient for $\phi_1$, which in view of the representation (\ref{phi_1-formula})
  is continuous with respect to $H^{2,2}$. In our situation this approximation argument is not sufficient because we need to control
  the error term (\ref{error-control-bis}) which requires boundedness in $H^{5,2}$.
  We now summarize the main steps of the proof of (\ref{upper-bound-OSC}), (\ref{weak-lower-bound}) and 
  (\ref{error-control-bis}) leaving out the details that can be found in Section \ref{heuristics}.
 
 \noindent
 \underline{Construction of the family $\{w_{s'}\}_{s'}$}:
  \begin{itemize}
  
     \item As in Section $\ref{heuristics}$, we fix a smooth mask $\hat w_0$ such that
    \begin{equation}\label{B.7}
    \hat w_0\quad\mbox{is supported in}\;\hat s\in\left[-\frac{1}{2},0\right]
    \quad\mbox{and nonvanishing in}\;\left(-\frac{1}{2},0\right)\,,
    \end{equation}
      and normalize it by $\int\hat w_0^2\, d\hat s=1$.
     \item As in Section \ref{heuristics}, we introduce the change of variables 
     \begin{equation}\label{B.1}
     s=\lambda \hat s  \quad \mbox{ and } \quad \hat w_{\lambda}\;=\;\lambda^{-\frac 12}\hat w_0\,,
     \end{equation}
     and rewrite the corresponding multiplier as
       \begin{equation}\label{B.6}
      \phi_{\lambda}\;=\;-\frac{2}{\lambda^2}\hat w_0\frac{d\hat w_0}{d\hat s  }
              -\frac{1}{\lambda^3}\hat w_0\frac{d^2\hat w_0}{d\hat s^2}
              +\frac{2}{\lambda^4}\hat w_0\frac{d^3\hat w_0}{d\hat s^3}
              +\frac{1}{\lambda^5}\hat w_0\frac{d^4\hat w_0}{d\hat s^4}\,,
	\end{equation}
       a form that highlights the desired dominance of the term 
      $-\frac{2}{\lambda^2}\hat w_0\frac{d\hat w_0}{d\hat s  }
      =-\frac{1}{\lambda^2}\frac{d\hat w_0^2}{d\hat s}$ for $\lambda\gg 1$ and 
      (heuristically) suggests the shape of the probability measure $\rho(s')$ as increasing over scales of order $1$ and 
      eventually decreasing over scales of order $S_1$ (see (\ref{form})).
       We note for later reference that in these variables, (\ref{error-control-bis})
       assumes the form 
       $$\int_{-\infty}^{\infty}\exp(-3s')\int_{-\infty}^0\exp(-5\lambda \hat s)
       \left[\frac{1}{\lambda^{10}}\left(\frac{d^5 \hat w_0}{d\hat s^5}\right)^2+\cdots +\hat w_0^2\right]d\hat s\rho(ds')\lesssim 1\,,$$
       and because of (\ref{B.7}) for $\lambda\geq 1$ follows from 
       \begin{equation}\label{error-control-bisbis}
       \int_{-\infty}^{\infty}\exp\left(-3s'+\frac 52\lambda\right)\int_{-\infty}^{\infty}\left(\left(\frac{d^5 \hat w_0}{d\hat s^5}\right)^2+
       \cdots +\hat w_0^2\right)d\hat s\rho(ds')\lesssim 1\,.
       \end{equation}

     \item In oder to obtain $\phi_1(s'')\sim -1$ over an $s''$-interval of
     length of the order $1$ followed by $\phi_1(s'')\lesssim\frac{1}{S_1}$, 
     we choose as in Section \ref{heuristics}
     \begin{equation}\label{B.4}
     \lambda=s', 
     \end{equation}
     meaning that $\lambda$ is small in the foot regions and 
     large in the plateau region (see argument after (\ref{form}), leading to this choice).
     Eventually we will need to modify (\ref{B.4}) for moderate and small $s'$, cf.\ (\ref{B.32}).
      With the choice of (\ref{B.4}),  (\ref{B.6}) turns into
     \begin{eqnarray}\label{B.9}
     \phi_{s'}(s)
     &=&-\frac{1}{(s')^2}\left(\frac{d\hat w_0^2}{d\hat s  }\right)\left(\frac{s}{s'}\right)
     -\frac{1}{(s')^3}\left(\hat w_0\frac{d^2\hat w_0}{d\hat s^2}\right)\left(\frac{s}{s'}\right)\nonumber\\
     &+&\frac{2}{(s')^4}\left(\hat w_0\frac{d^3\hat w_0}{d\hat s^3}\right)\left(\frac{s}{s'}\right)
     +\frac{1}{(s')^5}\left(\hat w_0\frac{d^4\hat w_0}{d\hat s^4}\right)\left(\frac{s}{s'}\right)\,.
     \end{eqnarray}
     
     \item As in Section \ref{heuristics}, replacing in (\ref{phi_1}) the integration over $s'$ by the integration
      over the argument $\hat s=\frac{s}{\lambda}$ of $\hat w_0$ according to the 
     nonlinear change of variable
     \begin{equation}\label{B.8}
      \hat s\stackrel{(\ref{B.1})}{=}\frac{s''-s'}{\lambda}=\frac{s''-s'}{s'}=\frac{s''}{s'}-1\;\quad
      \Longleftrightarrow\quad s'\;=\;\frac{s''}{1+\hat s}\,,
    \end{equation}
    $\phi_1$ can be written as 
     \begin{eqnarray}\label{B.first-rep}
     \phi_1(s'')&=&
     -\int_{-\infty}^\infty\frac{1}{(1+\hat s)^2}\,\hat w_0^2
       \frac{d\rho}{ds'}\,d\hat s-\frac{1}{(s'')^2}\int_{-\infty}^\infty(1+\hat s)\,\hat w_0\frac{d^2\hat w_0}{d\hat s^2}
       \rho\,d\hat s\\
       &+&\frac{2}{(s'')^3}\int_{-\infty}^\infty(1+\hat s)^2\,\hat w_0\frac{d^3\hat w_0}{d\hat s^3  }
       \rho\,d\hat s+\frac{1}{(s'')^4}\int_{-\infty}^\infty(1+\hat s)^3\,\hat w_0\frac{d^4\hat w_0}{d\hat s^4  }
     \rho\,d\hat s\,.\nonumber
     \end{eqnarray}
   \end{itemize}
   \underline{Construction of the probability measure $\rho(s')\,ds'$ supported in $s'\in(-\infty,\ln H]:$}
   \begin{itemize}
   \item As in Section \ref{heuristics}, we formulate an intermediate goal: Find a measure $0\leq \trho(s')\leq 1$ but not supported in $(-\infty,\ln H)$ such that 
   $\widetilde{\phi}_1(s'')=\int\phi_{s'}(s''-s')\trho(s')\, ds'$ satisfies
      \begin{equation}\label{B.25}
    -\widetilde{\phi}_1(s'')\;\in\;
    \left\{\begin{array}{lcccccc}
    [\frac{1}{C}\frac{1}{(s'')^2},C\frac{1}{(s'')^2}]&\mbox{for}&S_0           &\le&s''&   &\\
    (0,C]                                            &\mbox{for}&\frac{1}{2}S_0&<  &s''&\le&S_0\\
    \{0\}                                            &\mbox{for}&              &   &s''&\le&\frac{1}{2}S_0
    \end{array}\right\}\,.
    \end{equation}
    \item From the representation (\ref{B.first-rep}) and the assumption that $\rho$ varies slowly,
     we learn that $\widetilde{\phi}_1$ is negative if $\frac{d\trho}{ds'}\gg \frac{1}{(s')^2}$.
    In Section \ref{heuristics} this motivated the following Ansatz for $\trho$ in the range $1\ll s'\ll S_1$: 
    we fix a smooth mask $\trho_0(\hat s')$ such 
   \begin{equation}\label{B.17}
    \trho_0\;=\;0\;\;\mbox{for}\;\hat s'\le 0,\quad
     \frac{d\trho_0}{d\hat s'}\;>\;0\;\;\mbox{for}\;0<\hat s'\le 2,\quad
     \trho_0\;=\;1-\frac{1}{\hat s'}\;\;\mbox{for}\;2\le \hat s' \,.
     \end{equation}
     For $S_0\gg 1$, consider the rescaled version
     \begin{equation}\label{B.24}
     \trho(S_0(\hat s'+1))\;=\;\trho_0(\hat s')\,.
     \end{equation}
     Eventually, for (\ref{error-control-bisbis}) and departing from the argument in Section \ref{heuristics},
     we will have to modify $\trho$ for moderate and small $s'$, cf.\ (\ref{B.28}).
     \item Finally, this $\trho$ does not decrease to zero on the large scales $s'\sim S_1$, which has to be done 
     by cutting it off as in Section \ref{heuristics}, cf.\ (\ref{reduced-goal}). This allows to pass from the intermediate 
     goal (\ref{B.25}) to its final version (\ref{upper-bound-OSC}). 
      \end{itemize}

   Exactly as in the proof of (\ref{reduced-goal}), we distinguish the regions of small, intermediate and large $s''$
    (note that for $s''\in (-\infty,\frac{S_0}{2}]$ all the integrals in (\ref{B.first-rep}) vanish
    because the supports of $\hat w_0$ and $\rho$ do not intersect).
  In Section \ref{heuristics}  
    we established
   \begin{equation}\label{B.26}
   \widetilde{\phi}_1\;\sim\;
   -\frac{1}{S_0}\frac{1}{(s'')^2}
   \quad\mbox{uniformly in}\;s''\ge 3S_0\quad\mbox{for}\;S_0\gg 1\,
   \end{equation}
  and
   \begin{equation}\label{B.27}
   \widetilde{\phi}_1\;\sim\;
   -\frac{1}{S_0}
   \quad\mbox{uniformly in}\; s''\in\left[\frac{3}{4}S_0,3S_0\right]\quad\mbox{for}\;S_0\gg 1\,.
   \end{equation}

    As we have seen in Section \ref{heuristics}, in the  \underline{range of small $s''$}, i.\ e.\
    $s''\;\in\;\left(\frac{1}{2}S_0,\frac{3}{4}S_0\right)$,
  the behavior of $\phi_{s'}$ near the left edge $-\frac{1}{2}$ of
   is dominated by the $\frac{1}{\lambda^5}\hat w_0\frac{d^4\hat w_0}{d\hat s^4}$-term
    and thus automatically is {\it strictly positive}. 
   In Section \ref{heuristics}, we 
   solved this problem by  {\it giving up smoothness} of $\hat w_0$ near the
   left edge $-\frac{1}{2}$ of its support $[-\frac{1}{2},0]$ and eventually using an approximation 
   argument in $H^{2,2}$. As discussed earlier, we cannot use this approximation in the present situation,
   since we need to keep the error term (\ref{error-control-bisbis}) under control.
   The remainder of this section is devoted to the way out to this dilemma and it consists of three steps (the first one 
   is the same as in (\ref{small-s''}) and we report it just for the sake of clarity).
  \begin{itemize}
  \item In the first stage, we {\it give up smoothness} of $\hat w_0$ near the
  left edge $-\frac{1}{2}$ of its support $[-\frac{1}{2},0]$. In fact, as in Section \ref{heuristics}, we shall first
  assume that $\hat w_0$ is of the specific form
  \begin{equation}\label{B.22}
  \hat w_0\;=\;\frac{1}{2}\left(\hat s+\frac{1}{2}\right)^2\quad\mbox{for}\;\hat s\in\left[-\frac{1}{2},-\frac{1}{4}\right]\,.
  \end{equation}
  This means that $\hat w_0$ has a bounded but discontinuous second derivative. 
  Our non-smooth Ansatz together with  (\ref{B.6}) implies
   \begin{equation}\label{B.23}
 \phi_{s'}\;=\;-\frac{1}{(s')^2}\left(\frac{s}{s'}+\frac{1}{2}\right)^3
               -\frac{1}{2}\frac{1}{(s')^3}\left(\frac{s}{s'}+\frac{1}{2}\right)^2\;<\;0
 \quad\mbox{for}\;s\in\left(-\frac{s'}{2},-\frac{s'}{4}\right]\,.
 \end{equation}
  As shown in Section \ref{heuristics}, the corresponding $\widetilde{\phi}_1$ is, as desired, strictly 
  negative on $s''\in(\frac{S_0}{2},\frac{3}{4}S_0]$ for all $S_0$ sufficiently large.
  We fix a sufficiently large but universal $S_0$ such that together with (\ref{B.26})\& (\ref{B.27}) we obtain 
  \begin{equation*}
  -\widetilde{\phi}_1(s'')\;\in\;
  \left\{\begin{array}{lcccccc}
  [\frac{1}{C}\frac{1}{(s'')^2},C\frac{1}{(s'')^2}]&\mbox{for}&S_0           &\le&s''&   &\\
  (0,C]                                            &\mbox{for}&\frac{1}{2}S_0&<  &s''&\le&S_0\\
  \{0\}                                            &\mbox{for}&              &   &s''&\le&\frac{1}{2}S_0
  \end{array}\right\},
  \end{equation*}
  for some generic universal constant $C$.
  \item In the second stage, \textit{we modify the definition (\ref{B.24}) of $\trho(s')$} by adding a small-amplitude
  and fast-decaying (exponential) tail for $s'\downarrow-\infty$. More precisely,
  we make the Ansatz
  \begin{equation}\label{B.28}
  \ttrho\;=\;\trho+\varepsilon\delta\trho\quad\mbox{with}\quad
  \delta\trho\;:=\;\exp\left(\frac{s'}{S_2}\right)\eta_0\left(\frac{s'}{S_0}\right)\,,
  \end{equation}	
  where $\eta_0(\hat s')$ is the mask of a smooth cut-off function with
  \begin{equation}\label{B.31}
  \eta_0\;=\;1\quad\mbox{for}\;\hat s'\le 2\quad\mbox{and}\quad
  \eta_0\;=\;0\quad\mbox{for}\;\hat s'\ge 3\,.
  \end{equation}
  Here $0<S_2\ll1$ is some small length-scale and $\varepsilon\ll 1$ is some small amplitude
  to be chosen below. Recall that $S_0$ is the universal constant fixed in the first stage. Since $\ttrho$
  is no longer supported on $s'\in[S_0,\infty)$ but is positive on the entire line,
  we need to extend our definition of the function $\hat w_{s'}$ from $s'\ge S_0$ to all $s'$. In
  view of (\ref{B.1}) we just have to extend the definition (\ref{B.4}) of
  the rescaling parameter $\lambda(s')$ to
  \begin{equation}\label{B.32}
  \lambda\;=\;\left\{\begin{array}{ccc}
  s' &\mbox{for}&s'\ge S_0\\
  S_0&\mbox{for}&s'\le S_0
  \end{array}\right\}\,.
  \end{equation}
  We will show that we can first choose a universal $0<S_2\ll 1$ and then a universal
  $0<\varepsilon\ll 1$ such that we obtain for 
  $\pp(s''):=\int_{-\infty}^\infty\phi_{\lambda(s')}(s''-s')\ttrho(s')\,ds'$
  the following estimates
  \begin{equation}\label{B.29}
  -\pp(s'')\;\in\;
  \left\{\begin{array}{lcccccc}
  \,[ \frac{1}{C}\frac{1}{(s'')^2},C\frac{1}{(s'')^2}]&\mbox{for}&S_0           &\le&s''&   &\\
  \,[ \frac{1}{C},C]                                  &\mbox{for}&\frac{1}{2}S_0&<  &s''&\le&S_0\\
  \,[ \frac{1}{C}\exp(\frac{s''}{S_2}),C\exp(\frac{s''}{S_2})]
                                                  &\mbox{for}&              &   &s''&\le&\frac{1}{2}S_0
  \end{array}\right\},
  \end{equation}
for some generic universal constant $C$. The gain with respect to $\widetilde{\phi}_1$ is that $\pp$
is strictly negative also for $s''\leq \frac 12 S_0\,$ which will allow us to pass to the third stage.
\item In a third stage, \textit{we smoothen out $\hat w_0$}: We define a sequence of smooth functions 
$\{\tilde{w}_{0}^{\alpha}\}_{\alpha\downarrow 0}$ which approximate $\hat w_0$
in such a way that the corresponding $\ppa$ still satisfies (\ref{B.29}). This takes care of (\ref{error-control-bisbis}):
Since for fixed $\alpha>0$ to be chosen later, (\ref{error-control-bisbis}) with $\hat w_0$ replaced by $\hat w_0^{\alpha}$
reduces to
 $$\int_{-\infty}^{\infty}\exp\left(-3s'+\frac 52\lambda(s')\right)\ttrho(s')\, ds'\lesssim 1\,.$$
 For $s'\leq S_0$ this follows from $\ttrho(s')\stackrel{(\ref{B.24})\& (\ref{B.28})}{=}\varepsilon \exp\left(\frac{s'}{S_2}\right)$
 and $\lambda(s')\stackrel{(\ref{B.32})}{=}S_0$ because of $S_2\ll 1$. For $s'\geq S_0$, this follows from $\ttrho\lesssim 1 $
 and $\lambda(s')\stackrel{(\ref{B.32})}{=}s'$.
\end{itemize}
%
%

  We turn to the details for the \underline{second stage}, i.\ e.\ the effect of the modification $\ttrho(s')$ of $\trho(s')$. Consider the perturbation
  $\delta\widetilde{\phi}_1$ of the multiplier $\widetilde{\phi}_1$:
  \begin{equation}\label{B*}
  \delta\widetilde{\phi}_1(s'')\;:=\;\int_{-\infty}^\infty\phi_{\lambda(s')}(s''-s')\delta\trho(s')\,ds'
  \;\stackrel{(\ref{B.32})}{=}\;\int_{-\infty}^\infty\phi_{\lambda(s')}(s''-s')\exp\left(\frac{s'}{S_2}\right)\eta_0\left(\frac{s'}{S_0}\right)\,ds'\,.
  \end{equation}
  In order to show that the unperturbed (\ref{B.25}) upgrades to (\ref{B.29}), it is enough to establish
  \begin{equation}\label{B.30}
  -\delta\widetilde{\phi}_1(s'')\;\in\;
  \left\{\begin{array}{lcccccc}
  \{0\}                                       
   &\mbox{for}&3S_0          &\le&s''&   &\\
  \,[-C,C]                                       
   &\mbox{for}& S_0          &\le&s''&\le&3S_0\\
  \,[\frac{1}{C},C]
   &\mbox{for}&\frac{1}{2}S_0&\le&s''&\le& S_0\\
  \,[\frac{1}{C}\exp(\frac{s''}{S_2}),C\exp(\frac{s''}{S_2})]
   &\mbox{for}&              &   &s''&\le&\frac{1}{2}S_0
  \end{array}\right\},
  \end{equation}
  for some sufficiently small {\it but fixed} $S_2$, where $C$ denotes a universal constant.
  Indeed, choosing $\varepsilon\ll 1$, we see from $\pp=\widetilde{\phi}_1+\varepsilon\delta\widetilde{\phi}_1$
  that (\ref{B.30}) upgrades (\ref{B.25}) to (\ref{B.29}). 

  We start the argument of (\ref{B.30}) with the range of large $s''$, i.\ e.\ $s''\ge 3S_0$, and consider the
  integral $\delta\widetilde{\phi}_1(s'')=\int_{-\infty}^\infty\phi_{\lambda(s')}(s''-s')
  \exp(-\frac{s'}{S_2})\eta_0(\frac{s'}{S_0})\,ds'$. Because of our choice
  (\ref{B.31}) of the cut-off $\eta_0$, the second factor $\exp(-\frac{s'}{S_2})\eta_0(\frac{s'}{S_0})$
  is supported in $s'\in(-\infty,3S_0]$. We note that in view of our choice (\ref{B.32}) of the
  scaling factor $\lambda$, $\hat w_{s'}(s)$ and thus $\phi_{\lambda(s')}(s)$ are supported in
  $s\in[-\frac{1}{2}S_0,0]$ for $s'\le S_0$ and in $s\in[-\frac{1}{2}s',0]$ for $s'\ge S_0$. Hence
  $(s',s'')\mapsto\phi_{\lambda(s')}(s''-s')$ is supported in 
  $s''\in[s'-\frac{1}{2}S_0,s']$ for $s'\le S_0$ and in $s''\in[\frac{1}{2}s',s']$ for $s'\ge S_0$,
  or --- equivalently --- in $s'\in[s'',s''+\frac{1}{2}S_0]$ for $s''\le\frac{S_0}{2}$ and
  in $s'\in[s'',2s'']$ for $s''\ge\frac{S_0}{2}$. Since $s''\ge 3S_0$, we are in the
  latter case and $\phi_{\lambda(s')}(s''-s')$ is supported in $s'\in[s'',2s'']\subset[3S_0,\infty)$. 
  Hence both factors $\exp(-\frac{s'}{S_2})\eta_0(\frac{s'}{S_0})$ and $\phi_{\lambda(s')}(s''-s')=\phi_{s'}(s''-s')$
  have disjoint support in $s'$ and thus the integral (\ref{B*}) in $s'$ vanishes. This establishes the first line in (\ref{B.30}).
  
  We now turn to the very small $s''$, i.\ e.\ $s''\le\frac{S_0}{2}$ in (\ref{B.30}). By the above,
  $s'\mapsto\phi_{\lambda(s')}(s''-s')$ is supported in $s'\in [s'',s''+\frac{S_0}{2}]\subset(-\infty,S_0]$;
   in this $s'$-range we have  for the cut-off function $\eta_0(\frac{s'}{S_0})=1$, and  $\phi_{\lambda(s')}=\phi_{S_0}$.
    Hence the definition (\ref{B*}) simplifies to
  \begin{equation}\label{B.33}
  \frac{1}{\varepsilon}\delta\widetilde{\phi}_1(s'')=
  \int_{-\infty}^\infty\phi_{S_0}(s''-s')\exp\left(\frac{s'}{S_2}\right)\,ds'
  =\exp\left(\frac{s''}{S_2}\right)\int_{-\infty}^\infty\phi_{S_0}(s)\exp\left(-\frac{s}{S_2}\right)\,ds\,.
  \end{equation}
   We note that by (\ref{B.23}) we have
  \begin{equation}\nonumber
  \phi_{S_0}<0 \mbox{ for }\;s\in\left(-\frac{S_0}{2},-\frac{S_0}{4}\right)\quad
  \mbox{and supported in}\;s\in\left[-\frac{S_0}{2},0\right]\,.
  \end{equation}
  This allows us to use Laplace's method for $S_2\ll 1$ in the integral in (\ref{B.33}),
  \begin{eqnarray}\nonumber
  \lefteqn{-\int_{-\infty}^\infty\phi_{S_0}(s)\exp\left(-\frac{s}{S_2}\right)\,ds}\nonumber\\
  &\approx&
  -\int_{-\infty}^{-\frac{S_0}{4}}\phi_{S_0}(s)\exp\left(-\frac{s}{S_2}\right)\,ds\nonumber\\
  &\stackrel{(\ref{B.23})}{=}&
  \int_{-\frac{S_0}{2}}^{-\frac{S_0}{4}}\left(\frac{1}{ S_0^2}\left(\frac{s}{S_0}+\frac{1}{2}\right)^3
  +\frac{1}{2S_0^3}\left(\frac{s}{S_0}+\frac{1}{2}\right)^2\right)\exp\left(-\frac{s}{S_2}\right)\,ds\nonumber\\
  &\approx&
  \int_{-\frac{S_0}{2}}^{\infty}\frac{1}{2S_0^3}\left(\frac{s}{S_0}+\frac{1}{2}\right)^2
  \exp\left(-\frac{s}{S_2}\right)\,ds\nonumber\\
  &=&\frac{1}{S_0^2}\int_{-\frac{1}{2}}^\infty
  \frac{1}{2}\left(\hat s+\frac{1}{2}\right)^2\exp\left(-\frac{S_0}{S_2}\hat s\right)\,d\hat s\nonumber\\
  &=&\frac{S_2^3}{S_0^5}\exp\left(\frac{1}{2}\frac{S_0}{S_2}\right)\,.\nonumber
  \end{eqnarray}
  Plugging this into (\ref{B.33}) yields as claimed in (\ref{B.30})
  \begin{equation}\label{B.34}
  -\delta\widetilde{\phi}_1(s'')\;
  \approx\;\frac{S_2^3}{S_0^5}\exp\left(\frac{1}{2}\frac{S_0}{S_2}\right)\exp\left(\frac{s''}{S_2}\right)
  \quad\mbox{uniformly in}\;s''\le\frac{S_0}{2}\;\mbox{for}\;S_2\ll 1\,.
  \end{equation}
  We now treat the intermediary small values $\frac{S_0}{2}\le s''\le S_0$ in (\ref{B.30}). 
  This time, the function
  $s'\mapsto\phi_{\lambda(s')}(s''-s')$ is supported in $s'\in[s'',2s'']\subset(-\infty,2S_0]$,
  so that also in this $s'$-range we have for the cut-off function $\eta_0(\frac{s'}{S_0})=1$. Hence the
  representation simplifies to
  \begin{eqnarray}\nonumber
  \delta\widetilde{\phi}_1(s'')&=&
  \int_{-\infty}^\infty\phi_{\lambda(s')}(s''-s')\exp\left(\frac{s'}{S_2}\right)\,ds'\,.\nonumber
  \end{eqnarray}
  On this integral, we can again use the Laplace's method for $S_2\ll 1$: By (\ref{B.23}) we have
  for the continuous function $(s',s'')\mapsto\phi_{\lambda(s')}(s''-s')$
  \begin{equation*}
  \phi_{\lambda(s')}(s''-s')\;\left\{\begin{array}{ccc}
  <&0&\mbox{for}\;s'\in(\frac{3}{2}s'',2s'')\\
  =&0&\mbox{for}\;s'\not\in(s'',2s'')
  \end{array}\right\}\,.
  \end{equation*}
  Hence we obtain as claimed in (\ref{B.30})
  \begin{equation}\label{B.35}	
  \delta\widetilde{\phi}_1(s'')\;<\;0\quad
  \mbox{uniformly in}\;s''\in\left[\frac{S_0}{2},S_0\right]\;\mbox{for}\;S_2\ll 1\,.
  \end{equation}
  We finally address the remaining intermediary range, that is, $S_0\le s''\le 3S_0$.
  We clearly have by continuity of $(s',s'')\mapsto\phi_{\lambda(s')}(s''-s')$ and $\eta_0(\hat s')$ that
  \begin{equation}\label{B.36}
  \delta\widetilde{\phi}_1(s'')=
  \int_{-\infty}^\infty\phi_{\lambda(s')}(s''-s)\exp\left(\frac{s'}{S_2}\right)\eta_0\left(\frac{s'}{S_0}\right)\,ds'
  \end{equation}
  is uniformly bounded for $s''\in[S_0,3S_0]$.
  Estimate (\ref{B.30}) now follows from (\ref{B.34}), (\ref{B.35}) \& (\ref{B.36}) 
  for a choice of sufficiently small $S_2$.
  
  We now turn to the details for the \underline{third stage}. We approximate $\hat w_0$, which is non-smooth
  at the left edge of its support, cf.\ (\ref{B.22}), by a sequence of smooth $\hat w^{\alpha}_0$
  in such a way that the corresponding $\phi_{\lambda(s')}$ and $\phi^{\alpha}_{\lambda(s')}$ 
  are close in $L^1$. More precisely, we select a smooth function $F(w)$ with
  \begin{equation*}
  F(w)\;=\;0\quad\mbox{for}\;w\le 0\quad\mbox{and}\quad
  F(w)\;=\;w\quad\mbox{for}\;w\ge 1\,.
  \end{equation*}
  For a small parameter $0<\alpha\ll 1$ we now define $\hat w_0^{\alpha}(\hat s)$ via
  \begin{equation*}
  \hat w_0^{\alpha}\;:=\;\alpha^2 F\left(\frac{\hat w_0}{\alpha^2}\right)
  \;\stackrel{(\ref{B.22})}{=}\;\alpha^2 F\left(\frac{(\hat s+\frac{1}{2})^2}{2\alpha^2}\right)
  \quad\mbox{for}\;\;\hat s\in\left[-\frac{1}{2},-\frac{1}{4}\right]\,;
  \end{equation*}
  for $\hat s\not\in[-\frac{1}{2},-\frac{1}{4}]$, $\hat w_0^{\alpha}$ is set equal to $\hat w_0$.
  Clearly, the so defined $\hat w_0^{\alpha}$ is smooth on the whole line.
%
%
  We want to show that the convex combination of multipliers
  \begin{equation*}
  \ppa(s'')=\int_{-\infty}^\infty\phi^{\alpha}_{\lambda(s')}(s''-s')\ttrho(s')\,ds'\,,
  \end{equation*}
  still satisfies (\ref{B.29}), that is
  \begin{equation}\label{B.37}
  -\ppa(s'')\;\in\;
  \left\{\begin{array}{lcccccc}
  \,[ \frac{1}{C}\frac{1}{(s'')^2},C\frac{1}{(s'')^2}]&\mbox{for}&           S_0&\le&s''&   &\\
  \,[ \frac{1}{C},C]                                  &\mbox{for}&\frac{1}{2}S_0&<  &s''&\le&S_0\\
  \,[ \frac{1}{C}\exp(\frac{s''}{S_2}),C\exp(\frac{s''}{S_2})]
						    &\mbox{for}&              &   &s''&\le&\frac{1}{2}S_0
  \end{array}\right\}\,,
  \end{equation}
  for some choice of $0<\alpha\ll 1$ and a generic universal constant $C$.
  For this purpose we consider the difference of the combination of multipliers, that is,
  $\dppa:=\ppa-\pp$. In order to pass from (\ref{B.29}) to (\ref{B.37}), it is sufficient to establish
  \begin{equation}\label{B-new}
  |\dppa(s'')|\lesssim \alpha\left\{\begin{array}{lcccccc}
                                                       \frac{1}{(s'')^2} & \mbox{ for }& 3S_0&\leq& s''\,,\\
                                                       1             & \mbox{ for }& \frac 12 S_0&\leq& s''\leq 3S_0\,,\\
                                                       \exp\left(\frac{s''}{S_2}\right) & \mbox{ for }& s''&\leq& \frac{1}{2}S_0
                                                      \end{array}\right\}\,.
  \end{equation}
  To this aim, we first observe that
%
%
  \begin{equation}\label{B.40}
   \left|\hat w_0^{\alpha}\frac{d^k\hat w_0^{\alpha}}{d\hat s^k}-\hat   w_0 \frac{d^k\hat   w_0}{d\hat s^k}\right|\lesssim \alpha^{4-k} \qquad \mbox{ with } k=0,\cdots, 4\,.
  \end{equation}

  which follows from the fact that
  \begin{equation}\label{extraline}
    \mbox{ all these differences are supported on the interval }\qquad  \hat s\in\left[-\frac{1}{2},-\frac{1}{2}+\sqrt{2}\alpha\right]\,,
  \end{equation}
  and that on this interval, the two terms forming the difference are by themselves of the claimed size.

  We first treat the case of large $s''$-values in (\ref{B.37}), that is, of $s''\ge 3S_0$.
  In this case, $s'\mapsto\phi_{\lambda(s')}(s''-s')$ and $s'\mapsto{\phi}^{\alpha}_{\lambda(s')}(s''-s')$
  are supported in $s'\in[s'',2s'']$. 
  In particular, $s'\ge S_0$ so that $\lambda\stackrel{(\ref{B.32})}{=}s'$. Hence (\ref{B.first-rep}) takes the form 
  \begin{eqnarray}\nonumber
  \dppa(s'')&=&
  -\int_{-\infty}^\infty\frac{1}{(1+\hat s)^2}( (\hat w_0^{\alpha})^2-\hat w_0^2)\frac{d\ttrho}{ds'}\,d\hat s\nonumber\\
  &&-\frac{1}{(s'')^2}\int_{-\infty}^\infty(1+\hat s)
  \left(\hat w_0^{\alpha}\frac{d^2\hat w_0^{\alpha}}{d\hat s^2}-\hat w_0\frac{d^2\hat w_0}{d\hat s^2}\right)\ttrho\,d\hat s\nonumber\\
  &&+\frac{2}{(s'')^3}\int_{-\infty}^\infty(1+\hat s)^2
  \left(\hat w_0^{\alpha}\frac{d^2\hat w_0^{\alpha}}{d\hat s^3}-\hat w_0\frac{d^3\hat w_0}{d\hat s^3}\right)\ttrho\,d\hat s\nonumber\\
  &&+\frac{1}{(s'')^4}\int_{-\infty}^\infty(1+\hat s)^3
  \left(\hat w_0^{\alpha}\frac{d^2\hat w_0^{\alpha}}{d\hat s^4}-\hat w_0\frac{d^4\hat w_0}{d\hat s^4}\right)\ttrho\,d\hat s\,.\nonumber
  \end{eqnarray}
  In particular, we also have $s'\ge 3S_0$ so that
  $\ttrho(s')\stackrel{(\ref{B.28})}{=}\trho(s')\stackrel{(\ref{B.17})}{=}1-\frac{1}{\frac{s'}{S_0}-1}=1-\frac{S_0}{s'-S_0}\leq 1$, and thus
  $\frac{d\tt\rho}{ds'}=\frac{S_0}{(s'-S_0)^2}\leq \left(\frac{3}{2}\right)^2 S_0\frac{1}{(s'')^2}$ since $s'\geq s''\geq 3S_0$.
   Hence the above representation yields
  \begin{eqnarray}\nonumber
  |\dppa(s'')|&\lesssim&
  \frac{S_0}{(s'')^2}\int_{-\infty}^\infty\left|(\hat w_0^{\alpha})^2-\hat w_0^2\right|\,d\hat s\nonumber\\
  &&+\frac{1}{(s'')^2}\int_{-\infty}^\infty(1+\hat s)
  \left|\hat w_0^{\alpha}\frac{d^2\hat w_0^{\alpha}}{d\hat s^2}-\hat w_0\frac{d^2\hat w_0}{d\hat s^2}\right|\,d\hat s\nonumber\\
  &&+\frac{1}{(s'')^3}\int_{-\infty}^\infty(1+\hat s)^2
  \left|\hat w_0^{\alpha}\frac{d^3\hat w_0^{\alpha}}{d\hat s^3}-\hat w_0\frac{d^3\hat w_0}{d\hat s^3}\right|\,d\hat s\nonumber\\
  &&+\frac{1}{(s'')^4}\int_{-\infty}^\infty(1+\hat s)^3
  \left|\hat w_0^{\alpha}\frac{d^4\hat w_0^{\alpha}}{d\hat s^4}-\hat w_0\frac{d^4\hat w_0}{d\hat s^4}\right|\,d\hat s\,.\nonumber
  \end{eqnarray}
  Using (\ref{extraline}) and inserting the estimate (\ref{B.40}) for $k=0,2,3,4$ we obtain,
  as claimed in (\ref{B-new}),
  \begin{equation}\label{B.43}
  |\dppa(s'')|\;\leq\;C\alpha\left(\frac{\alpha^4}{(s'')^2}+\frac{\alpha^2}{(s'')^2}
  +\frac{\alpha}{(s'')^3}+\frac{1}{(s'')^4}\right)\;\lesssim\;C\frac{\alpha}{(s'')^2}
  \quad\mbox{for}\;s''\ge 3S_0.
  \end{equation}
  We now address the small $s''$-values, that is, $s''\le\frac{S_0}{2}$.
  In this case, $s'\mapsto\phi_{\lambda(s')}(s''-s')$ and $s'\mapsto\phi_{\lambda(s')}^{\alpha}(s''-s')$
  are supported in $s'\in[s'',s''+\frac{S_0}{2}]$.
  In particular, $s'\le S_0$ so that $\lambda\stackrel{(\ref{B.32})}{=}S_0$. Hence by (\ref{B.6}) and (\ref{phi_1})
   we obtain the representation
  \begin{eqnarray*}
  \dppa(s'')&=&
  -\frac{2}{S_0^2}\int_{-\infty}^\infty
  \left(\hat w_0^{\alpha}\frac{d\hat w_0^{\alpha}}{d\hat s}-\hat w_0^{\alpha}\frac{d\hat w_0^{\alpha}}{d\hat s}\right)\ttrho\,d\hat s\\
  &&-\frac{1}{S_0^3}\int_{-\infty}^\infty
  \left(\hat w_0^{\alpha}\frac{d^2\hat w_0^{\alpha}}{d\hat s^2}-\hat w_0\frac{d^2\hat w_0}{d\hat s^2}\right)\ttrho\,d\hat s\\
  &&+\frac{2}{S_0^4}\int_{-\infty}^\infty
  \left(\hat w_0^{\alpha}\frac{d^3\hat w_0^{\alpha}}{d\hat s^3}-\hat w_0\frac{d^3\hat w_0}{d\hat s^3}\right)\ttrho\,d\hat s\\
  &&+\frac{1}{S_0^5}\int_{-\infty}^\infty
  \left(\hat w_0^{\alpha}\frac{d^4\hat w_0^{\alpha}}{d\hat s^4}-\hat w_0\frac{d^4\hat w_0}{d\hat s^4}\right)\ttrho\,d\hat s\,.
  \end{eqnarray*}
  Moreover, $s'\le S_0$ implies $\rho(s')=0$ (cf.\  (\ref{B.17})\&(\ref{B.24})), $\eta_0(\frac{s'}{S_0})=1$ and (cf.\ (\ref{B.31})) thus
  $\widetilde\rho(s')=\varepsilon\exp(\frac{s'}{S_2})$. 
  In terms of $\hat s$ given by
  $s'=s''-S_0 \hat s$, this translates into
  $\widetilde\rho(s')=\varepsilon\exp(\frac{s''}{S_2})\exp(-\frac{S_0}{S_2}\hat s)$.
  Hence the above representation specifies to
  \begin{eqnarray}\nonumber
  \dppa(s'')&=&
  -\varepsilon\frac{2\exp(\frac{s''}{S_2})}{S_0^2 }\int_{-\infty}^\infty
  \left(\hat w_0^{\alpha}\frac{d\hat w_0^{\alpha}}{d\hat s}-\hat w_0\frac{d\hat w_0}{d\hat s}\right)\exp\left(-\frac{S_0}{S_2}\hat s\right)\,d\hat s\nonumber\\
  &&-\varepsilon\frac{\exp(\frac{s''}{S_2})}{S_0^3}\int_{-\infty}^\infty
  \left(\hat w_0^{\alpha}\frac{d^2\hat w_0^{\alpha}}{d\hat s^2}-\hat w_0\frac{d^2\hat w_0}{d\hat s^2}\right)\exp\left(-\frac{S_0}{S_2}\hat s\right)\,d\hat s\nonumber\\
  &&+\varepsilon\frac{2\exp(\frac{s''}{S_2})}{S_0^4}\int_{-\infty}^\infty
  \left(\hat w_0^{\alpha}\frac{d^3\hat w_0^{\alpha}}{d\hat s^3}-\hat w_0\frac{d^3\hat w_0}{d\hat s^3}\right)\exp\left(-\frac{S_0}{S_2}\hat s\right)\,d\hat s\nonumber\\
  &&+\varepsilon\frac{\exp(\frac{s''}{S_2})}{S_0^4}\int_{-\infty}^\infty
  \left(\hat w_0^{\alpha}\frac{d^4\hat w_0^{\alpha}}{d\hat s^4}-\hat w_0\frac{d^4\hat w_0}{d\hat s^4}\right)\exp\left(-\frac{S_0}{S_2}\hat s\right)\,d\hat s\,.\nonumber
  \end{eqnarray}
  Inserting the estimate (\ref{B.40}) for $k=1,2,3,4$ and using (\ref{extraline}) we obtain as claimed in (\ref{B-new})
  \begin{equation}\label{B.44}
  |\dppa(s'')|\;\lesssim\;\alpha\varepsilon\exp\left(\frac{s''}{S_2}\right)(\alpha^3+\alpha^2+\alpha^1+1)
  \;\lesssim\;\alpha\exp\left(\frac{s''}{S_2}\right)
  \quad\mbox{for}\;s''\le\frac{S_0}{2}.
  \end{equation}

  We finally address the intermediate values of $s''$, that is, $\frac{S_0}{2}\le s''\le 3S_0$.
  Splitting the $ds'$-integrals into $s'\in[S_0,\infty)$ and $s'\in(-\infty,S_0]$ in order
  to treat $\lambda\stackrel{(\ref{B.32})}{=}\max\{s',S_0\}$, we obtain
  \begin{eqnarray}\nonumber
  \dppa(s'')&=&
  -\frac{2}{ s''   }\int_{-\infty}^{\frac{s''}{S_0}-1}
  \left(\hat w_0^{\alpha}\frac{d  \hat w_0^{\alpha}}{d\hat s  }-\hat w_0\frac{d  \hat w_0}{d\hat s  }\right)\ttrho\,d\hat s\nonumber\\
  &&-\frac{1}{(s'')^2}\int_{-\infty}^{\frac{s''}{S_0}-1}(1+\hat s)
  \left(\hat w_0^{\alpha}\frac{d^2\hat w_0^{\alpha}}{d\hat s^2}-\hat w_0\frac{d^2\hat w_0}{d\hat s^2}\right)\ttrho\,d\hat s\nonumber\\
  &&+\frac{2}{(s'')^3}\int_{-\infty}^{\frac{s''}{S_0}-1}(1+\hat s)^2
  \left(\hat w_0^{\alpha}\frac{d^3\hat w_0^{\alpha}}{d\hat s^3}-\hat w_0\frac{d^3\hat w_0}{d\hat s^3}\right)\ttrho\,d\hat s\nonumber\\
  &&+\frac{1}{(s'')^4}\int_{-\infty}^{\frac{s''}{S_0}-1}(1+\hat s)^3
  \left(\hat w_0^{\alpha}\frac{d^4\hat w_0^{\alpha}}{d\hat s^4}-\hat w_0\frac{d^4\hat w_0}{d\hat s^4}\right)\ttrho\,d\hat s\nonumber\\
  &&-\frac{2}{S_0^2  }\int_{\frac{s''}{S_0}-1}^\infty
  \left(\hat w_0^{\alpha}\frac{d\hat w_0^{\alpha}}{d\hat s}-\hat w_0\frac{d\hat w_0}{d\hat s}\right)\ttrho\,d\hat s\nonumber\\
  &&-\frac{1}{S_0^3}\int_{\frac{s''}{S_0}-1}^\infty
  \left(\hat w_0^{\alpha}\frac{d^2\hat w_0^{\alpha}}{d\hat s^2}-\hat w_0\frac{d^2\hat w_0}{d\hat s^2}\right)\ttrho\,d\hat s\nonumber\\
  &&+\frac{2}{S_0^4}\int_{\frac{s''}{S_0}-1}^\infty
  \left(\hat w_0^{\alpha}\frac{d^3\hat w_0^{\alpha}}{d\hat s^3}-\hat w_0\frac{d^3\hat w_0}{d\hat s^3}\right)\ttrho\,d\hat s\nonumber\\
  &&+\frac{1}{S_0^5}\int_{\frac{s''}{S_0}-1}^\infty
  \left(\hat w_0^{\alpha}\frac{d^4\hat w_0^{\alpha}}{d\hat s^4}-\hat w_0\frac{d^4\hat w_0}{d\hat s^4}\right)\ttrho\,d\hat s\,.\nonumber
  \end{eqnarray}
  Since $|\ttrho|\le 1$ and since $|\frac{1}{s''}|\le\frac{2}{S_0}$, we obtain from using (\ref{extraline}) and 
  inserting the estimate (\ref{B.40}) for $k=1,2,3,4$:	
  \begin{equation}\label{B.45}
  |\dppa(s'')|\;\lesssim\;\alpha(\alpha^3+\alpha^2+\alpha^1+1)\;\lesssim\;\alpha
  \quad\mbox{for}\;s''\in\left[\frac{S_0}{2},2S_0\right]\,.
  \end{equation}
  Now (\ref{B.43}), (\ref{B.44}) and (\ref{B.45}) establish (\ref{B-new}).
  
  As in the proof of (\ref{upper-bound-OSC}) in Section \ref{heuristics}, in order to obtain (\ref{upper-bound-OSC}) in the range $s''\geq S_1$,
   we need to cut-off the measure
  $\widetilde{\rho}$ (defined in (\ref{B.17})\&(\ref{B.24}) and modified in (\ref{B.28})) in the range $\frac{S_1}{2}\leq s'\leq S_1$ so that 
   \begin{equation*}
    \phi_1(s'')=\int_{-\infty}^{\infty}\phi^{\alpha}_{\lambda(s')}(s''-s')\widetilde{\rho}(s')\eta\left(\frac{s'}{S_1}\right)\, ds'\,
   \end{equation*}
   satisfies (\ref{upper-bound-OSC}).
    In this region ($s'\geq \frac{S_1}{2}$ or $s''\geq \frac{S_1}{4}$) the modification (\ref{B.28})
   of $\rho$ and (\ref{B.32}) of $\lambda$, are not effective. So we may directly quote the argument of Section \ref{heuristics} for the modification of $\widetilde{\rho}$ through $\eta$. Note that this argument is unaffected by having replaced $\hat w_0$ by $\hat w_0^{\alpha}$. This concludes the proof of (\ref{upper-bound-OSC}) and (\ref{weak-lower-bound}).


  \subsection{Approximate positivity in the boundary layers: proof of Lemma \ref{Lemma3}}\label{Proof-Lemma3}

  The approximate non-negativity of $\hat\xi_0$, cf.\ (\ref{convolution}), is
  lost in the boundary layer $s\ll -1$, cf.\ (\ref{1}). However, in this subsection we
  show that $\hat\xi_0$ cannot be too negative in the boundary layer
  {\it provided $\hat\xi_0$ is sufficiently small in the
  transition region $|s|\lesssim 1$}. 
  We recall the statement of Lemma \ref{Lemma3}: For all $S_2\gg 1$ and $\varepsilon \leq 1$ we have 
  \begin{equation*}
  -\int_{-S_2}^{-1}\hat \xi_0 ds\lesssim \frac{1}{\varepsilon}\int_{-1}^{0}\hat\xi_0 ds+\frac{1}{\varepsilon}+\int_{-S_2}^{-S_2+1}|\hat\xi_0|ds+\varepsilon \exp(5S_2)\,.
  \end{equation*}
%
  With the standard rescaling (cf.\ proof of Theorem \ref{Th1})
  \begin{equation*}
  s\;\leadsto\;s+S_0,\quad
  \hat\xi\;\leadsto\;\exp(-3S_0)\hat\xi
  \quad\mbox{and thus also}\quad
  \hat\xi_0\;\leadsto\;\exp(-3S_0)\hat\xi_0\,,
  \end{equation*}
  the above turns into
    \begin{equation}\label{approx-pos-BL}
 -\int_{-S_2-S_0}^{-S_0-1}\hat\xi_0\,ds
  \lesssim
  \frac{1}{\varepsilon}\int_{-S_0-1}^{-S_0}\hat\xi_0\,ds
  +\frac{1}{\varepsilon}\exp(3S_0)	
  +\int_{-S_2-S_0}^{-S_2-S_0+1}|\hat\xi_0|\,ds
  +\varepsilon\exp(5S_2+3S_0)\,.
  \end{equation}
  Hence it is enough to establish the latter for \textit{some} $S_0$. In fact, we shall show
  that for \textit{all} $S_0\gg 1$ and $S_1\gg S_0$
  \begin{equation}\label{4.42}
 -\int_{-S_1}^{-S_0-1}\hat\xi_0\,ds
  \lesssim
  \frac{1}{\varepsilon}\int_{-S_0-1}^{-S_0}\hat\xi_0\,ds
  +\frac{1}{\varepsilon}\exp(5S_0)	
  +\int_{-S_1}^{-S_1+1}|\hat\xi_0|\,ds
  +\varepsilon\exp(5S_1)\,,
  \end{equation}
  where we write $S_1=S_2+S_0$\,.
  Indeed, fixing an order-one $S_0$ which is sufficiently large so that (\ref{4.42}) is valid, we obtain (\ref{approx-pos-BL}).
  %
  Multiplying both sides of (\ref{start}) by $\phi_0(s')$ (see definition (\ref{convolution-kernel})) and integrating in $(-\infty, \infty)$,
  we deduce
  \begin{equation}\label{4.20}
  \int_{-\infty}^\infty\hat\xi_0\,\phi\,ds\ge
  -C\left\{\int_0^1\left[\frac{d}{d\hat z}\left(-\frac{d^2}{d\hat z^2}+1\right)^2w\right]^2
  +\left[\left(-\frac{d^2}{d\hat z^2}+1\right)^2w\right]^2\,d\hat z\right\}\,,
  \end{equation}
  for any smooth $w$, supported in $\hat z\in[0,1]$ and satisfying the
  boundary conditions $w=\frac{dw}{d\hat z}=\left(-\frac{d^2}{d\hat z^2}+1\right)^2w=0$ at $\hat z=0$, where
  as before we use the abbreviation
  \begin{equation*}
  \phi\;:=\;w\left(-\frac{d^2}{d\hat z^2}+1\right)^2w\,
  \end{equation*}
  for the multiplier. 
  This time, $w$ will {\it not} be compactly supported
  in $\hat z\in(0,1]$ (only in $[0,1]$) so that the boundary conditions matters. Using
  the fact that the function $\hat z\sinh\hat z$ satisfies these
  boundary conditions, we enforce them for $w$ by the Ansatz
  \begin{equation}\label{4.18}
  w\;=\;(\hat z\sinh\hat z)\hat w\quad\mbox{with}\quad\hat w=const\;\;\mbox{for}\;\;\hat z\ll 1\,.
  \end{equation}
  As in the previous subsections, it is more telling to express (\ref{4.20})
  in terms of the $s$-variable $s=\ln \hat z$. Appealing to the representations 
  \begin{eqnarray}\label{box}
  &&(-\partial^2_{\hat z}+1)^2
  \hat z\sinh\hat z\notag\\
  &&=\hat z^{-2}\Big(\hat z^{-1}\sinh\hat z\,(\partial_s-2)(\partial_s-1)
  +4\cosh\hat z(\partial_s-1)
  +4\hat z\sinh\hat z\Big)
 (\partial_s+1)\partial_s\,
  \end{eqnarray}
  and
  \begin{eqnarray}\label{boxbox}
  &&\partial_{\hat z}(-\partial^2_{\hat z}+1)^2	
  \hat z\sinh\hat z\notag\\
  &&=\hat z^{-3}\Big(
  \hat z^{-1}\sinh\hat z\,(\partial_s-3)(\partial_s-2)(\partial_s-1)
  +5\cosh\hat z(\partial_s-2)(\partial_s-1)\notag\\
  &&\;\;+8\hat z\sinh\hat z(\partial_s-1)
  +4\hat z^2\cosh\hat z\Big)
  (\partial_s+1)\partial_s
  \end{eqnarray}
 (their proofs are reported in the Appendix to Section \ref{Proof-Lemma3}) we obtain
  \begin{equation}\label{4.20bis}
  \int_{-\infty}^\infty\hat\xi_0\phi\,ds\ge
  -C\int_{-\infty}^\infty\exp(-5s)
  \left[\left(\frac{d^5\hat w}{ds^5}\right)^2+\cdots+\left(\frac{d\hat w}{ds}\right)^2\right]\,ds\,,
  \end{equation}
  where according to the formula
  \begin{eqnarray*}
  &&\hat z\sinh\hat z\left(-\frac{d^2}{d\hat z^2}+1\right)^2\hat z\sinh\hat z\\
  &&=\left(\frac{d}{ds}+1\right)\frac{d}{ds}\left(\frac{\sinh\hat z}{\hat z}\right)^2
  \left(\frac{d}{ds}+1\right)\frac{d}{ds}-2\left(\frac{d}{ds}+1\right)\frac{d}{ds}\,,\nonumber
  \end{eqnarray*}
  (the argument for the formula above is given at the end of this section, see (\ref{4.11})), the multiplier is given by
  \begin{equation}\label{4.27bis}
  \phi\;=\;\hat w
  \left(\frac{d}{ds}+1\right)\frac{d}{ds}\left[\left(\frac{\sinh\hat z}{\hat z}\right)^2
  \left(\frac{d}{ds}+1\right)\frac{d}{ds}-2\right]\hat w\,.
  \end{equation}

  We now make the following Ansatz for $\hat w$:
  \begin{equation}\label{4.19}
  \hat w\;=\;\frac{1}{\sqrt{\varepsilon}}\hat w_0+\sqrt{\varepsilon}\hat w_1\,,
  \end{equation}
  with the constraints
  \begin{equation}\label{4.22}
  \hat w_0\;=\;\left\{\begin{array}{llll}
  1&\mbox{for}&s\le& -S_0-1\\
  0&\mbox{for}&s\ge& -S_0\end{array}\right\},\quad
  \hat w_1\;=\;\left\{\begin{array}{llll}
  const&\mbox{for}&s\le& -S_1\\
  0&\mbox{for}&s\ge& -S_0-1\end{array}\right\}\,,
  \end{equation}
  so that (\ref{4.18}) is satisfied.
  We don't want to specify the value of the constant appearing in the definition of $w$ since it will
  not appear in the future estimates.
  The merit of the Ansatz (\ref{4.19}) is that, because $\frac{d\hat w_0}{ds}$
  and $\frac{d\hat w_1}{ds}$ have disjoint support,
  the multiplier $\phi$, cf.\ (\ref{4.27bis}), splits into three parts
  \begin{equation}\label{4.41}
  \phi\;=\;\frac{1}{\varepsilon}\phi_0+\phi_{01}+\varepsilon\phi_1\,,
  \end{equation}
  where
  \begin{eqnarray}\nonumber
  \phi_0&:=&\hat w_0
  \left(\frac{d}{ds}+1\right)\frac{d}{ds}\left[\left(\frac{\sinh\hat z}{\hat z}\right)^2
  \left(\frac{d}{ds}+1\right)\frac{d}{ds}-2\right]\hat w_0,
  \nonumber\\ 
  \phi_{01}&:=&\hat w_0\left(\frac{d}{ds}+1\right)\frac{d}{ds}\left[\left(\frac{\sinh\hat z}{\hat z}\right)^2
  \left(\frac{d}{ds}+1\right)\frac{d}{ds}-2\right]\hat w_1,
  \label{4.28}\\
  \phi_1&:=&\hat w_1\left(\frac{d}{ds}+1\right)\frac{d}{ds}\left[\left(\frac{\sinh\hat z}{\hat z}\right)^2
  \left(\frac{d}{ds}+1\right)\frac{d}{ds}-2\right]\hat w_1\,.\nonumber
  \end{eqnarray}
  As a related side effect of the disjoint support of the functions
  $\frac{d\hat w_0}{ds}$ and
  $\frac{d\hat w_1}{ds}$, the error term in (\ref{4.20bis})
  splits into two parts:	
  \begin{eqnarray}
  \lefteqn{\int_{-\infty}^\infty\exp(-5s)
  \left[\left(\frac{d^5\hat w}{ds^5}\right)^2+\cdots+\left(\frac{d\hat w}{ds}\right)^2\right]\,ds}\nonumber\\
  &=&\frac{1}{\varepsilon}
  \int_{-\infty}^\infty\exp(-5s)
  \left[\left(\frac{d^5\hat w_0}{ds^5}\right)^2+\cdots+\left(\frac{d\hat w_0}{ds}\right)^2\right]\,ds\label{4.30bis}\\
  &+&\varepsilon
  \int_{-\infty}^\infty\exp(-5s)
  \left[\left(\frac{d^5\hat w_1}{ds^5}\right)^2+\cdots+\left(\frac{d\hat w_1}{ds}\right)^2\right]\,ds\,.\label{4.30}
  \end{eqnarray}
  Hence in the sequel, we will have to consider five terms:
  \begin{itemize}
  \item three multiplier terms: 
  $\frac{1}{\varepsilon}\int_{-\infty}^\infty\hat\xi_0\phi_0\,ds$,
  $\int_{-\infty}^\infty\hat\xi_0\phi_{01}\,ds$, and
  $\varepsilon\int_{-\infty}^\infty\hat\xi_0\phi_1\,ds$,
  \item two error terms: the $\hat w_0$-error term (\ref{4.30bis})
  and the $\hat w_1$-error term (\ref{4.30}).
  \end{itemize}
  Below, we will construct $\hat w_1$ such that the mixed expression $\phi_{01}$,
  cf.\ (\ref{4.28}),  
  in the multiplier $\phi$ gives rise to the left-hand side of (\ref{4.42}).

  Before, we address the multiplier and the error term that only involve $\hat w_0$.
  Clearly, $\hat w_0$ can be chosen to satisfy $S_0$-independent bounds:
  $\sup_{s\in\mathbb{R}}|\hat w_0|,\;\cdots\;,
  \sup_{s\in\mathbb{R}}\left|\frac{d^5\hat w_0}{ds^5}\right|\;\lesssim\;1\,.$
  Hence in view of (\ref{4.22}), we obtain for the $\hat w_0$-error term (\ref{4.30bis})
  \begin{equation}\label{4.36}
  \frac{1}{\varepsilon}\int_{-\infty}^\infty\exp(-5s)
  \left[\left(\frac{d^5\hat w_0}{ds^5}\right)^2+\cdots+\left(\frac{d\hat w_0}{ds}\right)^2\right]\,ds
  \;\lesssim\;\frac{1}{\varepsilon}\exp(5S_0).
  \end{equation}
  Moreover, in view of (\ref{4.22}), we obtain
  \begin{equation}\label{4.31}
  |\phi_0|\;\le\;
  \left\{
   \begin{array}{lrrlll}
  0      \quad        &\mbox{for }&\;         &s& \le  -S_0-1\,&\\
  C_0      \quad      &\mbox{for }&\;-S_0-1\le&s&\le -S_0\,&\\
  0         \quad     &\mbox{for }&\;  -S_0 \leq      &s&\,&
  \end{array}\right\}\,,
  \end{equation}
  where we momentarily want to remember the value of the universal constant $C_0$.
  Since, by (\ref{1}) in Lemma \ref{Lemma1}, there exists a specific constant $C_1$ such that 
  $\int_{-S_0-1}^{-S_0}\hat\xi_0\,ds+C_1\exp(3S_0)\geq 0$,
  we obtain from (\ref{4.31})
\begin{eqnarray}\nonumber
\int_{-S_0-1}^{-S_0}\hat\xi_0(\phi_0-C_0)\,ds
&=&
\int_{-S_0-1}^{-S_0}(-\hat\xi_0)(-\phi_0+C_0)\,ds\nonumber\\
&\stackrel{(\ref{1})}{\le}&
C_1\int_{-S_0-1}^{-S_0}\exp(-3s)(-\phi_0+C_0)\,ds\nonumber\\
&\stackrel{(\ref{4.31})}{\le}&
2C_1C_0\int_{-S_0-1}^{-S_0}\exp(-3s)\,ds\nonumber\\
&\le&C_1C_0\exp(3S_0)\,,\label{same-argument}
\end{eqnarray}
so that for the $\phi_0$-multiplier term we obtain 
\begin{eqnarray}
\frac{1}{\varepsilon}\int_{-\infty}^\infty\hat\xi_0\phi_0\,ds
&\stackrel{(\ref{4.31})}{=}&
\frac{1}{\varepsilon}\int_{-S_0-1}^{-S_0}\hat\xi_0\phi_0\,ds\nonumber\\
&\le&\frac{1}{\varepsilon}\left(C_0\int_{-S_0-1}^{-S_0}\hat\xi_0\,ds+C_1C_0\exp(3S_0)\right)\,\nonumber\\
&\lesssim&\frac{1}{\varepsilon}\left(\int_{-S_0-1}^{-S_0}\hat\xi_0\,ds+C_1\exp(3S_0)\right)\,.\label{4.37}
\end{eqnarray}

  We now specify $\hat w_1$ with the goal that $\phi_{01}$,
  cf.\ (\ref{4.28}), gives rise to the l.\ h.\ s.\ of  (\ref{4.42}). 
  This motivates the construction of a universal function $\hat w_2$ with the property that
  \begin{equation}\label{4.26}
  \left(\frac{d}{ds}+1\right)\frac{d}{ds}\left[\left(\frac{\sinh\hat z}{\hat z}\right)^2
  \left(\frac{d}{ds}+1\right)\frac{d}{ds}-2\right]\hat w_2\;=\;1
  \quad\mbox{for}\;s\ll -1\,,
  \end{equation}
  which will be carried out below in such a way that
  \begin{equation}\label{4.29}
  \frac{|\hat w_2|}{|s|+1},\;
  \left|\frac{d\hat w_2}{ds}\right|,\;\cdots\;
  \left|\frac{d^5\hat w_2}{ds^5}\right|\;\lesssim\;1\,.
  \end{equation}
  Equipped with $\hat w_2$, we now make the Ansatz of blending $\hat w_2$ to $\hat w_2(-S_1)$ for $s<-S_1$
  and to zero for $s<-S_0-1$:
  \begin{equation}\label{4.25}
  \hat w_1(s)\;=\;\eta(s+S_1)\eta(-(s+S_0+1))\hat w_2(s)+(1-\eta(s+S_1))\hat w_2(-S_1)\,,
  \end{equation}
  where $\eta$ is a universal cut-off function with 
  \begin{equation}\label{4.24}
  \eta(s)=\left\{\begin{array}{ll}
  0 & \mbox{ for }  s\le 0\\
  1 & \mbox{ for } s\ge 1\
  \end{array}\right\}\,,
  \end{equation}
  so that (\ref{4.22}) is satisfied. The main merit of Ansatz (\ref{4.25}) \& (\ref{4.24}) 
  is that it makes use of (\ref{4.26})
  which by definition (\ref{4.28}) yields
  \begin{equation}\label{4.32}
  \phi_{01}=\left\{
 \begin{array}{lllll}
  0 &\mbox{ for }    \qquad     &s&\le  -S_1\,& \\
 1  &\mbox{ for }   \qquad -S_1+1\le &s&\le -S_0-2\,&\\
  0  &\mbox{ for }      \qquad      -S_0-1\le  &s&\,&
  \end{array}\right\}\,.
  \end{equation}
  Furthermore, the estimates (\ref{4.29}) turn into
  \begin{equation}\label{4.33}
  |\phi_{01}|,\;
  \frac{|\hat w_1|}{S_1},\;\left|\frac{d\hat w_1}{ds}\right|\,,\cdots
  \,,\left|\frac{d^5\hat w_1}{ds^5}\right|\;\leq\;C_0\,.
  \end{equation}
  In particular, we obtain for the $\phi_{01}$-multiplier term
   \begin{eqnarray}
  &C_0\displaystyle{\int_{-\infty}^\infty}&\hat\xi_0\phi_{01}\,ds\nonumber\\
  &\stackrel{(\ref{4.32})}{=}&
  C_0\int_{-S_1}^{-S_0-1}\hat\xi_0\,ds
  +\int_{-S_1}^{-S_1+1}\hat\xi_0(\phi_{01}-C_0)\,ds
  +\int_{-S_0-2}^{-S_0-1}\hat\xi_0(\phi_{01}-C_0)\,ds\nonumber\\
  &\stackrel{(\ref{4.33})}{\le}&
  C_0\int_{-S_1}^{-S_0-1}\hat\xi_0\,ds+2C_0\int_{-S_1}^{-S_1+1}|\hat\xi_0|\,ds
  +C_1C_0\exp(3(S_0+1))\label{4.38}\,,
  \end{eqnarray}
  where for $\int_{-S_0-2}^{-S_0-1}\hat\xi_0(\phi_{01}-1)\,ds$, we have
  used the same argument as in (\ref{same-argument}). 
%

  Because of $\phi_1=\hat w_1\phi_{01}$ another consequence
  of (\ref{4.33}) and (\ref{4.32}) is 
  \begin{equation*}
  |\phi_1|\;\lesssim\;\left\{
  \begin{array}{llrlll}
  0     &\mbox{for }&  \qquad &s&\le -S_1\\
  S_1   &\mbox{for}&    \qquad-S_1\le &s&\le -S_0-1\\
  0     &\mbox{for}&    \qquad -S_0-1\le &s& 
  \end{array}\right\}\,.
  \end{equation*}
  By the same argument that leads to (\ref{4.37}), this implies
  for the $\phi_1$-multiplier term	
  \begin{equation}\label{4.39}
  \varepsilon\int_{-\infty}^\infty\hat\xi_0\phi_1\,ds
  \lesssim\varepsilon S_1\left(\int_{-S_1}^{-S_0-1}\hat\xi_0\,ds+C_1\exp(3(S_1-1))\right)\,.
  \end{equation}
  We finally address the $\hat w_1$-error term (\ref{4.30}): It follows from (\ref{4.22})
  and (\ref{4.33}) that 
  \begin{equation}\label{4.40}
  \varepsilon\int_{-\infty}^\infty\exp(-5s)
  \left[\left(\frac{d^5\hat w_1}{ds^5}\right)^2+\cdots+\left(\frac{d\hat w_1}{ds}\right)^2\right]\,ds
  \;\le\;C\varepsilon\exp(5S_1)\,.
  \end{equation}

  We now collect the five estimates (\ref{4.36}), (\ref{4.37}), (\ref{4.38}), (\ref{4.39}), and (\ref{4.40}).
   Via (\ref{4.41}) and (\ref{4.30}) we obtain from (\ref{4.20bis}) that
  \begin{eqnarray*}
  \lefteqn{-C_0\int_{-S_1}^{-S_0-1}\hat\xi_0\,ds}\nonumber\\
  &\lesssim&\frac{1}{\varepsilon}\exp(5S_0) +\frac{1}{\varepsilon}\left(\int_{-S_0-1}^{-S_0}\hat\xi_0\,ds+C_1\exp(3S_0)\right)\\
  &+&C_0\int_{-S_1}^{-S_1+1}|\hat\xi_0|\,ds
  +C_0C_1\exp(3S_0)\\
  &+&\varepsilon S_1\left(\int_{-S_1}^{-S_0-1}\hat\xi_0\,ds+C_1\exp(3S_1)\right)
  +\varepsilon\exp(5S_1)\,,
  \end{eqnarray*}
  where we recall that $C_1$ was chosen such that the terms in the parentheses are non-negative. Hence we may discard the term
  $\int_{-S_1}^{-S_0-1}\hat \xi_0\, ds$ on the r.\ h.\ s.\ : If it is negative we may omit it; if it is positive, then the estimate
  comes for free. Dividing by $C_0$ we thus obtain 
  \begin{eqnarray*}
  -\int_{-S_1}^{-S_0-1}\hat \xi_0\, ds &\lesssim&\frac{1}{\varepsilon}\exp(5S_0)+
  \frac{1}{\varepsilon}\left(\int_{-S_0-1}^{-S_0}\hat \xi_0\, ds+C_1\exp(3S_0)\right)\\
  &+&\int_{-S_1}^{-S_1+1}|\hat\xi_0|\, ds+\exp(3S_0)+
  \varepsilon S_1\exp(3S_1)+\varepsilon \exp(5S_1)\,.
  \end{eqnarray*}
%
 which implies (\ref{4.42}) because $\varepsilon\leq 1$ and $S_1\gg S_0\gg 1$.

  We derive now the operator-valued formula
  \begin{eqnarray}\label{4.11}
  &&\hat z\sinh\hat z\left(-\frac{d^2}{d\hat z^2}+1\right)^2\hat z\sinh\hat z\nonumber\\
  &&=\left(\frac{d}{ds}+1\right)\frac{d}{ds}\left(\frac{\sinh\hat z}{\hat z}\right)^2
  \left(\frac{d}{ds}+1\right)\frac{d}{ds}-2\left(\frac{d}{ds}+1\right)\frac{d}{ds}\,,
  \end{eqnarray}
  that is a non-homogeneous generalization of 
  $\hat z^2\frac{d^4}{d\hat z^4}\hat z^2
  =(\frac{d}{ds}+2)(\frac{d}{ds}+1)\frac{d}{ds}(\frac{d}{ds}-1)$ (cf.\ (\ref{poly})).
  The fairly simple structure of this formula is not a surprise:
  Since the functions $\sinh\hat z$ and $\hat z\sinh\hat z$ are in the kernel of
  $(-\frac{d^2}{d\hat z^2}+1)^2$, the functions $1$ and $\hat z^{-1}$ are in
  the kernel of $(-\frac{d^2}{d\hat z^2}+1)^2\hat z\sinh\hat z$.
  In $s$ coordinates, these functions are $1$ and $\exp(-s)$,
  respectively. This explains the {\it right} factor $(\frac{d}{ds}+1)\frac{d}{ds}$ 
  on the r.\ h.\ s.\ of (\ref{4.11}). On the other hand, the {\it adjoint}
  of the l.\ h.\ s.\ of (\ref{4.11}) {\it w.\ r.\ t. to the measure}
  $\frac{d\hat z}{\hat z}=ds$ is given by 
  $\hat z^2\sinh\hat z(-\frac{d^2}{d\hat z^2}+1)^2\sinh\hat z$
  and thus has a kernel containing $1$ and $\hat z=\exp(s)$. Hence the
  adjoint of the r.\ h.\ s.\ of (\ref{4.11}) w.\ r.\ t.\ $ds$
  has to contain the right factor $(\frac{d}{ds}-1)\frac{d}{ds}$,
  which means that the operator itself should contain the {\it left} 
  factor $(\frac{d}{ds}+1)\frac{d}{ds}$.

  We claim that the formula (\ref{4.11}) can be factorized into the two formulas
  \begin{eqnarray}
  \left(\frac{d^2}{d\hat z^2}-1\right)\hat z\sinh\hat z
  &=&\left(\frac{\sinh\hat z}{\hat z}\frac{d}{ds}+2\cosh\hat z\right)\left(\frac{d}{ds}+1\right),\label{4.11bis}\\
  \hat z\sinh\hat z \left(\frac{d^2}{d\hat z^2}-1\right)
  &=&\frac{d}{ds}
  \left[\left(\frac{d}{ds}+1\right)\frac{\sinh\hat z}{\hat z}-2\cosh\hat z\right]\,.\label{4.12}
  \end{eqnarray}
  Indeed, the composition of (\ref{4.11bis}) and (\ref{4.12}) yields
  \begin{eqnarray}\nonumber
  \lefteqn{\hat z\sinh\hat z \left(-\frac{d^2}{d\hat z^2}+1\right)^2\hat z\sinh\hat z}\nonumber\\
  &=&
  \frac{d}{ds}\left(\frac{d}{ds}+1\right)\left(\frac{\sinh\hat z}{\hat z}\right)^2
  \frac{d}{ds}\left(\frac{d}{ds}+1\right)\nonumber\\
  &&-2\frac{d}{ds}\cosh\hat z\frac{\sinh\hat z}{\hat z}\frac{d}{ds}\left(\frac{d}{ds}+1\right)
  +2\frac{d}{ds}\left(\frac{d}{ds}+1\right)\cosh\hat z\frac{\sinh\hat z}{\hat z}\left(\frac{d}{ds}+1\right)
  \nonumber\\
  &&-4\frac{d}{ds}(\cosh\hat z)^2\left(\frac{d}{ds}+1\right)\nonumber\\
  &=&
  \frac{d}{ds}\left(\frac{d}{ds}+1\right)\left(\frac{\sinh\hat z}{\hat z}\right)^2
  \frac{d}{ds}\left(\frac{d}{ds}+1\right)\nonumber\\
  &&+2\frac{d}{ds}\left(\frac{d}{ds}\cosh\hat z\frac{\sinh\hat z}{\hat z}\right)\left(\frac{d}{ds}+1\right)
  +2\frac{d}{ds}\cosh\hat z\frac{\sinh\hat z}{\hat z}\left(\frac{d}{ds}+1\right)
  \nonumber\\
  &&-4\frac{d}{ds}(\cosh\hat z)^2\left(\frac{d}{ds}+1\right)\nonumber\\
  &=&
  \frac{d}{ds}\left(\frac{d}{ds}+1\right)\left(\frac{\sinh\hat z}{\hat z}\right)^2
  \frac{d}{ds}\left(\frac{d}{ds}+1\right)\nonumber\\
  &&+2\frac{d}{ds}\left[\left(
  \frac{d}{ds}\cosh\hat z\frac{\sinh\hat z}{\hat z}\right)
  +\cosh\hat z\frac{\sinh\hat z}{\hat z}
  -2(\cosh\hat z)^2\right]\left(\frac{d}{ds}+1\right)\,,\label{4.13}
  \end{eqnarray}
  where $\left(\frac{d}{ds}\cosh\hat z\frac{\sinh\hat z}{\hat z}\right)$ denotes
  the multiplication with the $s$-derivative of the function
  $\cosh\hat z\frac{\sinh\hat z}{\hat z}$.
  This implies (\ref{4.11}) since because of
  \begin{equation*}
  \left(\frac{d}{ds}\cosh\hat z\frac{\sinh\hat z}{\hat z}\right)=
  \hat z\left(\frac{d}{d\hat z}\cosh\hat z\frac{\sinh\hat z}{\hat z}\right)=
  (\sinh\hat z)^2+(\cosh\hat z)^2-\cosh\hat z\frac{\sinh\hat z}{\hat z}\,,
  \end{equation*}
  the factor in the last term of (\ref{4.13}) simplifies to
  \begin{equation*}
  \left(\frac{d}{ds}\cosh\hat z\frac{\sinh\hat z}{\hat z}\right)
  +\cosh\hat z\frac{\sinh\hat z}{\hat z}
  -2(\cosh\hat z)^2=
  (\sinh\hat z)^2-(\cosh\hat z)^2=
  -1\,.
  \end{equation*}

  We now turn to the argument for (\ref{4.11bis}) and (\ref{4.12}). We first note that
  (\ref{4.11bis}) and (\ref{4.12}) reduce to
  \begin{eqnarray}
  \left(\frac{d^2}{d\hat z^2}-1\right)\hat z\exp\hat z
  &=&
  \left(\frac{\exp\hat z}{\hat z}\frac{d}{ds}+2\exp\hat z\right)\left(\frac{d}{ds}+1\right)\label{4.16bis}\\
  &=&
  \left(\exp\hat z\frac{d}{d\hat z}+2\exp\hat z\right)\left(\hat z \frac{d}{d\hat z}+1\right)\nonumber\\
  &=&
  \exp\hat z\left(\frac{d}{d\hat z}+2\right)\frac{d}{d\hat z}\hat z\qquad\mbox{and}\label{4.16}\\
  \hat z\exp\hat z \left(\frac{d^2}{d\hat z^2}-1\right)
  &=&\frac{d}{ds}
  \left[\left(\frac{d}{ds}+1\right)\frac{\exp\hat z}{\hat z}-2\exp\hat z\right]\label{4.17bis}\\
  &=&\hat z\frac{d}{d\hat z}
  \left[\left(\hat z \frac{d}{d\hat z}+1\right)\frac{\exp\hat z}{\hat z}-2\exp\hat z\right]\nonumber\\
  &=&\hat z\frac{d}{d\hat z}\left(\frac{d}{d\hat z}-2\right)\exp\hat z\,.\label{4.17}
  \end{eqnarray}
  Indeed, replacing $\hat z$ by $-\hat z$ in (\ref{4.16bis}), using the invariance
  of $\frac{d}{ds}=\hat z\frac{d}{d\hat z}$ under this change of variables, 
  and adding both identities yields (\ref{4.11bis}). Likewise, (\ref{4.17bis})
  yields (\ref{4.12}).
  The identities (\ref{4.16}) and (\ref{4.17}) can easily be checked using the commutator 
  relation $\frac{d}{d\hat z}\exp\hat z
  =\exp\hat z\left(\frac{d}{d\hat z}+1\right)$ on their left hand sides:
  \begin{eqnarray}\nonumber
  \left(\frac{d^2}{d\hat z^2}-1\right)\exp\hat z&=&\exp\hat z\left[\left(\frac{d}{d\hat z}+1\right)^2-1\right]
  =\exp\hat z\left(\frac{d}{d\hat z}+2\right)\frac{d}{d\hat z}\qquad\mbox{and}\nonumber\\
  \exp\hat z \left(\frac{d^2}{d\hat z^2}-1\right)&=&\left[\left(\frac{d}{d\hat z}-1\right)^2-1\right]\exp\hat z
  =\frac{d}{d\hat z}\left(\frac{d}{d\hat z}-2\right)\exp\hat z\,.\nonumber
  \end{eqnarray}
  This concludes the argument for (\ref{4.11}).

  We turn now to the construction of the function $\hat w_2$ with (\ref{4.26}) and (\ref{4.29}).
  We start by reducing (\ref{4.26}) to a second-order problem with bounded
  right-hand side: It is enough to construct a universal
  smooth $\hat v_2$ with
  \begin{equation}\label{4.2}
  \left[\frac{d}{ds}\left(\frac{\sinh\hat z}{\hat z}\right)^2
  \left(\frac{d}{ds}+1\right)-2\right]\,\hat v_2\;=\;1\quad\mbox{for}\;s\;\le\;-S_0
  \end{equation}
  and
  \begin{equation}\label{4.3}
  |\hat v_2|,\;\left|\frac{d\hat v_2}{ds}\right|,\;\cdots\;,\left|\frac{d^4\hat v_2}{ds^4}\right|
  \;\lesssim\;1\quad\mbox{for all}\;s\,.
  \end{equation}
  Indeed, consider the anti derivative $\hat w_2(s):=\int_{0}^{s}\hat v_2ds'$. Since
  $\frac{d\hat w_2}{ds}=\hat v_2$, the estimates
  (\ref{4.3}) turn into the estimates (\ref{4.29}). 
%
%
%
 Likewise (\ref{4.2}) turns into (\ref{4.26}) because of 
 $$\left(\frac{d}{ds}+1\right)\left[\frac{d}{ds}\left(\frac{\sinh \hat z}{\hat z}\right)^2\left(\frac{d}{ds}+1\right)-2\right]\frac{d}{ds}=
 \left(\frac{d}{ds}+1\right)\frac{d}{ds}\left[\left(\frac{\sinh \hat z}{\hat z}\right)^2\left(\frac{d}{ds}+1\right)\frac{d}{ds}-2\right]\,.$$
  We now extend (\ref{4.2}) to a problem on the entire line with
  nearly constant coefficients.
  Note that the coefficient $\left(\frac{\sinh\hat z}{\hat z}\right)^2$ is
  an entire, even function in $\hat z$ with value 1 at $\hat z=0$. Hence
  for every $S_0\gg 1$, we may write	
  \begin{equation*}
  \left(\frac{\sinh\hat z}{\hat z}\right)^2\;=\;1-a\quad\mbox{for all}\;s\le -S_0\,,
  \end{equation*}
  where
  \begin{equation}\label{4.8}
  \sup_{s\in\mathbb{R}}|a|,
  \sup_{s\in\mathbb{R}}\left|\frac{da}{ds}\right|,\cdots,
  \sup_{s\in\mathbb{R}}\left|\frac{d^3a}{ds^3}\right|
  \;\lesssim\;\exp(-2S_0)\,.
  \end{equation}
  Thus we construct a universal smooth $\hat v_2(s)$ with
  \begin{equation}\label{4.6}
  \left[\frac{d}{ds}(1-a)\left(\frac{d}{ds}+1\right)-2\right]\,\hat v_2
  \;=\;1\quad\mbox{for all}\;s
  \end{equation}
  and
  \begin{equation}\label{4.7}
  \sup_{s\in\mathbb{R}}|\hat v_2|,
  \sup_{s\in\mathbb{R}}\left|\frac{d\hat v_2}{ds}\right|,\cdots,
  \sup_{s\in\mathbb{R}}\left|\frac{d^4\hat v_2}{ds^4}\right|<\infty\,.
  \end{equation}
  We finally reformulate (\ref{4.6}) as a fixed point problem. Note that
  since $\frac{d}{ds}(\frac{d}{ds}+1)-2=(\frac{d}{ds}-1)(\frac{d}{ds}+2)$,
  the bounded solution of $\left[\frac{d}{ds}(\frac{d}{ds}+1)-2\right]\,\hat v
  \;=\;\hat f$ for some bounded continuous $\hat f$ is given by
  \begin{eqnarray}\nonumber
  \hat v(s)&=&-\int_{-\infty}^s\exp(2(s'-s))\int_{s'}^{\infty}
  \exp(s'-s'')\hat f(s'')\,ds''\,ds'\nonumber\\
  &=&
  -\frac{1}{3}\int_{-\infty}^\infty\exp(3\min\{s,s''\}-2s-s'')\hat f(s'')\,ds''
  \nonumber\\
  &=&:\;(T \hat f)(s)\,,
  \end{eqnarray}
  defining an operator $T$. From its above representation with
  the bounded and Lipschitz-continuous kernel $\exp(3\min\{s,s''\}-2s-s'')$ we read off
  that $T$ is a bounded operator from $C^0$ (the space of bounded continuous
  functions endowed with the sup norm) into $C^1$ and by the solution property
  of $T$ thus also into $C^2$. Note that (\ref{4.6}) can be reformulated as 
  \begin{eqnarray}\nonumber
  \lefteqn{\left[\frac{d}{ds}\left(\frac{d}{ds}+1\right)-2\right]\,\hat v_2}\nonumber\\
  &=&1+\frac{d}{ds}a\left(\frac{d}{ds}+1\right)\hat v_2\nonumber\\
  &=&1+\left[\left(\frac{d^2}{ds^2}+\frac{d}{ds}\right)a
           -\frac{d}{ds}\frac{da}{ds}\right]\hat v_2\nonumber\\
  &=&1+\left[\left(\frac{d}{ds}\left(\frac{d}{ds}+1\right)-2\right)a
           -\frac{d}{ds}\frac{da}{ds}
           +2 a\right]\hat v_2\,.
  \end{eqnarray}
  An application of the translation-invariant operator $T$
  (formally) yields
  \begin{equation}\label{4.9}
  \hat v_2\;=\;T\,1+\left(a-\frac{d}{ds}\,T\,\frac{da}{ds}+2 T\,a\right)\hat v_2\,.
  \end{equation}
  We view this equation as a fixed-point equation for $\hat v_2$ in the Banach 
  space $C^0$. As mentioned above,
  $T$ and even the composition $\frac{d}{ds}\,T$ are bounded operators (in $C^0$).
  In view of (\ref{4.8}), the multiplication with $a$ and with $\frac{da}{ds}$
  are operators with $C^0$-operator norm estimated by $C\exp(-2S_0)$.
  Hence for sufficiently large $S_0$, the operator 
  $a-\frac{d}{ds}\,T\,\frac{da}{ds}+2 T\,a$ has norm strictly less than
  one. Thus the contraction mapping theorem ensures the existence of
  a solution of (\ref{4.9}), that is, a $C^2$-solution $\hat v_2$ of (\ref{4.6})
  with $\sup_{s\in\mathbb{R}}|\hat v_2|$, 
  $\sup_{s\in\mathbb{R}}|\frac{d\hat v_2}{ds}|$,
  $\sup_{s\in\mathbb{R}}|\frac{d^2\hat v_2}{ds^2}|<\infty$.
  Finally, we obtain the rest of (\ref{4.7}) from (\ref{4.8})
  by a booth-strap argument.

\bigskip

\subsection{Proof of Lemma \ref{Lemma4}}

 Here we give the argument for (\ref{CLAIM-1}).
  We note that by definition (\ref{convolution}) of the convolution $\hat\xi_0$,
 the change of variables $s=\ln \hat z$ and $s'=-\ln k$ already used in (\ref{k-change}) \& (\ref{her.2}),
  and by definition (\ref{new-variables}) of $\hat\xi$ we have

\begin{eqnarray}\nonumber
\int_{-\infty}^{\ln H}\hat\xi_0\,ds'
&\stackrel{(\ref{convolution})}{=}&
\int_{-\infty}^{\ln H}\int_{-\infty}^\infty\hat\xi(s+s')\phi_0(s)\,ds\,ds'\nonumber\\
&\stackrel{(\ref{k-change})\&(\ref{her.2})}{=}&
\int_{\frac{1}{H}}^\infty\int_0^1\hat\xi\left(\frac{\hat z}{k}\right)\,
\phi_0(\hat z)\,\frac{d\hat z}{\hat z}\,\frac{dk}{k}\nonumber\\
&\stackrel{(\ref{new-variables})}{=}&
\int_0^1\int_{\frac{1}{H}}^\infty\xi\left(\frac{\hat z}{k}\right)\,\frac{dk}{k^2}\,\phi_0(\hat z)\,
\,d\hat z\nonumber\\
&=&
\int_0^1\int_0^{H\hat z}\xi(z)\,\frac{dz}{\hat z}\,\phi_0(\hat z)\,d\hat z\nonumber\\
&=&
\int_0^H\xi(z)\,\int_{\frac{z}{H}}^1\frac{1}{\hat z}\phi_0(\hat z)\,d\hat z\,dz\,.\nonumber
\end{eqnarray}
In view of this identity and the up-down symmetry (i.\ e.\ the symmetry of the
problem under $z\leadsto H-z$), (\ref{CLAIM-1}) will follow if we show that
$\int_0^H\xi\,dz=-1$ implies
 \begin{equation}\label{CLAIM-2}
 \int_{0}^{H}\xi(z)\left(\int_{\frac{z}{H}}^{1}\frac{1}{\hat z}\phi_0(\hat z)\,d\hat z+\int_{1-\frac{z}{H}}^{1}\frac{1}{\hat z}\phi_0(\hat z)\,d\hat z\right)dz\lesssim -1 \,.
 \end{equation}
 %
 With the normalization (\ref{convolution-kernel}) and our assumption $\int_0^H\xi\,dz=-1$, (\ref{CLAIM-2}) will 
 follow once we show 
 \begin{equation}\label{CLAIM-3}
 \int_{0}^{H}\xi(z)\left(1-\int_{\frac{z}{H}}^1\frac{1}{\hat z}\phi_0(\hat z)\,d\hat z-\int_{1-\frac{z}{H}}\frac{1}{\hat z}\phi_0(\hat z)\,d\hat z\right)dz\gtrsim -\frac{(\ln H)^{\frac{1}{45}}}{H^{\frac 23}}\,,
 \end{equation}
 for which we will use the second assumption on $\int_0^H\xi^2\, dz$: Claim (\ref{CLAIM-3}) clearly implies (\ref{CLAIM-2}) in the regime of $H\gg1$.
%
%
%
%
Let us reformulate (\ref{CLAIM-3}) as 
\begin{equation}\label{CLAIM-4}
\int_{0}^{H}\xi(z)\rho(z)\,dz\gtrsim -\frac{(\ln H)^{\frac{1}{45}}}{H^{\frac 23}}\,,
\end{equation}
where we introduced 
\begin{equation}\label{rho-rescaled}
 \rho(z):=\rho_{0}\left(\frac{z}{H}\right)\; \mbox{ with }\; \rho_0\left(\hat z\right):=
 1-\int_{\hat z}^{1}\frac{1}{\hat z'}\phi_0(\hat z')\,d\hat z'-\int_{1-\hat z}^{1}\frac{1}{\hat z'}\phi_0(\hat z')\,d\hat z'\,.
\end{equation}
%
The symmetry (\ref{convolution-kernel}) of $\phi_0(\hat z)\geq 0$ implies that
\begin{equation*}
\frac{d\rho_0}{d\hat z}(\hat z)\;=\;\left(\frac{1}{\hat z}-\frac{1}{1-\hat z}\right)\phi_0(\hat z)\;
\begin{cases}
\ge 0 &\mbox{ for } \hat z\le\frac{1}{2}\,,\\
\le 0 &\mbox{ for } \hat z\ge\frac{1}{2}\,,
\end{cases}
\end{equation*}
so that using the normalization (\ref{convolution-kernel}) we have
\begin{equation}\label{rho0-positivity}
\rho_0\ge 0\,,
\end{equation}
and
\begin{equation}\label{upper-bound-for-rho0}
\rho_0\le 1 \,.
\end{equation}
  Hence (\ref{CLAIM-4}) is yet another way of expressing approximate non-negativity of $\xi$, this time in and up-down symmetric way in the bulk.
  
  The strategy to establish (\ref{CLAIM-3}) is now to construct an even (but not necessary non-negative) mollification kernel $\phi(z)$ 
   of length-scale $\ell\ll H$ such that 
  \begin{eqnarray}
  &(\xi\ast\phi)(z)\gtrsim -\frac{1}{\ell^4}  &\mbox{ for }  z\in (\ell,  H-\ell),\label{from-stab}\\
  &\int_{0}^{H}(\phi\ast \rho-\rho)^2\, dz\lesssim \frac{\ell^4}{H^3} &\mbox{ for } \ell\ll H.\label{scaleH}
  \end{eqnarray}
  We first argue how (\ref{from-stab}) and (\ref{scaleH}) imply (\ref{CLAIM-3}).
  Indeed by the evenness of $\phi$ we have the representation
  $$\int_{-\infty}^{\infty}\xi\rho \,dz=\int_{-\infty}^{\infty}(\xi\ast\phi)\rho \,dz-\int_{-\infty}^{\infty}\xi(\rho\ast\phi-\rho)\,dz\,,$$
  from which, since $\rho\geq 0$ (cf.\ (\ref{rho0-positivity})), we get
  $$\int_{-\infty}^{\infty}\xi \rho \,dz\geq \inf_{z\in \rm{supp}\rho}(\xi\ast\phi)(z)\int_{-\infty}^{\infty}\rho \,dz-\left(\int_{-\infty}^{\infty}\xi^2 dz\int_{-\infty}^{\infty}(\rho\ast\phi-\rho)^2\,dz\right)^{\frac{1}{2}}\,.$$
  We note that since $\phi_0(\hat z)$ is supported in $[\frac{1}{4},\frac{3}{4}]$, cf.\ Lemma \ref{Lemma1}, $\rho_0(\hat z)$
  is supported in the same interval. Hence $\rho$ is supported in $[\frac{1}{4}H,\frac{3}{4}H]$. Hence 
  we may apply (\ref{from-stab}) as soon as $\ell\leq \frac{H}{4}$.
  Using (\ref{from-stab}) and (\ref{scaleH}) together with our assumption that $\int\xi^2\, dz\lesssim (\ln H)^{\frac{1}{15}}$ and $\int_0^H\rho \,dz\leq H$ (from (\ref{upper-bound-for-rho0}))
  we obtain the estimate 
  $$\int_{-\infty}^{\infty}\xi\rho \,dz\gtrsim -\frac{H}{\ell^4}-\left((\ln H)^{\frac{1}{15}}\frac{\ell^4}{H^3}\right)^{\frac{1}{2}}\,.$$
  The balancing choice of $\ell=\left(\frac{H^5}{(\ln H)^{\frac{1}{15}}}\right)^{\frac{1}{12}}$ turns this estimate into (\ref{CLAIM-3}).

  We now turn to the construction of the mollification kernel $\phi$. We select a (nonvanishing) smooth 
  and even $w_0(\hat z)$, compactly supported in $\hat z\in[-1,1]$, and consider the corresponding multiplier
  $$\phi_0=w_0\left(-\frac{d^2}{d\hat z^2}+1\right)^2 w_0\,.$$
  Notice that 
  $\int_{0}^{H}\phi_0\,d\hat z=\int_{-\infty}^{\infty}\left( \left(\frac{d^2w_0}{d\hat z^2}\right)^2+\left(\frac{dw_0}{d\hat z}\right)^2+w_0^2\right) \,d\hat z>0\,,$
  so that by changing $w_0$ by a multiplicative constant we may achieve
  $$\int_{-\infty}^{\infty}\phi_0 \,d\hat z=1.$$
  We change variables according to $z=\ell \hat z$ and rescale the mask $\phi_0$ by $\ell$ so as to preserve its integral 
  \begin{equation}\label{mask-rescaled}
 \ell\phi(\ell\hat z)=\phi_0(\hat z)\,,
  \end{equation}
  and note that also $\phi$ is a multiplier in the sense of
  \begin{equation}\label{phi-rescaled}
  \phi=w\left(-\frac{d^2}{d\hat z^2}+\frac{1}{\ell^2}\right)^2w\,,
  \end{equation}
  provided $w$ is the following rescaling of $w_0$:
  \begin{equation}\label{w-rescaled}
   \frac{1}{\ell^{\frac{3}{2}}}w(\ell\hat z)= w_0(\hat z)\,.
  \end{equation}
  For any translation $z'\in(\ell,H-\ell)$, the translated test function
  $z\mapsto w(z-z')$ is compactly supported in $z\in(0,H)$ and we may thus 
  apply the stability condition (\ref{OSCinFourier-reduct}) with $k=\frac{1}{\ell}$. 
  Because of (\ref{phi-rescaled}), this yields (\ref{from-stab}):
  \begin{eqnarray}\nonumber
  \lefteqn{\int_0^H\xi(z)\,\phi(z-z')\,dz}\nonumber\\
  &\ge&
  -\int_0^H\left[\ell^2\left(\frac{d}{dz}\left(-\frac{d^2}{dz^2}+\frac{1}{\ell^2}\right)^2w\right)^2
  +\left(\left(-\frac{d^2}{dz^2}+\frac{1}{\ell^2}\right)^2w\right)^2\right](z-z')\,dz\nonumber\\
  &\stackrel{(\ref{w-rescaled})}{=}&
  -\frac{1}{\ell^4}\int_{-1}^1\left[\left(\frac{d}{d\hat z}\left(-\frac{d^2}{d\hat z^2}+1\right)^2w_0\right)^2
  +\left(\left(-\frac{d^2}{d\hat z^2}+1\right)^2w_0\right)^2\right]\,d\hat z\;\sim\;-\frac{1}{\ell^4}\,.\nonumber
  \end{eqnarray}	
  We finally turn to (\ref{scaleH}). From the representation	
  \begin{eqnarray*}
  (\rho*\phi-\rho)(z')&=&
  \int_{-\infty}^\infty(\rho(z'-z)-\rho(z'))\,\phi(z)\,dz\\
  &\stackrel{\mbox{$\phi$ is even}}{=}&
  \frac{1}{2}\int_{-\infty}^\infty(\rho(z'+z)+\rho(z'-z)-2\rho(z'))\,\phi(z)\,dz\,,
  \end{eqnarray*}
  we obtain the inequality
  \begin{eqnarray*}
  |(\rho*\phi-\rho)(z')|&\le&
  \frac{1}{2}\sup\left|\frac{d^2\rho}{dz^2}\right|\int_{-\infty}^\infty z^2|\phi(z)|\,dz\\
  &\stackrel{(\ref{rho-rescaled}),(\ref{mask-rescaled})}{=}&\frac{1}{H^2}\sup\left|\frac{d^2\rho_0}{d\hat z^2}\right|
  \ell^2\int_{-\infty}^\infty\hat z^2|\phi_0(\hat z)|\,d\hat z\,,
  \end{eqnarray*}
  which yields (\ref{scaleH}) after integration in $z'\in[0,H]$.


\section{Appendix}
\subsection{Appendix for Section \ref{heuristics}}\label{Appendix1}
  Here, we argue how to derive (\ref{new-star}). Recall the change of variables (\ref{nonlinear-change})
  for $s'\longleftrightarrow \hat s$ with $s''$ as a fixed parameter.  
 If $p$, $\tilde p$
denote generic polynomials of degree $n$, we have
\begin{eqnarray}\nonumber
\frac{1}{(s')^m}\frac{d^n}{d\hat s^n}
&\stackrel{(\ref{nonlinear-change})}{=}&\frac{1}{(s'')^n}\frac{1}{(s')^{m-n}}(1+\hat s)^n\frac{d^n}{d\hat s^n}\nonumber\\
&=&\frac{1}{(s'')^n}\frac{1}{(s')^{m-n}}p((1+\hat s)\frac{d}{d\hat s})\nonumber\\
&\stackrel{(\ref{nonlinear-change})}{=}&\frac{1}{(s'')^n}\frac{1}{(s')^{m-n}}p(-s'\frac{d}{ds'})\nonumber\\
&=&\frac{1}{(s'')^n}\tilde p(s'\frac{d}{ds'})\frac{1}{(s')^{m-n}}\nonumber\\
&=&\frac{1}{(s'')^n}\sum_{k=0}^na_n\frac{d^k}{ds'^k}\frac{1}{(s')^{m-n-k}}\nonumber\\
&\stackrel{(\ref{nonlinear-change})}{=}&\sum_{k=0}^na_n\frac{1}{(s'')^{m-k}}\frac{d^k}{ds'^k}(1+\hat s)^{m-n-k}.\nonumber
\end{eqnarray}
This shows that the desired relations exist in principle, it remains to determine the coefficients $a_0,\cdots,a_n$.
We start with the case $m=n+1$ (which yields the shortest formula). To this purpose, we
again use $(1+\hat s)\frac{d}{d\hat s}=-s'\frac{d}{ds'}$, which we rewrite as
$\frac{d}{d\hat s}(1+\hat s)=-(s')^2\frac{d}{ds'}\frac{1}{s'}$.
The latter yields
\begin{equation*}
\left(\frac{d}{d\hat s}(1+\hat s)\right)^n\;=\;(-1)^n
s'\left(s'\frac{d}{ds'}\right)^n\frac{1}{s'}\quad\mbox{for every}\;n\in\mathbb{N}\,,
\end{equation*}
  which implies inductively
  \begin{equation*}
  \frac{d^n}{d\hat s^n}(1+\hat s)^n\;=\;(-1)^n
  (s')^{1+n}\frac{d^n}{ds'^n}\frac{1}{s'}\quad\mbox{for every}\;n\in\mathbb{N}\,,
  \end{equation*}
 and which we rewrite as (using again $s''=s'(1+\hat s)$)
  \begin{equation}\label{tranformation-formula}
  \frac{1}{(s')^{n+1}}\frac{d^n}{d\hat s^n}
  \;=\;(-1)^n\frac{d^n}{ds'^n}\frac{1}{s'}\frac{1}{(1+\hat s)^n}
  \;=\;(-1)^n\frac{1}{s''}\frac{d^n}{ds'^n}\frac{1}{(1+\hat s)^{n-1}}.
  \end{equation}
  In view of the first line on the r.\ h.\ s.\ of (\ref{phi-in-s-prime-2}), we need the
  latter transformation formula for $n=1,2,3,4$. 
  In view of the second line, we also need:
  \begin{eqnarray}\label{A1}
  \frac{1}{(s')^{4}}\frac{d}{d\hat s}&\stackrel{(\ref{tranformation-formula})}{=}&-\frac{1}{s''}\frac{1}{(s')^2}\frac{d}{ds'}\nonumber\\
  &=&-\frac{1}{s''}(\frac{d}{ds'}+\frac{2}{s'})\frac{1}{(s')^2}\nonumber\\
  &=&-\frac{1}{(s'')^3}\frac{d}{ds'}(1+\hat s)^2-\frac{2}{(s'')^4}(1+\hat s)^3\,,
  \end{eqnarray}
  \begin{eqnarray}\label{A2}
  \frac{1}{(s')^{5}}\frac{d^2}{d\hat s^2}
  &\stackrel{(\ref{tranformation-formula})}{=}&\frac{1}{s''}\frac{1}{(s')^2}\frac{d^2}{ds'^2}\frac{1}{1+\hat s}\nonumber\\
  &=&\frac{1}{s''}\left(\frac{d^2}{ds'^2}+4\frac{d}{ds'}\frac{1}{s'}+6\frac{1}{(s')^2}\right)\frac{1}{(s')^2}\frac{1}{1+\hat s}\nonumber\\
  &=&\frac{1}{(s'')^3}\frac{d^2}{ds'^2}(1+\hat s)+\frac{4}{(s'')^4}\frac{d  }{ds'  }(1+\hat s)^2+\frac{6}{(s'')^5}(1+\hat s)^3\,.
  \end{eqnarray}

  \subsection{Appendix for Subsection \ref{Proof-Lemma3}}
  In this subsection we derive the formulas (\ref{box}) and (\ref{boxbox}).
  The main step is to establish
  \begin{eqnarray}
  &&\hat z^2(\partial_{\hat z}^4-2\partial_{\hat z}^2+1)
  \hat z\sinh\hat z\label{1.1}\\
  &&=\Big(\hat z^{-1}\sinh\hat z\,(\partial_s-2)(\partial_s-1)
  +4\cosh\hat z(\partial_s-1)
  +4\hat z\sinh\hat z\Big)
  (\partial_s+1)\partial_s\,.\nonumber
  \end{eqnarray}
%
%
  Let us give a motivation for formula (\ref{1.1}): The factor $(\partial_s+1)\partial_s$
  has to be there since $\hat z^{-1}=e^{-s}$ and $1$ are in the kernel
  of $(\partial_{\hat z}^4-2\partial_{\hat z}^2+1)\hat z\sinh\hat z$,
  which in turn follows from the fact that $\sinh\hat z$ and $\hat z\sinh\hat z$ are in the kernel of
  $\partial_{\hat z}^4-2\partial_{\hat z}^2+1$. Note that
  for $\hat z\ll 1$,
  \begin{equation*}
  \hat z^{-1}\sinh\hat z=1+O(\hat z^2),\quad
  \cosh\hat z=1+O(\hat z^2),\quad
    \hat z\sinh\hat z=O(\hat z^2)\,,
  \end{equation*}
  so that for $\hat z\ll 1$, (\ref{1.1}) collapses to the identity already used in (\ref{poly})
  \begin{equation}\label{1.2}
  \hat z^2\partial_{\hat z}^4\hat z^2
  \;=\;
  (\partial_s+2)(\partial_s+1)\partial_s(\partial_s-1)\,.
  \end{equation}
  This identity is easily seen to be true because both differential
  operators are of fourth order and
  are homogeneous of degree zero in $\hat z$,
  because the four
  functions $\hat z^{-2}=e^{-2s}$, $\hat z^{-1}=e^{-s}$,
  $1$, and $\hat z=e^s$ are in the kernel of both differential operators,
  and because on $\hat z^2=e^{2s}$, both operators give $4!\hat z^2=
  4!e^{2s}$. 
%

%

  Let us give the argument for (\ref{1.1}). Because of the transformation
  properties under $\hat z\leadsto -\hat z$, it suffices to show
  \begin{eqnarray*}
  &&\hat z^2(\partial_{\hat z}^4-2\partial_{\hat z}^2+1)\hat z \exp(\hat z)\\
  &&=\left[\hat z^{-1}\exp(\hat z)\,(\partial_s-2)(\partial_s-1)
  +4\exp(\hat z)(\partial_s-1)
  +4\hat z\exp(\hat z)\right](\partial_s+1)\partial_s\,,\nonumber
  \end{eqnarray*}
  which we rearrange as
  \begin{eqnarray}
  \lefteqn{\hat z^3(\partial_{\hat z}^4
  -2\partial_{\hat z}^2+1)\hat z \exp(\hat z)}\nonumber\\
  &=&\exp(\hat z)\left[(\partial_s-2)(\partial_s-1)
  +4\hat z(\partial_s-1)
  +4\hat z^2\right](\partial_s+1)\partial_s\,.\label{1.4}
  \end{eqnarray}
  We note that because of $\partial_{\hat z}\exp(\hat z)=
  \exp(\hat z)(\partial_{\hat z}+1)$, we have
  \begin{eqnarray}\nonumber
  (\partial_{\hat z}^4-2\partial_{\hat z}^2+1)\exp(\hat z)
  &=&
  \exp(\hat z)\left[(\partial_{\hat z}+1)^4-2(\partial_{\hat z}+1)^2+1\right]\nonumber\\
  &=&
  \exp(\hat z)(\partial_{\hat z}^4+4\partial_{\hat z}^3+4\partial_{\hat z}^2)\,,\nonumber
  \end{eqnarray}
  so that	
  \begin{equation}\nonumber
  \hat z^2(\partial_{\hat z}^4-2\partial_{\hat z}^2+1)\hat z\exp(\hat z)
  =
  \exp(\hat z)\left[\hat z^3\partial_{\hat z}^4\hat z
  +4\hat z(\hat z^2\partial_{\hat z}^3\hat z)
  +4\hat z^2(\hat z\partial_{\hat z}^2\hat z)\right]\,.
  \nonumber
  \end{equation}
  Now (\ref{1.4}) follows by inserting the formulas
  \begin{eqnarray}
  \hat z\partial_{\hat z}^2\hat z&=&(\partial_s+1)\partial_s,\nonumber\\
  \hat z^2\partial_{\hat z}^3\hat z&=&(\partial_s+1)\partial_s(\partial_s-1),
  \nonumber\\
  \hat z^3\partial_{\hat z}^4\hat z&=&(\partial_s+1)\partial_s
  (\partial_s-1)(\partial_s-2)\,.\label{1.3}
  \end{eqnarray}
  These formulas can easily seen to be true; let us address (\ref{1.3}):
  Both sides are differential operators of order 4 that are homogeneous
  of degree 0 in $\hat z$; the kernel of both operators is spanned by
  the four functions $\hat z^{-1}=e^{-s}$, $1$, $\hat z=e^{s}$, and
  $\hat z^2=e^{2s}$; On $\hat z^3=e^{3s}$, both operators yield
  $4!\hat z^3=4!e^{3s}$. 

  Formulas (\ref{box}) and (\ref{boxbox}) easily follow from (\ref{1.1}).
  Formula (\ref{box}) is an immediate consequence of (\ref{1.1}). 
  Formula (\ref{boxbox}) follows from (\ref{box}) using the identities
  $\partial_{\hat z}=\hat z^{-1}\partial_s$ and
  \begin{eqnarray}
  \partial_{\hat z}\left(\hat z^{-3}\sinh\hat z\right)&=&
  \hat z^{-3}(\cosh\hat z-3\hat z^{-1}\sinh\hat z),\nonumber\\
  \partial_{\hat z}\left(4\hat z^{-2}\cosh\hat z\right)&=&
  \hat z^{-3}(4\hat z\sinh\hat z-8\cosh\hat z),\nonumber\\
  \partial_{\hat z}\left(4\hat z^{-1}\sinh\hat z\right)&=&
  \hat z^{-3}(4\hat z^2\cosh\hat z-4\hat z\sinh\hat z)\,,\nonumber
  \end{eqnarray}
  which lead as desired to
  \begin{eqnarray}
  &&\partial_{\hat z}(\partial_{\hat z}^4-2\partial_{\hat z}^2+1)
  \hat z\sinh\hat z\nonumber\\
  &&=\hat z^{-3}\left[(
  \hat z^{-1}\sinh\hat z\,\left((\partial_s-2)(\partial_s-1)\partial_s
  -3(\partial_s-2)(\partial_s-1)\right)\right.
  \nonumber\\
  &&\quad+\cosh\hat z\,\left(4(\partial_s-1)\partial_s
  +(\partial_s-2)(\partial_s-1)-8(\partial_s-1)\right)
  \nonumber\\
  &&\quad+\hat z\sinh\hat z\,\left(4\partial_s+4(\partial_s-1)-4\right)
  \left.\quad+\hat z^2\cosh\hat z\,4\right]
  \times(\partial_s+1)\partial_s\nonumber\\
  &&=\hat z^{-3}\Big(
   \hat z^{-1}\sinh\hat z\,(\partial_s-3)(\partial_s-2)(\partial_s-1)
   +5\cosh\hat z(\partial_s-2)(\partial_s-1)\nonumber\\
   &&\;\;+8\hat z\sinh\hat z(\partial_s-1)
   +4\hat z^2\cosh\hat z\Big)
   (\partial_s+1)\partial_s\,.
  \end{eqnarray}

  \subsection{Notations}\label{notations}

  \underline{The spatial vector}:
  $$x=(y,z)\in [0,L)^{d-1}\times [0,H] \,,$$
  where $H$ denotes the height of the container and $L$ is the lateral horizontal cell-size.
  
  \noindent
  \underline{Vertical velocity component}:
  \[w:=u\cdot e_z \qquad \mbox{ where } \qquad u=u(y,z,t)\,.\]

  \noindent
  \underline{Background profile}:
  \[\tau:[0,H]\rightarrow \R \quad\mbox{ such that } \quad \tau(0)=1 \mbox{ and } \tau(H)=1\,,\]
  \[\tau=\tau(z)\,,\qquad  \qquad \qquad \xi:=\frac{d\tau}{dz}\,.\]

  \noindent
  \underline{Long-time and horizontal average}:
  $$\langle f \rangle:=\limsup_{t_0\uparrow\infty}\frac{1}{t_0}\int_0^{t_0}\frac{1}{L^{d-1}}\int_{[0,L)^{d-1}}f(t,y)dydt.$$
  
  \noindent
  \underline{Gradient}:
  \[\nabla f=\left(\begin{array}{c}
                  \nabla_{y} \\ \partial_z
                 \end{array} 
  \right) f\,.\]
 
 \noindent
  \underline{Laplacian}:
  \[\Delta f=\Delta_{y} f+\partial_z^2 f\,.\]

  \noindent
  \underline{Horizontal Fourier transform:}
  \[\F f(k,z)=\frac{1}{L^{d-1}}\int_{[0,L)^{d-1}} e^{-ik\cdot y}f(y,z)dy\,,\]
  where $k\in \frac{2\pi}{L}\Z^{d-1}$ is the dual variable of $y$.
  
  \noindent
  \underline{Real part of an imaginary number }:
  \rm{Re} stands for the real part of a complex number. 
  
  \noindent
  \underline{Complex conjugate }:
  $\overline{\F w}$ and $\overline{\F \theta}$ are the complex conjugates of the (complex valued) functions $\F w$ and $\F \theta$. 

  \noindent
  \underline{Universal and specific constants:}
  We call \textit{universal constant} a constant $C$ such that  $0<C<\infty$ and it only depends on $d$  but not on $H$, on $L$ and on the initial data.
  Throughout the paper $A\lesssim B$ means $A\leq CB$ with $C$ a universal constant. Likewise
  a condition $A\ll B$ means that there exists a possibly large universal constant $C$ such that $A\leq \frac{1}{C}B$. 
  We indicate specific constants with $C_0, C_1, C_2,\cdots$.


  \section*{Acknowledgement}
  The authors thank Charlie Doering for the many insightful discussion and for pointing out the connection of our result with the 
  one of Ierley, Kerswell and Plasting.

  C.N. was partially supported by the IMPRS of MPI MIS (Leipzig), the Max-Planck Institute for Mathematics in the Sciences in Leipzig,
  and by the University of Basel.
 
  \bibliographystyle{unsrt}
  \bibliography{BiblioLB}

  \end{document}